\documentclass[leqno,12pt]{article}
\usepackage{amsmath,amssymb,rotating,a4wide,color}
\usepackage{wasysym}
\usepackage{hyperref}
\usepackage{enumerate}
\usepackage{mathtools}
\usepackage{float}
\usepackage{colortbl}
\usepackage{makecell}
\usepackage[font=footnotesize,labelfont=bf,margin=2cm]{caption}
\usepackage[isbn=false,maxnames=15]{biblatex}
\DefineBibliographyStrings{english}{%
	pages = {}, page = {}
}

\usepackage{listings}
\usepackage{color}

\definecolor{mygreen}{rgb}{0,0.6,0}
\definecolor{mygray}{rgb}{0.5,0.5,0.5}
\definecolor{mymauve}{rgb}{0.58,0,0.82}
\definecolor{backcolour}{rgb}{0.95,0.95,0.92}

\usepackage{array}
\newcolumntype{M}{>{\centering\arraybackslash}m{\dimexpr.25\linewidth-2\tabcolsep}}

\usepackage{tikz}
\usetikzlibrary{circuits.logic.US,circuits.logic.IEC,fit}
\usetikzlibrary{cd}
\usetikzlibrary{backgrounds}
\usetikzlibrary{math}
\usepackage{subcaption}
\usepackage{lscape}

\usepackage[all,cmtip]{xy}

\usepackage{verbatim}

\newdir^{ (}{{}*!/-3pt/\dir^{(}}
\newdir_{ (}{{}*!/-3pt/\dir_{(}}

\entrymodifiers={+!!<0pt,\fontdimen22\textfont2>}

\definecolor{cbblue}{HTML}{8da0cb}
\definecolor{cborange}{HTML}{fc8d62}
\definecolor{cbgreen}{HTML}{66c2a5}

\def\geqone{g\kern1pt{=}\kern1pt1}
\def\geqtwo{g\kern1pt{=}\kern1pt2}

\def\FigTwoAllA{\begin{tikzpicture}[scale=0.7]
\draw [thick,black] plot  coordinates {(-4/8,-1/8) (24/8,6/8) };
\begin{scope}[shift={(0,0)}]
\draw [thick,black] plot  coordinates {(-1/32,1/8) (1/32,-1/8) };
\end{scope}
\begin{scope}[shift={(1/2,1/8)}]
\draw [thick,black] plot  coordinates {(-1/32,1/8) (1/32,-1/8) };
\end{scope}
\begin{scope}[shift={(2/2,2/8)}]
\draw [thick,black] plot  coordinates {(-1/32,1/8) (1/32,-1/8) };
\end{scope}
\begin{scope}[shift={(3/2,3/8)}]
\draw [thick,black] plot  coordinates {(-1/32,1/8) (1/32,-1/8) };
\end{scope}
\begin{scope}[shift={(4/2,4/8)}]
\draw [thick,black] plot  coordinates {(-1/32,1/8) (1/32,-1/8) };
\end{scope}
\begin{scope}[shift={(5/2,5/8)}]
\draw [thick,black] plot  coordinates {(-1/32,1/8) (1/32,-1/8) };
\end{scope}
\end{tikzpicture}}

\def\FigTwoAllB{\begin{tikzpicture}[scale=0.7]
\draw [thick,black] plot  coordinates {(-48/24,-32/24) (3/8,2/8) };
\draw [thick,black] plot  coordinates {(6/8,-9/8) (-4/24,6/24) };
\fill [black] (0,0) circle (0.08);
\begin{scope}[shift={(-3/8,-2/8)}]
\draw [thick,black] plot  coordinates {(-2/32,3/32) (2/32,-3/32) };
\end{scope}
\begin{scope}[shift={(-6/8,-4/8)}]
\draw [thick,black] plot  coordinates {(-2/32,3/32) (2/32,-3/32) };
\end{scope}
\begin{scope}[shift={(-9/8,-6/8)}]
\draw [thick,black] plot  coordinates {(-2/32,3/32) (2/32,-3/32) };
\end{scope}
\begin{scope}[shift={(-12/8,-8/8)}]
\draw [thick,black] plot  coordinates {(-2/32,3/32) (2/32,-3/32) };
\end{scope}
\begin{scope}[shift={(2/8,-3/8)}]
\draw [thick,black] plot  coordinates {(3/32,2/32) (-3/32,-2/32) };
\end{scope}
\begin{scope}[shift={(4/8,-6/8)}]
\draw [thick,black] plot  coordinates {(3/32,2/32) (-3/32,-2/32) };
\end{scope}
\end{tikzpicture}}

\def\FigTwoAllE{\begin{tikzpicture}[scale=0.7]
\draw (7/4,0/4) node {};
\draw [thick,black] plot  coordinates {(-11/8,-11/8) (1/4,1/4) };
\draw [thick,black] plot  coordinates {(11/8,-11/8) (-1/4,1/4) };
\begin{scope}[shift={(0,0)}]
\fill [white] (0,0) circle (0.1);
\draw [thick] (0,0) circle (0.1);
\end{scope}
\begin{scope}[shift={(1/3,-1/3)}]
\draw [thick,black] plot  coordinates {(-1/12,-1/12) (1/12,1/12) };
\end{scope}
\begin{scope}[shift={(2/3,-2/3)}]
\draw [thick,black] plot  coordinates {(-1/12,-1/12) (1/12,1/12) };
\end{scope}
\begin{scope}[shift={(3/3,-3/3)}]
\draw [thick,black] plot  coordinates {(-1/12,-1/12) (1/12,1/12) };
\end{scope}
\begin{scope}[shift={(-1/3,-1/3)}]
\draw [thick,black] plot  coordinates {(-1/12,1/12) (1/12,-1/12) };
\end{scope}
\begin{scope}[shift={(-2/3,-2/3)}]
\draw [thick,black] plot  coordinates {(-1/12,1/12) (1/12,-1/12) };
\end{scope}
\begin{scope}[shift={(-3/3,-3/3)}]
\draw [thick,black] plot  coordinates {(-1/12,1/12) (1/12,-1/12) };
\end{scope}
\end{tikzpicture}}

\def\FigTwoAllC{\begin{tikzpicture}[scale=0.7]
\draw [thick,black] plot  coordinates {(-10/8,-10/8) (1/4,1/4) };
\draw [thick,black] plot  coordinates {(10/8,-10/8) (-1/4,1/4) };
\draw [thick,black] plot  coordinates {(6/8,-10/8) (18/8,2/8) };
\fill [black] (0,0) circle (0.08);
\fill [black] (4/4,-4/4) circle (0.08);
\begin{scope}[shift={(-1/2,-1/2)}]
\draw [thick,black] plot  coordinates {(-1/12,1/12) (1/12,-1/12) };
\end{scope}
\begin{scope}[shift={(-5/6,-5/6)}]
\draw [thick,black] plot  coordinates {(-1/12,1/12) (1/12,-1/12) };
\end{scope}
\begin{scope}[shift={(1/3,-1/3)}]
\draw [thick,black] plot  coordinates {(-1/12,-1/12) (1/12,1/12) };
\end{scope}
\begin{scope}[shift={(2/3,-2/3)}]
\draw [thick,black] plot  coordinates {(-1/12,-1/12) (1/12,1/12) };
\end{scope}
\begin{scope}[shift={(2-1/6,-1/6)}]
\draw [thick,black] plot  coordinates {(-1/12,1/12) (1/12,-1/12) };
\end{scope}
\begin{scope}[shift={(2-3/6,-3/6)}]
\draw [thick,black] plot  coordinates {(-1/12,1/12) (1/12,-1/12) };
\end{scope}
\end{tikzpicture}}

\def\FigTwoAllF{\begin{tikzpicture}[scale=0.7]
\draw [thick,black] plot  coordinates {(-10/8,-10/8) (1/4,1/4) };
\draw [thick,black] plot  coordinates {(10/8,-10/8) (-1/4,1/4) };
\draw [thick,black] plot  coordinates {(6/8,-10/8) (18/8,2/8) };
\begin{scope}[shift={(0,0)}]
\fill [white] (0,0) circle (0.1);
\draw [thick] (0,0) circle (0.1);
\end{scope}
\fill [black] (4/4,-4/4) circle (0.08);
\begin{scope}[shift={(-1/3,-1/3)}]
\draw [thick,black] plot  coordinates {(-1/12,1/12) (1/12,-1/12) };
\end{scope}
\begin{scope}[shift={(-2/3,-2/3)}]
\draw [thick,black] plot  coordinates {(-1/12,1/12) (1/12,-1/12) };
\end{scope}
\begin{scope}[shift={(-3/3,-3/3)}]
\draw [thick,black] plot  coordinates {(-1/12,1/12) (1/12,-1/12) };
\end{scope}
\begin{scope}[shift={(1/2,-1/2)}]
\draw [thick,black] plot  coordinates {(-1/12,-1/12) (1/12,1/12) };
\end{scope}
\begin{scope}[shift={(2-1/6,-1/6)}]
\draw [thick,black] plot  coordinates {(-1/12,1/12) (1/12,-1/12) };
\end{scope}
\begin{scope}[shift={(2-3/6,-3/6)}]
\draw [thick,black] plot  coordinates {(-1/12,1/12) (1/12,-1/12) };
\end{scope}
\end{tikzpicture}}

\def\FigTwoAllD{\begin{tikzpicture}[scale=0.7]
\draw [thick,black] plot  coordinates {(4/8,-1/8) (-20/8,5/8) };
\begin{scope}[shift={(0,0)}]
\draw [thick,black] plot  coordinates {(2/24,8/24) (-5/16,-20/16) };
\fill [black] (0,0) circle (0.08);
\begin{scope}[shift={(-3/32,-12/32)}]
\draw [thick,black] plot  coordinates {(-1/8,1/32) (1/8,-1/32) };
\end{scope}
\begin{scope}[shift={(-6/32,-24/32)}]
\draw [thick,black] plot  coordinates {(-1/8,1/32) (1/8,-1/32) };
\end{scope}
\end{scope}
\begin{scope}[shift={(-4/4,1/4)}]
\draw [thick,black] plot  coordinates {(2/24,8/24) (-5/16,-20/16) };
\fill [black] (0,0) circle (0.08);
\begin{scope}[shift={(-3/32,-12/32)}]
\draw [thick,black] plot  coordinates {(-1/8,1/32) (1/8,-1/32) };
\end{scope}
\begin{scope}[shift={(-6/32,-24/32)}]
\draw [thick,black] plot  coordinates {(-1/8,1/32) (1/8,-1/32) };
\end{scope}
\end{scope}
\begin{scope}[shift={(-8/4,2/4)}]
\draw [thick,black] plot  coordinates {(2/24,8/24) (-5/16,-20/16) };
\fill [black] (0,0) circle (0.08);
\begin{scope}[shift={(-3/32,-12/32)}]
\draw [thick,black] plot  coordinates {(-1/8,1/32) (1/8,-1/32) };
\end{scope}
\begin{scope}[shift={(-6/32,-24/32)}]
\draw [thick,black] plot  coordinates {(-1/8,1/32) (1/8,-1/32) };
\end{scope}
\end{scope}
\end{tikzpicture}}

\def\FigTwoAllG{\begin{tikzpicture}[scale=0.7]
\draw [thick,black] plot  coordinates {(-10/8,-10/8) (1/4,1/4) };
\draw [thick,black] plot  coordinates {(8/8,-8/8) (-1/4,1/4) };
\draw [thick,black] plot  coordinates {(4/8,-8/8) (14/8,2/8) };
\draw [thick,black] plot  coordinates {(22/8,-10/8) (5/4,1/4) };
\fill [black] (0,0) circle (0.08);
\begin{scope}[shift={(3/4,-3/4)}]
\fill [white] (0,0) circle (0.1);
\draw [thick] (0,0) circle (0.1);
\end{scope}
\fill [black] (6/4,0) circle (0.08);
\begin{scope}[shift={(-1/2,-1/2)}]
\draw [thick,black] plot  coordinates {(-1/12,1/12) (1/12,-1/12) };
\end{scope}
\begin{scope}[shift={(-5/6,-5/6)}]
\draw [thick,black] plot  coordinates {(-1/12,1/12) (1/12,-1/12) };
\end{scope}
\begin{scope}[shift={(3/8,-3/8)}]
\draw [thick,black] plot  coordinates {(-1/12,-1/12) (1/12,1/12) };
\end{scope}
\begin{scope}[shift={(3/2-3/8,-3/8)}]
\draw [thick,black] plot  coordinates {(-1/12,1/12) (1/12,-1/12) };
\end{scope}
\begin{scope}[shift={(2,-1/2)}]
\draw [thick,black] plot  coordinates {(-1/12,-1/12) (1/12,1/12) };
\end{scope}
\begin{scope}[shift={(2+2/6,-5/6)}]
\draw [thick,black] plot  coordinates {(-1/12,-1/12) (1/12,1/12) };
\end{scope}
\end{tikzpicture}}

\def\includefig#1#2{
\vcenter{\hbox{\makebox(110pt,55pt){
\fbox{\makebox(0pt,45pt)[lt]{\makebox(20pt,20pt){(\small #1)}}
\makebox(90pt,45pt){#2}}}}}}

\def\FigTwoAll{
$\xymatrix@C+0pt@R+0pt{
\includefig{A}{\FigTwoAllA} \ar[d] \ar[dr] \\
\includefig{B}{\FigTwoAllB} \ar[d] \ar[dr] & \includefig{E}{\FigTwoAllE} \ar[d] \\
\includefig{C}{\FigTwoAllC}  \ar[d] \ar[dr] & \includefig{F}{\FigTwoAllF} \ar[d] \\
\includefig{D}{\hskip15pt \FigTwoAllD} & \includefig{G}{\FigTwoAllG} }$}

\def\FigTwoCaseAa{\begin{tikzpicture}[baseline=25pt]
\draw (7/4,0/4) node {};
\draw [thick,black] plot  coordinates {(-15/8,-15/8) (3/4,3/4) };
\draw [thick,cbblue] plot  coordinates {(12/8,-10/8) (-1/8,3/8) };
\fill [black] (1/8,1/8) circle (0.06);
\draw (1/8,1/8-1/16) node [below,red] {$\frac{1}{3}$};
\draw (11/8-2/3,-11/8+2/4) node [below] {\hbox{\footnotesize$\geqone$}};
\begin{scope}[shift={(2/4,2/4)}]
\draw [thick,cbblue] plot  coordinates {(11/8,-11/8) (-1/4,1/4) };
\fill [black] (0,0) circle (0.06);
\draw (0/8,0/8+1/16) node [above,red] {$\frac{1}{3}$};
\draw (11/8,-11/8+5/12) node [above] {\hbox{\footnotesize$\geqone$}};
\end{scope}
\begin{scope}[shift={(-1/4,-1/4)}]
\draw [thick,black] plot  coordinates {(-1/12,1/12) (1/12,-1/12) };
\end{scope}
\begin{scope}[shift={(-2/4,-2/4)}]
\draw [thick,black] plot  coordinates {(-1/12,1/12) (1/12,-1/12) };
\end{scope}
\begin{scope}[shift={(-3/4,-3/4)}]
\draw [thick,black] plot  coordinates {(-1/12,1/12) (1/12,-1/12) };
\end{scope}
\begin{scope}[shift={(-4/4,-4/4)}]
\draw [thick,black] plot  coordinates {(-1/12,1/12) (1/12,-1/12) };
\end{scope}
\begin{scope}[shift={(-5/4,-5/4)}]
\draw [thick,black] plot  coordinates {(-1/12,1/12) (1/12,-1/12) };
\end{scope}
\begin{scope}[shift={(-6/4,-6/4)}]
\draw [thick,black] plot  coordinates {(-1/12,1/12) (1/12,-1/12) };
\end{scope}
\end{tikzpicture}}

\def\FigTwoCaseAb{\begin{tikzpicture}[baseline=23pt]
\draw (7/4,0/4) node {};
\draw [thick,black] plot  coordinates {(-14/8,-14/8) (4/8,4/8) };
\draw [thick,cbgreen] plot  coordinates {(15/8,-11/8) (0,2/4) };
\fill [black] (1/4,1/4) circle (0.06);
\draw (2/8,2/8+1/16) node [above,red] {$\frac{\delta}{4}$};
\begin{scope}[shift={(3/4,-3/4)}]
\draw [thick,cbblue] plot  coordinates {(-5/8,-5/8) (2/4,2/4) };
\fill [black] (1/4,1/4) circle (0.06);
\draw (2/8,3/8) node [above,red] {$\frac{4-5\delta}{12}$};
\draw (-5/8-0.20,-5/8) node [below] {\hbox{\footnotesize$\geqone$}};
\end{scope}
\begin{scope}[shift={(5/4,-5/4)}]
\draw [thick,cbblue] plot  coordinates {(-5/8,-5/8) (2/4,2/4) };
\fill [black] (1/4,1/4) circle (0.06);
\draw (3/8,2/8) node [right,red] {$\frac{4-5\delta}{12}$};
\draw (-5/8+0.05,-5/8) node [right] {\hbox{\footnotesize$\geqone$}};
\end{scope}
\begin{scope}[shift={(1/8,1/8)}]
\begin{scope}[shift={(-1/4,-1/4)}]
\draw [thick,black] plot  coordinates {(-1/12,1/12) (1/12,-1/12) };
\end{scope}
\begin{scope}[shift={(-2/4,-2/4)}]
\draw [thick,black] plot  coordinates {(-1/12,1/12) (1/12,-1/12) };
\end{scope}
\begin{scope}[shift={(-3/4,-3/4)}]
\draw [thick,black] plot  coordinates {(-1/12,1/12) (1/12,-1/12) };
\end{scope}
\begin{scope}[shift={(-4/4,-4/4)}]
\draw [thick,black] plot  coordinates {(-1/12,1/12) (1/12,-1/12) };
\end{scope}
\begin{scope}[shift={(-5/4,-5/4)}]
\draw [thick,black] plot  coordinates {(-1/12,1/12) (1/12,-1/12) };
\end{scope}
\begin{scope}[shift={(-6/4,-6/4)}]
\draw [thick,black] plot  coordinates {(-1/12,1/12) (1/12,-1/12) };
\end{scope}
\end{scope}
\end{tikzpicture}}

\def\FigTwoCaseAc{\begin{tikzpicture}[baseline=20pt]
\draw (7/4,0/4) node {};
\draw [thick,black] plot  coordinates {(-14/8,-14/8) (4/8,4/8) };
\begin{scope}[shift={(1/4,1/4)}]
\draw [thick,cbblue] plot  coordinates {(11/8,-11/8) (-1/4,1/4) };
\fill [black] (0,0) circle (0.06);
\draw (0/8,0/8+1/16) node [above,red] {$\frac{1}{5}$};
\draw (11/8-2/3,-11/8+1/3) node [below] {\hbox{\footnotesize$\geqtwo$}};
\end{scope}
\begin{scope}[shift={(1/8,1/8)}]
\begin{scope}[shift={(-1/4,-1/4)}]
\draw [thick,black] plot  coordinates {(-1/12,1/12) (1/12,-1/12) };
\end{scope}
\begin{scope}[shift={(-2/4,-2/4)}]
\draw [thick,black] plot  coordinates {(-1/12,1/12) (1/12,-1/12) };
\end{scope}
\begin{scope}[shift={(-3/4,-3/4)}]
\draw [thick,black] plot  coordinates {(-1/12,1/12) (1/12,-1/12) };
\end{scope}
\begin{scope}[shift={(-4/4,-4/4)}]
\draw [thick,black] plot  coordinates {(-1/12,1/12) (1/12,-1/12) };
\end{scope}
\begin{scope}[shift={(-5/4,-5/4)}]
\draw [thick,black] plot  coordinates {(-1/12,1/12) (1/12,-1/12) };
\end{scope}
\begin{scope}[shift={(-6/4,-6/4)}]
\draw [thick,black] plot  coordinates {(-1/12,1/12) (1/12,-1/12) };
\end{scope}
\end{scope}
\end{tikzpicture}}

\def\includefigA#1{
\makebox(130pt,90pt){#1}}

\def\FigTwoCaseA{
\begin{tabular}{|c|c|c|}
	\hline
{\large\strut}	(A1) & (A2) & (A3) \\
\hline
{\Large\strut} $\delta=0$ & $0<\delta<4/5$ & $\delta\ge4/5$ \\
\hline
\includefigA{\FigTwoCaseAa} & 
\includefigA{\FigTwoCaseAb} & 
\includefigA{\FigTwoCaseAc} \\
\hline
\end{tabular}}

\def\FigTwoCaseBi{\begin{tikzpicture}[baseline=13pt]
\draw [thick,black] plot  coordinates {(-16/12,-16/8) (1/6,1/4) };
\begin{scope}[shift={(-2/9,-3/9)}]
\draw [thick,black] plot  coordinates {(-3/32,2/32) (3/32,-2/32) };
\end{scope}
\begin{scope}[shift={(-4/9,-6/9)}]
\draw [thick,black] plot  coordinates {(-3/32,2/32) (3/32,-2/32) };
\end{scope}
\begin{scope}[shift={(-6/9,-9/9)}]
\draw [thick,black] plot  coordinates {(-3/32,2/32) (3/32,-2/32) };
\end{scope}
\begin{scope}[shift={(-8/9,-12/9)}]
\draw [thick,black] plot  coordinates {(-3/32,2/32) (3/32,-2/32) };
\end{scope}
\draw [thick,cborange] plot  coordinates {(-1/2,0) (4/2,0) };
\fill [black] (0,0) circle (0.06);
\draw (0,4/16)  node [left,red] {$\frac{\alpha}{4}$};
\draw [thick,cbblue] plot  coordinates {(3/4,1/3) (3/4,-6/4) };
\draw (3/4,-3/2) node [right] {\hbox{\footnotesize$\geqone$}};
\fill [black] (3/4,0) circle (0.06);
\draw (3/4+2/16,-6/16)  node [left,red] {$\frac{4-\alpha}{12}$};
\begin{scope}[shift={(3/2,0)}]
\draw [thick,black] plot  coordinates {(8/12,-8/8) (-1/6,1/4) };
\fill [black] (0,0) circle (0.06);
\draw (0,4/16)  node [right,red] {$\frac{\alpha}{4}$};
\begin{scope}[shift={(2/9,-3/9)}]
\draw [thick,black] plot  coordinates {(3/32,2/32) (-3/32,-2/32) };
\end{scope}
\begin{scope}[shift={(4/9,-6/9)}]
\draw [thick,black] plot  coordinates {(3/32,2/32) (-3/32,-2/32) };
\end{scope}
\end{scope}
\draw [thick,cbblue] plot  coordinates {(-10/9-5/6,-15/9+5/4) (-10/9+1/6,-15/9-1/4) };
\fill [black] (-10/9,-15/9) circle (0.06);
\draw (-10/9-1/16,-15/9) node [left,red] {$\frac{1}{3}$};
\draw (-10/9-5/6,-15/9+5/4+1/6) node [right] {\hbox{\footnotesize$\geqone$}};
\end{tikzpicture}}

\def\FigTwoCaseBii{\begin{tikzpicture}[baseline=13pt]
\draw [thick,black] plot  coordinates {(-16/12,-16/8) (1/6,1/4) };
\begin{scope}[shift={(-2/9,-3/9)}]
\draw [thick,black] plot  coordinates {(-3/32,2/32) (3/32,-2/32) };
\end{scope}
\begin{scope}[shift={(-4/9,-6/9)}]
\draw [thick,black] plot  coordinates {(-3/32,2/32) (3/32,-2/32) };
\end{scope}
\begin{scope}[shift={(-6/9,-9/9)}]
\draw [thick,black] plot  coordinates {(-3/32,2/32) (3/32,-2/32) };
\end{scope}
\begin{scope}[shift={(-8/9,-12/9)}]
\draw [thick,black] plot  coordinates {(-3/32,2/32) (3/32,-2/32) };
\end{scope}
\draw [thick,cborange] plot  coordinates {(-1/2,0) (4/2,0) };
\draw (3/4,0) node [below] {\hbox{\footnotesize$\geqone$}};
\fill [black] (0,0) circle (0.06);
\draw (0,4/16)  node [left,red] {\scriptsize$1$};
\begin{scope}[shift={(3/2,0)}]
\draw [thick,black] plot  coordinates {(8/12,-8/8) (-1/6,1/4) };
\fill [black] (0,0) circle (0.06);
\draw (0,4/16)  node [right,red] {\scriptsize$1$};
\begin{scope}[shift={(2/9,-3/9)}]
\draw [thick,black] plot  coordinates {(3/32,2/32) (-3/32,-2/32) };
\end{scope}
\begin{scope}[shift={(4/9,-6/9)}]
\draw [thick,black] plot  coordinates {(3/32,2/32) (-3/32,-2/32) };
\end{scope}
\end{scope}
\draw [thick,cbblue] plot  coordinates {(-10/9-5/6,-15/9+5/4) (-10/9+1/6,-15/9-1/4) };
\fill [black] (-10/9,-15/9) circle (0.06);
\draw (-10/9-1/16,-15/9) node [left,red] {$\frac{1}{3}$};
\draw (-10/9-5/6,-15/9+5/4+1/6) node [right] {\hbox{\footnotesize$\geqone$}};
\end{tikzpicture}}

\def\FigTwoCaseBiii{\begin{tikzpicture}[baseline=13pt]
\draw [thick,black] plot  coordinates {(-16/12,-16/8) (1/6,1/4) };
\begin{scope}[shift={(-2/9,-3/9)}]
\draw [thick,black] plot  coordinates {(-3/32,2/32) (3/32,-2/32) };
\end{scope}
\begin{scope}[shift={(-4/9,-6/9)}]
\draw [thick,black] plot  coordinates {(-3/32,2/32) (3/32,-2/32) };
\end{scope}
\begin{scope}[shift={(-6/9,-9/9)}]
\draw [thick,black] plot  coordinates {(-3/32,2/32) (3/32,-2/32) };
\end{scope}
\begin{scope}[shift={(-8/9,-12/9)}]
\draw [thick,black] plot  coordinates {(-3/32,2/32) (3/32,-2/32) };
\end{scope}
\begin{scope}[shift={(3/4,0)}]
\draw [thick,cborange] plot [smooth, tension=1] coordinates {(-3/3,2/12) (1/4,-6/12) (-1/3,-14/12) };
\draw [thick,cborange] plot [smooth, tension=1] coordinates {(3/3,2/12) (-1/4,-6/12) (1/3,-14/12) };
\fill [black] (0,-0.26) circle (0.06);
\fill [black] (0,-0.99) circle (0.06);
\draw (0,-0.26+1/16) node [above,red] {\scriptsize $\alpha\text{-}4$};
\draw (0,-0.99-1/16) node [below,red] {\scriptsize $\alpha\text{-}4$};
\fill [black] (-3/4+0.14/3,0.14/2) circle (0.06);
\fill [black] (3/4-0.14/3,0.14/2) circle (0.06);
\draw (3/4-2/16,6/16)  node [right,red] {\scriptsize $1$};
\draw (-3/4+2/16,6/16)  node [left,red] {\scriptsize $1$};
\end{scope}
\begin{scope}[shift={(3/2,0)}]
\draw [thick,black] plot  coordinates {(8/12,-8/8) (-1/6,1/4) };
\begin{scope}[shift={(2/9,-3/9)}]
\draw [thick,black] plot  coordinates {(3/32,2/32) (-3/32,-2/32) };
\end{scope}
\begin{scope}[shift={(4/9,-6/9)}]
\draw [thick,black] plot  coordinates {(3/32,2/32) (-3/32,-2/32) };
\end{scope}
\end{scope}
\draw [thick,cbblue] plot  coordinates {(-10/9-5/6,-15/9+5/4) (-10/9+1/6,-15/9-1/4) };
\fill [black] (-10/9,-15/9) circle (0.06);
\draw (-10/9-1/16,-15/9) node [left,red] {$\frac{1}{3}$};
\draw (-10/9-5/6,-15/9+5/4+1/6) node [right] {\hbox{\footnotesize$\geqone$}};
\end{tikzpicture}}

\def\FigTwoCaseBvii{\begin{tikzpicture}[baseline=13pt]
\draw [thick,black] plot  coordinates {(-14/12,-14/8) (1/6,1/4) };
\begin{scope}[shift={(-2/9,-3/9)}]
\draw [thick,black] plot  coordinates {(-3/32,2/32) (3/32,-2/32) };
\end{scope}
\begin{scope}[shift={(-4/9,-6/9)}]
\draw [thick,black] plot  coordinates {(-3/32,2/32) (3/32,-2/32) };
\end{scope}
\begin{scope}[shift={(-6/9,-9/9)}]
\draw [thick,black] plot  coordinates {(-3/32,2/32) (3/32,-2/32) };
\end{scope}
\begin{scope}[shift={(-8/9,-12/9)}]
\draw [thick,black] plot  coordinates {(-3/32,2/32) (3/32,-2/32) };
\end{scope}
\begin{scope}[shift={(3/4,0)}]
\draw [thick,cborange] plot [smooth, tension=1] coordinates {(-3/3,2/12) (1,-8/12) (-1/2,-20/12) };
\draw [thick,cborange] plot [smooth, tension=.8] coordinates {(1/2+3/3,2/12) (5/12,-8/12) (11/12,-16/12) };
\fill [black] (0.69,-0.36) circle (0.06);
\draw (0.69,-0.36) node [above,red] { $\frac{3\alpha\text{-}8}{3}$};
\fill [black] (0.66,-1.16) circle (0.06);
\draw (0.66,-1.16-2/16) node [below,red] { $\frac{3\alpha\text{-}8}{3}$};
\draw (-1/2,-20/12) node [above] {\hbox{\footnotesize$\geqone$}};
\end{scope}
\fill [black] (0.195/3,0.195/2) circle (0.06);
\draw (0.195/3,0.195/2-1/16) node [below,red] { $\frac{1}{3}$};
\begin{scope}[shift={(2,0)}]
\draw [thick,black] plot  coordinates {(8/12,-8/8) (-1/6,1/4) };
\fill [black] (-0.035/3,0.035/2) circle (0.06);
\draw (-0.035/3,0.035/2) node [below,red] { \scriptsize $1$};
\begin{scope}[shift={(2/9,-3/9)}]
\draw [thick,black] plot  coordinates {(3/32,2/32) (-3/32,-2/32) };
\end{scope}
\begin{scope}[shift={(4/9,-6/9)}]
\draw [thick,black] plot  coordinates {(3/32,2/32) (-3/32,-2/32) };
\end{scope}
\end{scope}
\end{tikzpicture}}

\def\FigTwoCaseBvi{\begin{tikzpicture}[baseline=13pt]
\draw [thick,black] plot  coordinates {(-1-14/12,-14/8) (-1+1/6,1/4) };
\begin{scope}[shift={(-3/36,-4/36)}]
\begin{scope}[shift={(-1-2/9,-3/9)}]
\draw [thick,black] plot  coordinates {(-3/32,2/32) (3/32,-2/32) };
\end{scope}
\begin{scope}[shift={(-1-4/9,-6/9)}]
\draw [thick,black] plot  coordinates {(-3/32,2/32) (3/32,-2/32) };
\end{scope}
\begin{scope}[shift={(-1-6/9,-9/9)}]
\draw [thick,black] plot  coordinates {(-3/32,2/32) (3/32,-2/32) };
\end{scope}
\begin{scope}[shift={(-1-8/9,-12/9)}]
\draw [thick,black] plot  coordinates {(-3/32,2/32) (3/32,-2/32) };
\end{scope}
\end{scope}
\draw [thick,cborange] plot [smooth, tension=1] coordinates {(4/12,1/4) (-6/12,-1/4) (-16/12,1/4) };
\begin{scope}[shift={(3/4,0)}]
\draw [thick,cborange] plot [smooth, tension=1] coordinates {(-3/3,2/12) (1/4,-6/12) (-1/3,-14/12) };
\draw [thick,cborange] plot [smooth, tension=1] coordinates {(3/3,2/12) (-1/4,-6/12) (1/3,-14/12) };
\fill [black] (0,-0.26) circle (0.06);
\draw (0,-0.26+1/16)  node [above,red] {\scriptsize $\alpha\texttt{+}\delta \text{-}4$};
\fill [black] (0,-0.99) circle (0.06);
\draw (0,-0.99-1/16)  node [below,red] {\scriptsize $\alpha\texttt{+}\delta \text{-}4$};
\fill [black] (3/4-0.14/3,0.14/2) circle (0.06);
\draw (3/4-0.14/3+2/16,0.14/2+1/16)  node [above,red] {\scriptsize $1$};
\end{scope}
\fill [black] (-5/6-0.6/3,1/4-0.6/2) circle (0.06);
\draw (-5/6-0.6/3-1/16,1/4-0.6/2+1/16)  node [above,red] {$\frac{\delta}{4}$};
\fill [black] (0.14,0.036) circle (0.06);
\draw (0.14,0.036)  node [above,red] {$\frac{4-3\delta}{4}$};
\begin{scope}[shift={(3/2,0)}]
\draw [thick,black] plot  coordinates {(8/12,-8/8) (-1/6,1/4) };
\begin{scope}[shift={(2/9,-3/9)}]
\draw [thick,black] plot  coordinates {(3/32,2/32) (-3/32,-2/32) };
\end{scope}
\begin{scope}[shift={(4/9,-6/9)}]
\draw [thick,black] plot  coordinates {(3/32,2/32) (-3/32,-2/32) };
\end{scope}
\end{scope}
\draw [thick,cbblue] plot  coordinates {(-6/12,-15/9+7/4) (-6/12,-15/9+1/4) };
\fill [black] (-6/12,-.25) circle (0.06);
\draw (-6/12+1/16,-.5-2/16)  node [left,red] {$\frac{4\text{-}3\delta}{12}$};
\draw (-6/6,-10/6) node [right] {\hbox{\footnotesize$\geqone$}};
\end{tikzpicture}}

\def\FigTwoCaseBiv{\begin{tikzpicture}[baseline=13pt]
\draw [thick,black] plot  coordinates {(-14/12,-14/8) (1/6,1/4) };
\begin{scope}[shift={(-2/9,-3/9)}]
\draw [thick,black] plot  coordinates {(-3/32,2/32) (3/32,-2/32) };
\end{scope}
\begin{scope}[shift={(-4/9,-6/9)}]
\draw [thick,black] plot  coordinates {(-3/32,2/32) (3/32,-2/32) };
\end{scope}
\begin{scope}[shift={(-6/9,-9/9)}]
\draw [thick,black] plot  coordinates {(-3/32,2/32) (3/32,-2/32) };
\end{scope}
\begin{scope}[shift={(-8/9,-12/9)}]
\draw [thick,black] plot  coordinates {(-3/32,2/32) (3/32,-2/32) };
\end{scope}
\draw [thick,cborange] plot  coordinates {(-1/4,1/6) (6/4,-6/6) };
\draw [thick,cborange] plot  coordinates {(5/2+1/4,1/6) (5/2-6/4,-6/6) };
\fill [black] (0,0) circle (0.06);
\draw (-3/16,0)  node [left,red] {$\frac{\delta}{4}$};
\fill [black] (5/4,-5/6) circle (0.06);
\draw (5/4,-5/6-2/16)  node [below,red] {$\frac{\alpha\text{-}2\delta}{4}$};
\begin{scope}[shift={(3/5,-2/5)}]
\draw [thick,cbblue] plot coordinates {(1/12,1/3) (-7/24,-7/6) };
\fill [black] (0,0) circle (0.06);
\draw (0,-6/16)  node [left,red] {$\frac{4\text{-}3\delta}{12}$};
\draw (-6.5/24,-6.5/6) node [below] {\hbox{\footnotesize$\geqone$}};
\end{scope}
\begin{scope}[shift={(5/2-3/5,-2/5)}]
\draw [thick,cbblue] plot coordinates {(-1/12,1/3) (7/24,-7/6) };
\fill [black] (0,0) circle (0.06);
\draw (2/16,6/16)  node [left,red] {$\frac{4\text{-}\alpha\text{-}\delta}{12}$};
\draw (6.5/24,-6.5/6) node [below] {\hbox{\footnotesize$\geqone$}};
\end{scope}
\begin{scope}[shift={(5/2,0)}]
\draw [thick,black] plot  coordinates {(8/12,-8/8) (-1/6,1/4) };
\fill [black] (0,0) circle (0.06);
\draw (1/16,-1/16)  node [right,red] {$\frac{\alpha+\delta}{4}$};
\begin{scope}[shift={(2/9,-3/9)}]
\draw [thick,black] plot  coordinates {(3/32,2/32) (-3/32,-2/32) };
\end{scope}
\begin{scope}[shift={(4/9,-6/9)}]
\draw [thick,black] plot  coordinates {(3/32,2/32) (-3/32,-2/32) };
\end{scope}
\end{scope}
\end{tikzpicture}}

\def\FigTwoCaseBv{\begin{tikzpicture}[baseline=13pt]
\draw [thick,black] plot  coordinates {(-14/12,-14/8) (1/6,1/4) };
\begin{scope}[shift={(-2/9,-3/9)}]
\draw [thick,black] plot  coordinates {(-3/32,2/32) (3/32,-2/32) };
\end{scope}
\begin{scope}[shift={(-4/9,-6/9)}]
\draw [thick,black] plot  coordinates {(-3/32,2/32) (3/32,-2/32) };
\end{scope}
\begin{scope}[shift={(-6/9,-9/9)}]
\draw [thick,black] plot  coordinates {(-3/32,2/32) (3/32,-2/32) };
\end{scope}
\begin{scope}[shift={(-8/9,-12/9)}]
\draw [thick,black] plot  coordinates {(-3/32,2/32) (3/32,-2/32) };
\end{scope}
\draw [thick,cborange] plot  coordinates {(-1/4,1/6) (6/4,-6/6) };
\draw [thick,cborange] plot  coordinates {(5/2+1/4,1/6) (5/2-6/4,-6/6) };
\fill [black] (0,0) circle (0.06);
\draw (-3/16,0)  node [left,red] {$\frac{\delta}{4}$};
\fill [black] (5/4,-5/6) circle (0.06);
\draw (5/4,-5/6-2/16)  node [below,red] {$\frac{\alpha\text{-}2\delta}{4}$};
\begin{scope}[shift={(3/5,-2/5)}]
\draw [thick,cbblue] plot coordinates {(1/12,1/3) (-7/24,-7/6) };
\fill [black] (0,0) circle (0.06);
\draw (0,-6/16)  node [left,red] {$\frac{4\text{-}3\delta}{12}$};
\draw (-6.5/24,-6.5/6) node [below] {\hbox{\footnotesize$\geqone$}};
\end{scope}
\draw (1+4/4+1/12,-1+4/6) node [below] {\hbox{\footnotesize$\geqone$}};
\begin{scope}[shift={(5/2,0)}]
\draw [thick,black] plot  coordinates {(8/12,-8/8) (-1/6,1/4) };
\fill [black] (0,0) circle (0.06);
\draw (1/16,-1/16)  node [right,red] {$\frac{\alpha+\delta}{4}$};
\begin{scope}[shift={(2/9,-3/9)}]
\draw [thick,black] plot  coordinates {(3/32,2/32) (-3/32,-2/32) };
\end{scope}
\begin{scope}[shift={(4/9,-6/9)}]
\draw [thick,black] plot  coordinates {(3/32,2/32) (-3/32,-2/32) };
\end{scope}
\end{scope}
\end{tikzpicture}}

\def\FigTwoCaseBviii{\begin{tikzpicture}[baseline=15pt]
\draw [thick,black] plot  coordinates {(-14/12,-14/8) (1/6,1/4) };
\begin{scope}[shift={(-2/9,-3/9)}]
\draw [thick,black] plot  coordinates {(-3/32,2/32) (3/32,-2/32) };
\end{scope}
\begin{scope}[shift={(-4/9,-6/9)}]
\draw [thick,black] plot  coordinates {(-3/32,2/32) (3/32,-2/32) };
\end{scope}
\begin{scope}[shift={(-6/9,-9/9)}]
\draw [thick,black] plot  coordinates {(-3/32,2/32) (3/32,-2/32) };
\end{scope}
\begin{scope}[shift={(-8/9,-12/9)}]
\draw [thick,black] plot  coordinates {(-3/32,2/32) (3/32,-2/32) };
\end{scope}
\draw [thick,cborange] plot  coordinates {(-1/2,0) (5/2,0) };
\fill [black] (0,0) circle (0.06);
\draw (0-1/16,0) node [above,red] { $\frac{\alpha}{8}$};
\begin{scope}[shift={(4/2,0)}]
\draw [thick,black] plot  coordinates {(8/12,-8/8) (-1/6,1/4) };
\fill [black] (0,0) circle (0.06);
\draw (0+1/16,0) node [above,red] { $\frac{3\alpha}{8}$};
\begin{scope}[shift={(2/9,-3/9)}]
\draw [thick,black] plot  coordinates {(3/32,2/32) (-3/32,-2/32) };
\end{scope}
\begin{scope}[shift={(4/9,-6/9)}]
\draw [thick,black] plot  coordinates {(3/32,2/32) (-3/32,-2/32) };
\end{scope}
\end{scope}
\begin{scope}[shift={(2/3,0)}]
\draw [thick,cbblue] plot  coordinates {(0,1/3) (0,-6/4) };
\fill [black] (0,0) circle (0.06);
\draw (0-6/16,0) node [below,red] { $\frac{8\text{-}3\alpha}{24}$};
\draw (0,-6/4) node [left] {\hbox{\footnotesize$\geqone$\kern-1pt}};
\end{scope}
\begin{scope}[shift={(4/3,0)}]
\draw [thick,cbblue] plot  coordinates {(0,1/3) (0,-6/4) };
\fill [black] (0,0) circle (0.06);
\draw (0+6/16,0) node [below,red] { $\frac{8\text{-}3\alpha}{24}$};
\draw (0,-6/4) node [right] {\hbox{\footnotesize$\geqone$}};
\end{scope}
\end{tikzpicture}}

\def\FigTwoCaseBix{\begin{tikzpicture}[baseline=15pt]
\draw [thick,black] plot  coordinates {(-14/12,-14/8) (1/6,1/4) };
\begin{scope}[shift={(-2/9,-3/9)}]
\draw [thick,black] plot  coordinates {(-3/32,2/32) (3/32,-2/32) };
\end{scope}
\begin{scope}[shift={(-4/9,-6/9)}]
\draw [thick,black] plot  coordinates {(-3/32,2/32) (3/32,-2/32) };
\end{scope}
\begin{scope}[shift={(-6/9,-9/9)}]
\draw [thick,black] plot  coordinates {(-3/32,2/32) (3/32,-2/32) };
\end{scope}
\begin{scope}[shift={(-8/9,-12/9)}]
\draw [thick,black] plot  coordinates {(-3/32,2/32) (3/32,-2/32) };
\end{scope}
\draw [thick,cborange] plot  coordinates {(-1/2,0) (5/2,0) };
\fill [black] (0,0) circle (0.06);
\draw (0-1/16,0) node [above,red] { $\frac{\alpha}{8}$};
\begin{scope}[shift={(4/2,0)}]
\draw [thick,black] plot  coordinates {(8/12,-8/8) (-1/6,1/4) };
\fill [black] (0,0) circle (0.06);
\draw (0+1/16,0) node [above,red] { $\frac{3\alpha}{8}$};
\begin{scope}[shift={(2/9,-3/9)}]
\draw [thick,black] plot  coordinates {(3/32,2/32) (-3/32,-2/32) };
\end{scope}
\begin{scope}[shift={(4/9,-6/9)}]
\draw [thick,black] plot  coordinates {(3/32,2/32) (-3/32,-2/32) };
\end{scope}
\end{scope}
\begin{scope}[shift={(1/2,0)}]
\draw [thick,cbgreen] plot  coordinates {(0,1/3) (0,-7/4) };
\fill [black] (0,0) circle (0.06);
\draw (0+4/16,0) node [above,red] { $\frac{2\delta\text{-}\alpha}{8}$};
\end{scope}
\begin{scope}[shift={(1/2,-2/3)}]
\draw [thick,cbblue] plot  coordinates {(-1/3,0) (4/3,0)};
\fill [black] (0,0) circle (0.06);
\draw (0+9/16,-1/16) node [above,red] { $\frac{4\text{-}5\delta\texttt{+}\alpha}{12}$};
\draw (2/2+9/16,0) node [below] {\hbox{\footnotesize$\geqone$\kern-1pt}};
\end{scope}
\begin{scope}[shift={(1/2,-4/3)}]
\draw [thick,cbblue] plot  coordinates {(-1/3,0) (4/3,0)};
\fill [black] (0,0) circle (0.06);
\draw (0+9/16,-1/16) node [above,red] { $\frac{4\text{-}5\delta\texttt{+}\alpha}{12}$};
\draw (2/2,0) node [below] {\hbox{\footnotesize$\geqone$\kern-1pt}};
\end{scope}
\end{tikzpicture}}

\def\FigTwoCaseBx{\begin{tikzpicture}[baseline=15pt]
\draw [thick,black] plot  coordinates {(-14/12,-14/8) (1/6,1/4) };
\begin{scope}[shift={(-2/9,-3/9)}]
\draw [thick,black] plot  coordinates {(-3/32,2/32) (3/32,-2/32) };
\end{scope}
\begin{scope}[shift={(-4/9,-6/9)}]
\draw [thick,black] plot  coordinates {(-3/32,2/32) (3/32,-2/32) };
\end{scope}
\begin{scope}[shift={(-6/9,-9/9)}]
\draw [thick,black] plot  coordinates {(-3/32,2/32) (3/32,-2/32) };
\end{scope}
\begin{scope}[shift={(-8/9,-12/9)}]
\draw [thick,black] plot  coordinates {(-3/32,2/32) (3/32,-2/32) };
\end{scope}
\draw [thick,cborange] plot  coordinates {(-1/2,0) (5/2,0) };
\fill [black] (0,0) circle (0.06);
\draw (0-1/16,0) node [above,red] { $\frac{\alpha}{8}$};
\begin{scope}[shift={(4/2,0)}]
\draw [thick,black] plot  coordinates {(8/12,-8/8) (-1/6,1/4) };
\fill [black] (0,0) circle (0.06);
\draw (0+1/16,0) node [above,red] { $\frac{3\alpha}{8}$};
\begin{scope}[shift={(2/9,-3/9)}]
\draw [thick,black] plot  coordinates {(3/32,2/32) (-3/32,-2/32) };
\end{scope}
\begin{scope}[shift={(4/9,-6/9)}]
\draw [thick,black] plot  coordinates {(3/32,2/32) (-3/32,-2/32) };
\end{scope}
\end{scope}
\begin{scope}[shift={(1,0)}]
\draw [thick,cbblue] plot  coordinates {(0,1/3) (0,-6/4) };
\fill [black] (0,0) circle (0.06);
\draw (0-4/16,0) node [above,red] { $\frac{8\text{-}3\alpha}{40}$};
\draw (0,-6/4) node [right] {\hbox{\footnotesize$\geqtwo$}};
\end{scope}
\end{tikzpicture}}

\def\FigTwoCaseBxi{\begin{tikzpicture}[baseline=15pt]
\draw [thick,black] plot  coordinates {(-14/12,-14/8) (1/6,1/4) };
\begin{scope}[shift={(-2/9,-3/9)}]
\draw [thick,black] plot  coordinates {(-3/32,2/32) (3/32,-2/32) };
\end{scope}
\begin{scope}[shift={(-4/9,-6/9)}]
\draw [thick,black] plot  coordinates {(-3/32,2/32) (3/32,-2/32) };
\end{scope}
\begin{scope}[shift={(-6/9,-9/9)}]
\draw [thick,black] plot  coordinates {(-3/32,2/32) (3/32,-2/32) };
\end{scope}
\begin{scope}[shift={(-8/9,-12/9)}]
\draw [thick,black] plot  coordinates {(-3/32,2/32) (3/32,-2/32) };
\end{scope}
\draw [thick,cborange] plot  coordinates {(-1/2,0) (5/2,0) };
\fill [black] (0,0) circle (0.06);
\draw (0-1/16,0) node [above,red] { $\frac{\alpha}{8}$};
\begin{scope}[shift={(4/2,0)}]
\draw [thick,black] plot  coordinates {(8/12,-8/8) (-1/6,1/4) };
\fill [black] (0,0) circle (0.06);
\draw (0+1/16,0) node [above,red] { $\frac{3\alpha}{8}$};
\begin{scope}[shift={(2/9,-3/9)}]
\draw [thick,black] plot  coordinates {(3/32,2/32) (-3/32,-2/32) };
\end{scope}
\begin{scope}[shift={(4/9,-6/9)}]
\draw [thick,black] plot  coordinates {(3/32,2/32) (-3/32,-2/32) };
\end{scope}
\end{scope}
\draw (1,0) node [below] {\hbox{\footnotesize$\geqtwo$}};
\end{tikzpicture}}

\def\includefigB#1{
\makebox(130pt,75pt){#1}}

\def\FigTwoCaseB{
\begin{tabular}{|c|c|c|}
\hline
{\large\strut} (B1) & (B2) & (B3) \\
\hline
${\Large\strut}
\delta=0\ \land \ \alpha<4$ 
& $\delta=0\ \land \ \alpha=4$ 
& $\delta=0\ \land \ \alpha>4$ \\
\hline
\includefigB{\FigTwoCaseBi} & 
\includefigB{\FigTwoCaseBii} & 
\includefigB{\FigTwoCaseBiii} \\
\hline\hline
{\large\strut} (B4) & (B5) & (B6) \\
\hline
${\Large\strut}
0<2\delta<\alpha \ \land \ \alpha+\delta<4$ 
& $0<2\delta<\alpha\ \land \ \alpha+\delta=4$ 
& $0<\delta<\frac{4}{3}\ \land \ \alpha+\delta>4$ \\
\hline
\includefigB{\FigTwoCaseBiv} &
\includefigB{\FigTwoCaseBv} & 
\includefigB{\FigTwoCaseBvi} \\
\hline\hline
{\large\strut} (B7)& (B8)& (B9) \\ 
\hline
{\Large\strut} 
  $\delta\ge\frac{4}{3}\ \land \ \alpha>\frac{8}{3}$ 
&  $\alpha=2\delta<\frac{8}{3}$
& $\alpha<2\delta<\frac{8+2\alpha}{5}$ 
\\
\hline
\includefigB{\FigTwoCaseBvii} &
\includefigB{\FigTwoCaseBviii} &
\includefigB{\FigTwoCaseBix} 
 \\ 
\hline 
 \noalign{\global\arrayrulewidth=2.0pt}
\arrayrulecolor{white}\hline
 \noalign{\global\arrayrulewidth=0.4pt}
 \arrayrulecolor{black} \cline{1-2}
{\large\strut} (B10) & (B11) & \multicolumn{1}{c}{}  \\
\cline{1-2}
{\Large\strut}
$\alpha<\frac{8+2\alpha}{5}\leq 2 \delta$ 
& $\alpha=\frac{8}{3} \leq 2 \delta$ 
& \multicolumn{1}{c}{}  \\
\cline{1-2}
\includefigB{\FigTwoCaseBx} &
\includefigB{\FigTwoCaseBxi} & 
\multicolumn{1}{c}{} \\
\cline{1-2}

\end{tabular}}

\def\FigTwoCaseCi{\begin{tikzpicture}[baseline=13pt]
		\begin{scope}[shift={(3/4+1/2,1/2)}]
			\draw [thick, cborange] plot coordinates {(-3/4+0.14/3-3/8,0.14/2) (3/4-0.14/3+3/8,0.14/2) };
			\draw (0,-7/8) node [below] {\hbox{\footnotesize$\geqone$}};
			\draw [thick, cbblue] plot coordinates {(0,0.14/2+3/8) (0,0.14/2-8/8) } ;
			\fill [black] (-3/4+0.14/3,0.14/2) circle (0.06);
			\draw (-3/4+0.14/3-1/16,0.14/2) node [above,red] { $\frac{\alpha}{4}$};
			\fill [black] (3/4-0.14/3,0.14/2) circle (0.06);
						\draw(3/4-0.14/3+1/6,0.14/2) node [above,red] { $\frac{\alpha}{4}$};
			\fill [black] (0,0.14/2) circle (0.06);
						\draw(0-4/16,0.14/2) node [below,red] { $\frac{4\text{-}\alpha}{12}$};
			
		\end{scope}
		\begin{scope}[shift={(-2/2,-1)}, yscale=-1]
			\draw [thick, cborange] plot coordinates {(-3/4+0.14/3-3/8,0.14/2) (3/4-0.14/3+3/8,0.14/2) };
			\draw (0,-7/8) node [above] {\hbox{\footnotesize$\geqone$}};
			\draw [thick, cbblue] plot coordinates {(0,0.14/2+3/8) (0,0.14/2-8/8) } ;
			\fill [black] (-3/4+0.14/3,0.14/2) circle (0.06);
			\draw (-3/4+0.14/3-2/16,0.14/2) node [below,red] { $\frac{\beta}{4}$};
			\fill [black] (3/4-0.14/3,0.14/2) circle (0.06);
			\draw(3/4-0.14/3+1/6,0.14/2) node [below,red] { $\frac{\beta}{4}$};
			\fill [black] (0,0.14/2) circle (0.06);
			\draw(0-4/16,0.14/2) node [above,red] { $\frac{4\text{-}\beta}{12}$};
		\end{scope}
		
		\draw [thick, black, shorten >= -0.25cm, shorten <=-0.25cm] (0.14/3+1/2,0.14/2+1/2)--(-0.14/3 -1/4,-0.14/2-1);	
		\begin{scope} [shift={(-1/8+1/4+1/4*253/492,-1/4+1/4 )}]
			\draw [thick,black] plot  coordinates {(492/253*0.061,-1*0.061) (-492/253*0.061,1*0.061) };
	\end{scope} 
	
		\begin{scope} [shift={(-1/8+1/4-1/4*253/492,-1/4-1/4 )}]
			\draw [thick,black] plot  coordinates {(492/253*0.061,-1*0.061) (-492/253*0.061,1*0.061) };
	\end{scope}

		\begin{scope}[shift={(3/2+1/2,1/2)}]
			\begin{scope}[shift={(-1/48,0)}]
				\draw [thick,black] plot  coordinates {(4/12,-8/8) (-1/12,1/4) };
				\begin{scope}[shift={(1/9,-3/9)}]
					\draw [thick,black] plot  coordinates {(8/64,3/64) (-8/64,-3/64) };
				\end{scope}		
				\begin{scope}[shift={(2/9,-6/9)}]
					\draw [thick,black] plot  coordinates {(8/64,3/64) (-8/64,-3/64) };
				\end{scope}
			\end{scope}
		\end{scope}
		\begin{scope}[shift={(-3/2-1/4,-1)}, yscale=-1, xscale=-1]
			\begin{scope}[shift={(-1/48,0)}]
				\draw [thick,black] plot  coordinates {(4/12,-8/8) (-1/12,1/4) };
				\begin{scope}[shift={(1/9,-3/9)}]
					\draw [thick,black] plot  coordinates {(8/64,3/64) (-8/64,-3/64) };
				\end{scope}		
				\begin{scope}[shift={(2/9,-6/9)}]
					\draw [thick,black] plot  coordinates {(8/64,3/64) (-8/64,-3/64) };
				\end{scope}
			\end{scope}
		\end{scope}
\end{tikzpicture}}

\def\FigTwoCaseCii{\begin{tikzpicture}[baseline=13pt]
		\begin{scope}[shift={(3/4+1/2,0+1/2)}]
			\draw [thick, cborange] plot coordinates {(-3/4+0.14/3-3/8,0.14/2) (3/4-0.14/3+3/8,0.14/2) };
			\draw (0,0) node [below] {\hbox{\footnotesize$\geqone$}};
			\fill [black] (-3/4+0.14/3,0.14/2) circle (0.06);
			\draw (-3/4+0.14/3-2/16,0.14/2) node [above,red] { \scriptsize$1$};
			\fill [black] (3/4-0.14/3,0.14/2) circle (0.06);
				\draw(3/4-0.14/3+1/6,0.14/2) node [above,red] { \scriptsize$1$};
		\end{scope}
			\begin{scope}[shift={(-2/2,-1)}, yscale=-1]
			\draw [thick, cborange] plot coordinates {(-3/4+0.14/3-3/8,0.14/2) (3/4-0.14/3+3/8,0.14/2) };
			\draw (0,-7/8) node [above] {\hbox{\footnotesize$\geqone$}};
			\draw [thick, cbblue] plot coordinates {(0,0.14/2+3/8) (0,0.14/2-8/8) } ;
			\fill [black] (-3/4+0.14/3,0.14/2) circle (0.06);
			\draw (-3/4+0.14/3-2/16,0.14/2) node [below,red] { $\frac{\beta}{4}$};
			\fill [black] (3/4-0.14/3,0.14/2) circle (0.06);
			\draw(3/4-0.14/3+1/6,0.14/2) node [below,red] { $\frac{\beta}{4}$};
			\fill [black] (0,0.14/2) circle (0.06);
			\draw(0-4/16,0.14/2) node [above,red] { $\frac{4\text{-}\beta}{12}$};
		\end{scope}
	\draw [thick, black, shorten >= -0.25cm, shorten <=-0.25cm] (0.14/3+1/2,0.14/2+1/2)--(-0.14/3 -1/4,-0.14/2-1);	
	\begin{scope} [shift={(-1/8+1/4+1/4*253/492,-1/4+1/4 )}]
		\draw [thick,black] plot  coordinates {(492/253*0.061,-1*0.061) (-492/253*0.061,1*0.061) };
	\end{scope} 
	
	\begin{scope} [shift={(-1/8+1/4-1/4*253/492,-1/4-1/4 )}]
		\draw [thick,black] plot  coordinates {(492/253*0.061,-1*0.061) (-492/253*0.061,1*0.061) };
	\end{scope}

		\begin{scope}[shift={(3/2+1/2,1/2)}]
			\begin{scope}[shift={(-1/48,0)}]
				\draw [thick,black] plot  coordinates {(4/12,-8/8) (-1/12,1/4) };
				\begin{scope}[shift={(1/9,-3/9)}]
					\draw [thick,black] plot  coordinates {(8/64,3/64) (-8/64,-3/64) };
				\end{scope}		
				\begin{scope}[shift={(2/9,-6/9)}]
					\draw [thick,black] plot  coordinates {(8/64,3/64) (-8/64,-3/64) };
				\end{scope}
			\end{scope}
		\end{scope}
		\begin{scope}[shift={(-3/2-1/4,-1)}, yscale=-1, xscale=-1]
			\begin{scope}[shift={(-1/48,0)}]
				\draw [thick,black] plot  coordinates {(4/12,-8/8) (-1/12,1/4) };
				\begin{scope}[shift={(1/9,-3/9)}]
					\draw [thick,black] plot  coordinates {(8/64,3/64) (-8/64,-3/64) };
				\end{scope}		
				\begin{scope}[shift={(2/9,-6/9)}]
					\draw [thick,black] plot  coordinates {(8/64,3/64) (-8/64,-3/64) };
				\end{scope}
			\end{scope}
		\end{scope}
\end{tikzpicture}}

\def\FigTwoCaseCiii{\begin{tikzpicture}[baseline=13pt]
		\begin{scope}[shift={(3/4+1/2,0+1/2)}]
			\draw [thick, cborange] plot coordinates {(-3/4+0.14/3-3/8,0.14/2) (3/4-0.14/3+3/8,0.14/2) };
			\draw (0,0) node [below] {\hbox{\footnotesize$\geqone$}};
			\fill [black] (-3/4+0.14/3,0.14/2) circle (0.06);
			\draw (-3/4+0.14/3-2/16,0.14/2) node [above,red] { \scriptsize$1$};
			\fill [black] (3/4-0.14/3,0.14/2) circle (0.06);
			\draw(3/4-0.14/3+1/6,0.14/2) node [above,red] { \scriptsize$1$};
		\end{scope}
		\begin{scope}[shift={(-2/2,-1)}, yscale=-1]
		\draw [thick, cborange] plot coordinates {(-3/4+0.14/3-3/8,0.14/2) (3/4-0.14/3+3/8,0.14/2) };
		\draw (0,0) node [above] {\hbox{\footnotesize$\geqone$}};
		\fill [black] (-3/4+0.14/3,0.14/2) circle (0.06);
		\draw (-3/4+0.14/3-2/16,0.14/2) node [below,red] { \scriptsize$1$};
		\fill [black] (3/4-0.14/3,0.14/2) circle (0.06);
		\draw(3/4-0.14/3+1/6,0.14/2) node [below,red] { \scriptsize$1$};
		\end{scope}
	\draw [thick, black, shorten >= -0.25cm, shorten <=-0.25cm] (0.14/3+1/2,0.14/2+1/2)--(-0.14/3 -1/4,-0.14/2-1);	
	\begin{scope} [shift={(-1/8+1/4+1/4*253/492,-1/4+1/4 )}]
		\draw [thick,black] plot  coordinates {(492/253*0.061,-1*0.061) (-492/253*0.061,1*0.061) };
	\end{scope} 
	
	\begin{scope} [shift={(-1/8+1/4-1/4*253/492,-1/4-1/4 )}]
		\draw [thick,black] plot  coordinates {(492/253*0.061,-1*0.061) (-492/253*0.061,1*0.061) };
	\end{scope} 
		
		\begin{scope}[shift={(3/2+1/2,0+1/2)}]
			\begin{scope}[shift={(-1/48,0)}]
				\draw [thick,black] plot  coordinates {(4/12,-8/8) (-1/12,1/4) };
				\begin{scope}[shift={(1/9,-3/9)}]
					\draw [thick,black] plot  coordinates {(8/64,3/64) (-8/64,-3/64) };
				\end{scope}		
				\begin{scope}[shift={(2/9,-6/9)}]
					\draw [thick,black] plot  coordinates {(8/64,3/64) (-8/64,-3/64) };
				\end{scope}
			\end{scope}
		\end{scope}
		\begin{scope}[shift={(-3/2-1/4,-1)}, yscale=-1, xscale=-1]
			\begin{scope}[shift={(-1/48,0)}]
				\draw [thick,black] plot  coordinates {(4/12,-8/8) (-1/12,1/4) };
				\begin{scope}[shift={(1/9,-3/9)}]
					\draw [thick,black] plot  coordinates {(8/64,3/64) (-8/64,-3/64) };
				\end{scope}		
				\begin{scope}[shift={(2/9,-6/9)}]
					\draw [thick,black] plot  coordinates {(8/64,3/64) (-8/64,-3/64) };
				\end{scope}
			\end{scope}
		\end{scope}
\end{tikzpicture}}

\def\FigTwoCaseCiv{\begin{tikzpicture}[baseline=13pt]
		\begin{scope}[shift={(3/4+1/2,0+1/2)}]
			\draw [thick,cborange] plot [smooth, tension=1] coordinates {(-3/3,2/12) (1/4,-6/12) (-1/3,-14/12) };
			\draw [thick,cborange] plot [smooth, tension=1] coordinates {(3/3,2/12) (-1/4,-6/12) (1/3,-14/12) };
			\fill [black] (0,-0.26) circle (0.06);
				\draw (0,-0.26+1/16) node [above,red] { \scriptsize$\alpha \text{-}4$};
			\fill [black] (0,-0.99) circle (0.06);
						\draw (0,-0.99-1/16) node [below,red] { \scriptsize$\alpha \text{-}4$};
			\fill [black] (-3/4+0.14/3,0.14/2) circle (0.06);
				\draw (-3/4+0.14/3-2/16,0.14/2) node [above,red] { \scriptsize$1$};
			\fill [black] (3/4-0.14/3,0.14/2) circle (0.06);
						\draw(3/4-0.14/3+1/6,0.14/2) node [above,red] { \scriptsize$1$};
		\end{scope}
	\begin{scope}[shift={(-2/2,-1)}, yscale=-1]
	\draw [thick, cborange] plot coordinates {(-3/4+0.14/3-3/8,0.14/2) (3/4-0.14/3+3/8,0.14/2) };
	\draw (0,-7/8) node [above] {\hbox{\footnotesize$\geqone$}};
	\draw [thick, cbblue] plot coordinates {(0,0.14/2+3/8) (0,0.14/2-8/8) } ;
	\fill [black] (-3/4+0.14/3,0.14/2) circle (0.06);
	\draw (-3/4+0.14/3-2/16,0.14/2) node [below,red] { $\frac{\beta}{4}$};
	\fill [black] (3/4-0.14/3,0.14/2) circle (0.06);
	\draw(3/4-0.14/3+1/6,0.14/2) node [below,red] { $\frac{\beta}{4}$};
	\fill [black] (0,0.14/2) circle (0.06);
	\draw(0-4/16,0.14/2) node [above,red] { $\frac{4\text{-}\beta}{12}$};
\end{scope}
	\draw [thick, black, shorten >= -0.25cm, shorten <=-0.25cm] (0.14/3+1/2,0.14/2+1/2)--(-0.14/3 -1/4,-0.14/2-1);	
	\begin{scope} [shift={(-1/8+1/4+1/4*253/492,-1/4+1/4 )}]
		\draw [thick,black] plot  coordinates {(492/253*0.061,-1*0.061) (-492/253*0.061,1*0.061) };
	\end{scope} 
	
	\begin{scope} [shift={(-1/8+1/4-1/4*253/492,-1/4-1/4 )}]
		\draw [thick,black] plot  coordinates {(492/253*0.061,-1*0.061) (-492/253*0.061,1*0.061) };
	\end{scope} 
		\begin{scope}[shift={(3/2+1/2,0+1/2)}]
			\begin{scope}[shift={(-1/48,0)}]
				\draw [thick,black] plot  coordinates {(4/12,-8/8) (-1/12,1/4) };
				\begin{scope}[shift={(1/9,-3/9)}]
					\draw [thick,black] plot  coordinates {(8/64,3/64) (-8/64,-3/64) };
				\end{scope}		
				\begin{scope}[shift={(2/9,-6/9)}]
					\draw [thick,black] plot  coordinates {(8/64,3/64) (-8/64,-3/64) };
				\end{scope}
			\end{scope}
		\end{scope}
		\begin{scope}[shift={(-3/2-1/4,-1)}, yscale=-1, xscale=-1]
			\begin{scope}[shift={(-1/48,0)}]
				\draw [thick,black] plot  coordinates {(4/12,-8/8) (-1/12,1/4) };
				\begin{scope}[shift={(1/9,-3/9)}]
					\draw [thick,black] plot  coordinates {(8/64,3/64) (-8/64,-3/64) };
				\end{scope}		
				\begin{scope}[shift={(2/9,-6/9)}]
					\draw [thick,black] plot  coordinates {(8/64,3/64) (-8/64,-3/64) };
				\end{scope}
			\end{scope}
		\end{scope}
\end{tikzpicture}}

\def\FigTwoCaseCv{\begin{tikzpicture}[baseline=13pt]
		\begin{scope}[shift={(3/4+1/2,0+1/2)}]
				\draw [thick,cborange] plot [smooth, tension=1] coordinates {(-3/3,2/12) (1/4,-6/12) (-1/3,-14/12) };
			\draw [thick,cborange] plot [smooth, tension=1] coordinates {(3/3,2/12) (-1/4,-6/12) (1/3,-14/12) };
			\fill [black] (0,-0.26) circle (0.06);
			\draw (0,-0.26+1/16) node [above,red] { \scriptsize$\alpha \text{-}4$};
			\fill [black] (0,-0.99) circle (0.06);
			\draw (0,-0.99-1/16) node [below,red] { \scriptsize$\alpha \text{-}4$};
			\fill [black] (-3/4+0.14/3,0.14/2) circle (0.06);
			\draw (-3/4+0.14/3-2/16,0.14/2) node [above,red] { \scriptsize$1$};
			\fill [black] (3/4-0.14/3,0.14/2) circle (0.06);
			\draw(3/4-0.14/3+1/6,0.14/2) node [above,red] { \scriptsize$1$};
		\end{scope}
			\begin{scope}[shift={(-2/2,-1)}, yscale=-1]
			\draw [thick, cborange] plot coordinates {(-3/4+0.14/3-3/8,0.14/2) (3/4-0.14/3+3/8,0.14/2) };
			\draw (0,0) node [above] {\hbox{\footnotesize$\geqone$}};
			\fill [black] (-3/4+0.14/3,0.14/2) circle (0.06);
			\draw (-3/4+0.14/3-2/16,0.14/2) node [below,red] { \scriptsize$1$};
			\fill [black] (3/4-0.14/3,0.14/2) circle (0.06);
			\draw(3/4-0.14/3+1/6,0.14/2) node [below,red] { \scriptsize$1$};
		\end{scope}
		\draw [thick, black, shorten >= -0.25cm, shorten <=-0.25cm] (0.14/3+1/2,0.14/2+1/2)--(-0.14/3 -1/4,-0.14/2-1);	
		\begin{scope} [shift={(-1/8+1/4+1/4*253/492,-1/4+1/4 )}]
			\draw [thick,black] plot  coordinates {(492/253*0.061,-1*0.061) (-492/253*0.061,1*0.061) };
		\end{scope} 
		
		\begin{scope} [shift={(-1/8+1/4-1/4*253/492,-1/4-1/4 )}]
			\draw [thick,black] plot  coordinates {(492/253*0.061,-1*0.061) (-492/253*0.061,1*0.061) };
		\end{scope} 
	
		\begin{scope}[shift={(3/2+1/2,0+1/2)}]
			\begin{scope}[shift={(-1/48,0)}]
				\draw [thick,black] plot  coordinates {(4/12,-8/8) (-1/12,1/4) };
				\begin{scope}[shift={(1/9,-3/9)}]
					\draw [thick,black] plot  coordinates {(8/64,3/64) (-8/64,-3/64) };
				\end{scope}		
				\begin{scope}[shift={(2/9,-6/9)}]
					\draw [thick,black] plot  coordinates {(8/64,3/64) (-8/64,-3/64) };
				\end{scope}
			\end{scope}
		\end{scope}
		\begin{scope}[shift={(-3/2-1/4,-1)}, yscale=-1, xscale=-1]
			\begin{scope}[shift={(-1/48,0)}]
				\draw [thick,black] plot  coordinates {(4/12,-8/8) (-1/12,1/4) };
				\begin{scope}[shift={(1/9,-3/9)}]
					\draw [thick,black] plot  coordinates {(8/64,3/64) (-8/64,-3/64) };
				\end{scope}		
				\begin{scope}[shift={(2/9,-6/9)}]
					\draw [thick,black] plot  coordinates {(8/64,3/64) (-8/64,-3/64) };
				\end{scope}
			\end{scope}
		\end{scope}
\end{tikzpicture}}

\def\FigTwoCaseCvi{\begin{tikzpicture}[baseline=13pt]
		\begin{scope}[shift={(3/4+1/2,0+1/2)}]
				\draw [thick,cborange] plot [smooth, tension=1] coordinates {(-3/3,2/12) (1/4,-6/12) (-1/3,-14/12) };
			\draw [thick,cborange] plot [smooth, tension=1] coordinates {(3/3,2/12) (-1/4,-6/12) (1/3,-14/12) };
			\fill [black] (0,-0.26) circle (0.06);
			\draw (0,-0.26+1/16) node [above,red] { \scriptsize$\alpha \text{-}4$};
			\fill [black] (0,-0.99) circle (0.06);
			\draw (0,-0.99-1/16) node [below,red] { \scriptsize$\alpha \text{-}4$};
			\fill [black] (-3/4+0.14/3,0.14/2) circle (0.06);
			\draw (-3/4+0.14/3-2/16,0.14/2) node [above,red] { \scriptsize$1$};
			\fill [black] (3/4-0.14/3,0.14/2) circle (0.06);
			\draw(3/4-0.14/3+1/6,0.14/2) node [above,red] { \scriptsize$1$};
		\end{scope}
		\begin{scope}[shift={(-2/2,-1)}, yscale=-1]
			\draw [thick,cborange] plot [smooth, tension=1] coordinates {(-3/3,2/12) (1/4,-6/12) (-1/3,-14/12) };
			\draw [thick,cborange] plot [smooth, tension=1] coordinates {(3/3,2/12) (-1/4,-6/12) (1/3,-14/12) };
			\fill [black] (0,-0.26) circle (0.06);
						\draw (0,-0.26+1/16) node [below,red] { \scriptsize$\beta \text{-}4$};
			\fill [black] (0,-0.99) circle (0.06);
				\draw (0,-0.99-1/16) node [above,red] { \scriptsize$\beta \text{-}4$};
			\fill [black] (-3/4+0.14/3,0.14/2) circle (0.06);
						\draw (-3/4+0.14/3-2/16,0.14/2) node [below,red] { \scriptsize$1$};
			\fill [black] (3/4-0.14/3,0.14/2) circle (0.06);
						\draw(3/4-0.14/3+1/6,0.14/2) node [below,red] { \scriptsize$1$};
		\end{scope}
	\draw [thick, black, shorten >= -0.25cm, shorten <=-0.25cm] (0.14/3+1/2,0.14/2+1/2)--(-0.14/3 -1/4,-0.14/2-1);	
	\begin{scope} [shift={(-1/8+1/4+1/4*253/492,-1/4+1/4 )}]
		\draw [thick,black] plot  coordinates {(492/253*0.061,-1*0.061) (-492/253*0.061,1*0.061) };
	\end{scope} 
	
	\begin{scope} [shift={(-1/8+1/4-1/4*253/492,-1/4-1/4 )}]
		\draw [thick,black] plot  coordinates {(492/253*0.061,-1*0.061) (-492/253*0.061,1*0.061) };
	\end{scope} 
		
		\begin{scope}[shift={(3/2+1/2,0+1/2)}]
			\begin{scope}[shift={(-1/48,0)}]
				\draw [thick,black] plot  coordinates {(4/12,-8/8) (-1/12,1/4) };
				\begin{scope}[shift={(1/9,-3/9)}]
					\draw [thick,black] plot  coordinates {(8/64,3/64) (-8/64,-3/64) };
				\end{scope}		
				\begin{scope}[shift={(2/9,-6/9)}]
					\draw [thick,black] plot  coordinates {(8/64,3/64) (-8/64,-3/64) };
				\end{scope}
			\end{scope}
		\end{scope}
		\begin{scope}[shift={(-3/2-1/4,-1)}, yscale=-1, xscale=-1]
			\begin{scope}[shift={(-1/48,0)}]
				\draw [thick,black] plot  coordinates {(4/12,-8/8) (-1/12,1/4) };
				\begin{scope}[shift={(1/9,-3/9)}]
					\draw [thick,black] plot  coordinates {(8/64,3/64) (-8/64,-3/64) };
				\end{scope}		
				\begin{scope}[shift={(2/9,-6/9)}]
					\draw [thick,black] plot  coordinates {(8/64,3/64) (-8/64,-3/64) };
				\end{scope}
			\end{scope}
		\end{scope}
\end{tikzpicture}}

\def\includefigC#1{
\makebox(130pt,90pt){#1}}

\def\FigTwoCaseC{
\begin{tabular}{|c|c|c|}
	\hline
{\large\strut}	(C1) & (C2) & (C3) \\
\hline
{\Large\strut} $4>\alpha\ge\beta$ & $4=\alpha>\beta$ & $\alpha=4=\beta$ \\
\hline
\includefigC{\FigTwoCaseCi} & 
\includefigC{\FigTwoCaseCii} & 
\includefigC{\FigTwoCaseCiii} \\
\hline\hline
{\large\strut}	(C4) & (C5) & (C6) \\
\hline
{\Large\strut} $\alpha>4>\beta$ & $\alpha>4=\beta$ & $\alpha\ge\beta>4$ \\
\hline
\includefigC{\FigTwoCaseCiv} & 
\includefigC{\FigTwoCaseCv} & 
\includefigC{\FigTwoCaseCvi} \\
\hline
\end{tabular}}

\def\FigTwoCaseDi{\begin{tikzpicture}[baseline=13pt]
\draw [thick,black] plot  coordinates {(-1/2,0) (4/2,0) };
\begin{scope}[shift={(3/2,0)}]
	\draw [thick,black] plot  coordinates {(-4/6,-4/4) (1/6,1/4) };
	\fill [black] (0,0) circle (0.06);
	\draw (0,0+1/16) node [above,red] { \scriptsize$1$};
	\begin{scope}[shift={(-2/9,-3/9)}]
		\draw [thick,black] plot  coordinates {(-3/32,2/32) (3/32,-2/32) };
	\end{scope}
	\begin{scope}[shift={(-4/9,-6/9)}]
		\draw [thick,black] plot  coordinates {(-3/32,2/32) (3/32,-2/32) };
	\end{scope}
\end{scope}
\begin{scope}[shift={(3/4,0)}]
	\draw [thick,black] plot  coordinates {(-4/6,-4/4) (1/6,1/4) };
	\fill [black] (0,0) circle (0.06);
		\draw (0,0+1/16) node [above,red] { \scriptsize$1$};
	\begin{scope}[shift={(-2/9,-3/9)}]
		\draw [thick,black] plot  coordinates {(-3/32,2/32) (3/32,-2/32) };
	\end{scope}
	\begin{scope}[shift={(-4/9,-6/9)}]
		\draw [thick,black] plot  coordinates {(-3/32,2/32) (3/32,-2/32) };
	\end{scope}
\end{scope}
\begin{scope}[shift={(0,0)}]
	\draw [thick,black] plot  coordinates {(-4/6,-4/4) (1/6,1/4) };
	\fill [black] (0,0) circle (0.06);
		\draw (0,0+1/16) node [above,red] { \scriptsize$1$};
	\begin{scope}[shift={(-2/9,-3/9)}]
		\draw [thick,black] plot  coordinates {(-3/32,2/32) (3/32,-2/32) };
	\end{scope}
	\begin{scope}[shift={(-4/9,-6/9)}]
		\draw [thick,black] plot  coordinates {(-3/32,2/32) (3/32,-2/32) };
	\end{scope}
\end{scope}
\draw (-2/4,0) node [left] {\hbox{\footnotesize$\geqtwo$}};
\end{tikzpicture}}

\def\FigTwoCaseDii{\begin{tikzpicture}[baseline=13pt]
		\draw [thick,black] plot  coordinates {(-1/2,0) (3/2+1/8+1/100,0) };
		\draw [thick, black] (3/2+1/8,0) arc (270:450:0.25);
		\draw [thick, black] plot coordinates {(3/2+1/8+1/100,0.5) (2/2+1/8,0.5)};
		\begin{scope}[shift={(3/2,0)}]
			\draw [thick,cborange] plot  coordinates {(0,6/8) (0,0) };
			\draw [thick,cborange] plot [smooth, tension=1] coordinates {(0,0) (1/16,-4/8) (3/16+1/32,-7/8)  };
			\fill [black] (0,1/2) circle (0.06);
				\draw (0,1/2+3/16) node [left,red] { \scriptsize$\alpha\text{-}2$};
			\fill [black] (0,0) circle (0.06);
				\draw (0,0-3/16) node [left,red] { \scriptsize$\alpha\text{-}2$};
		\end{scope}
		\begin{scope}[shift={(3/2+1/16,-1/2)}]
			\draw [thick,black] plot  coordinates {(-4/6,-4/4) (1/6,1/4) };
			\fill [black] (0,0) circle (0.06);
				\draw (0,0+0/16) node [right,red] { \scriptsize$1$};
			\begin{scope}[shift={(-2/9,-3/9)}]
				\draw [thick,black] plot  coordinates {(-3/32,2/32) (3/32,-2/32) };
			\end{scope}
			\begin{scope}[shift={(-4/9,-6/9)}]
				\draw [thick,black] plot  coordinates {(-3/32,2/32) (3/32,-2/32) };
			\end{scope}
		\end{scope}
		\begin{scope}[shift={(3/4,0)}]
			\draw [thick,black] plot  coordinates {(-4/6,-4/4) (1/6,1/4) };
			\fill [black] (0,0) circle (0.06);
				\draw (0,0+1/16) node [above,red] { \scriptsize$1$};
			\begin{scope}[shift={(-2/9,-3/9)}]
				\draw [thick,black] plot  coordinates {(-3/32,2/32) (3/32,-2/32) };
			\end{scope}
			\begin{scope}[shift={(-4/9,-6/9)}]
				\draw [thick,black] plot  coordinates {(-3/32,2/32) (3/32,-2/32) };
			\end{scope}
		\end{scope}
		\begin{scope}[shift={(0,0)}]
			\draw [thick,black] plot  coordinates {(-4/6,-4/4) (1/6,1/4) };
			\fill [black] (0,0) circle (0.06);
				\draw (0,0+1/16) node [above,red] { \scriptsize$1$};
			\begin{scope}[shift={(-2/9,-3/9)}]
				\draw [thick,black] plot  coordinates {(-3/32,2/32) (3/32,-2/32) };
			\end{scope}
			\begin{scope}[shift={(-4/9,-6/9)}]
				\draw [thick,black] plot  coordinates {(-3/32,2/32) (3/32,-2/32) };
			\end{scope}
		\end{scope}
		\draw (-2/4,0) node [left] {\hbox{\footnotesize$\geqone$}};
\end{tikzpicture}}

\def\FigTwoCaseDiii{\begin{tikzpicture}[baseline=13pt]
			\draw [thick,black] plot  coordinates {(-1/2,0) (3/2+1/8+1/100,0) };
		\draw [thick, black] (3/2+1/8,0) arc (270:450:0.25);
		\draw [thick, black] plot coordinates {(3/2+1/8+1/100,0.5) (-1/2,0.5)};
		\begin{scope}[shift={(3/2,0)}]
\draw [thick,cborange] plot  coordinates {(0,6/8) (0,0) };
\draw [thick,cborange] plot [smooth, tension=1] coordinates {(0,0) (1/16,-4/8) (3/16+1/32,-7/8)  };
\fill [black] (0,1/2) circle (0.06);
\draw (0,1/2+3/16) node [left,red] { \scriptsize$\alpha\text{-}2$};
\fill [black] (0,0) circle (0.06);
\draw (0,0+3/16) node [left,red] { \scriptsize$\alpha\text{-}2$};
		\end{scope}
	\begin{scope}[shift={(3/2-6/4+3/4,0)}]
\draw [thick,cborange] plot  coordinates {(0,6/8) (0,0) };
\draw [thick,cborange] plot [smooth, tension=1] coordinates {(0,0) (1/16,-4/8) (3/16+1/32,-7/8)  };
\fill [black] (0,1/2) circle (0.06);
\draw (0,1/2+3/16) node [left,red] { \scriptsize$\beta\text{-}2$};
\fill [black] (0,0) circle (0.06);
\draw (0,0-3/16) node [left,red] { \scriptsize$\beta\text{-}2$};
	\end{scope}
		\begin{scope}[shift={(3/2+1/16,-1/2)}]
			\draw [thick,black] plot  coordinates {(-4/6,-4/4) (1/6,1/4) };
			\fill [black] (0,0) circle (0.06);
						\draw (0,0+0/16) node [right,red] { \scriptsize$1$};
			\begin{scope}[shift={(-2/9,-3/9)}]
				\draw [thick,black] plot  coordinates {(-3/32,2/32) (3/32,-2/32) };
			\end{scope}
			\begin{scope}[shift={(-4/9,-6/9)}]
				\draw [thick,black] plot  coordinates {(-3/32,2/32) (3/32,-2/32) };
			\end{scope}
		\end{scope}
	\begin{scope}[shift={(3/4+1/16,-1/2)}]
		\draw [thick,black] plot  coordinates {(-4/6,-4/4) (1/6,1/4) };
		\fill [black] (0,0) circle (0.06);
						\draw (0,0+0/16) node [right,red] { \scriptsize$1$};
		\begin{scope}[shift={(-2/9,-3/9)}]
			\draw [thick,black] plot  coordinates {(-3/32,2/32) (3/32,-2/32) };
		\end{scope}
		\begin{scope}[shift={(-4/9,-6/9)}]
			\draw [thick,black] plot  coordinates {(-3/32,2/32) (3/32,-2/32) };
		\end{scope}
	\end{scope}
		\begin{scope}[shift={(0+0/16,0)}]
			\draw [thick,black] plot  coordinates {(-4/6,-4/4) (1/6,1/4) };
			\fill [black] (0,0) circle (0.06);
							\draw (0,0+1/16) node [above,red] { \scriptsize$1$};
			\begin{scope}[shift={(-2/9,-3/9)}]
				\draw [thick,black] plot  coordinates {(-3/32,2/32) (3/32,-2/32) };
			\end{scope}
			\begin{scope}[shift={(-4/9,-6/9)}]
				\draw [thick,black] plot  coordinates {(-3/32,2/32) (3/32,-2/32) };
			\end{scope}
		\end{scope}
\end{tikzpicture}}

\def\FigTwoCaseDiv{\begin{tikzpicture}[baseline=13pt]
		\draw [thick,black] plot  coordinates {(-1/2,0) (3/2+3/8,0) };
	\draw [thick, black] plot coordinates {(3/2+3/8,0.5) (-1/2,0.5)};
	\begin{scope}[shift={(3/2,0)}]
\draw [thick,cborange] plot  coordinates {(0,6/8) (0,0) };
\draw [thick,cborange] plot [smooth, tension=1] coordinates {(0,0) (1/16,-4/8) (3/16+1/32,-7/8)  };
\fill [black] (0,1/2) circle (0.06);
\draw (0,1/2+3/16) node [left,red] { \scriptsize$\alpha\text{-}2$};
\fill [black] (0,0) circle (0.06);
\draw (0,0+3/16) node [left,red] { \scriptsize$\alpha\text{-}2$};
	\end{scope}
	\begin{scope}[shift={(3/2-3/4,0)}]
\draw [thick,cborange] plot  coordinates {(0,6/8) (0,0) };
\draw [thick,cborange] plot [smooth, tension=1] coordinates {(0,0) (1/16,-4/8) (3/16+1/32,-7/8)  };
\fill [black] (0,1/2) circle (0.06);
\draw (0,1/2+3/16) node [left,red] { \scriptsize$\beta\text{-}2$};
\fill [black] (0,0) circle (0.06);
\draw (0,0+3/16) node [left,red] { \scriptsize$\beta\text{-}2$};
	\end{scope}
	\begin{scope}[shift={(3/2+1/16,-1/2)}]
		\draw [thick,black] plot  coordinates {(-4/6,-4/4) (1/6,1/4) };
		\fill [black] (0,0) circle (0.06);
								\draw (0,0+0/16) node [right,red] { \scriptsize$1$};
		\begin{scope}[shift={(-2/9,-3/9)}]
			\draw [thick,black] plot  coordinates {(-3/32,2/32) (3/32,-2/32) };
		\end{scope}
		\begin{scope}[shift={(-4/9,-6/9)}]
			\draw [thick,black] plot  coordinates {(-3/32,2/32) (3/32,-2/32) };
		\end{scope}
	\end{scope}
	\begin{scope}[shift={(3/4+1/16,-1/2)}]
		\draw [thick,black] plot  coordinates {(-4/6,-4/4) (1/6,1/4) };
		\fill [black] (0,0) circle (0.06);
								\draw (0,0+0/16) node [right,red] { \scriptsize$1$};
		\begin{scope}[shift={(-2/9,-3/9)}]
			\draw [thick,black] plot  coordinates {(-3/32,2/32) (3/32,-2/32) };
		\end{scope}
		\begin{scope}[shift={(-4/9,-6/9)}]
			\draw [thick,black] plot  coordinates {(-3/32,2/32) (3/32,-2/32) };
		\end{scope}
	\end{scope}
\begin{scope}[shift={(0,0)}]
\draw [thick,cborange] plot  coordinates {(0,6/8) (0,0) };
\draw [thick,cborange] plot [smooth, tension=1] coordinates {(0,0) (1/16,-4/8) (3/16+1/32,-7/8)  };
\fill [black] (0,1/2) circle (0.06);
\draw (0,1/2+3/16) node [left,red] { \scriptsize$\gamma\text{-}2$};
\fill [black] (0,0) circle (0.06);
\draw (0,0+3/16) node [left,red] { \scriptsize$\gamma\text{-}2$};
\end{scope}
	\begin{scope}[shift={(0+1/16,-1/2)}]
		\draw [thick,black] plot  coordinates {(-4/6,-4/4) (1/6,1/4) };
		\fill [black] (0,0) circle (0.06);
								\draw (0,0+0/16) node [right,red] { \scriptsize$1$};
		\begin{scope}[shift={(-2/9,-3/9)}]
			\draw [thick,black] plot  coordinates {(-3/32,2/32) (3/32,-2/32) };
		\end{scope}
		\begin{scope}[shift={(-4/9,-6/9)}]
			\draw [thick,black] plot  coordinates {(-3/32,2/32) (3/32,-2/32) };
		\end{scope}
	\end{scope}
\end{tikzpicture}}

\def\FigTwoCaseDx{\begin{tikzpicture}[baseline=13pt]
		\begin{scope}[shift={(3/4+1/2,0+1/2)}]
				\draw [thick,cborange] plot [smooth, tension=1] coordinates {(-3/3,2/12) (1/4,-6/12) (-1/3,-14/12) };
			\draw [thick,cborange] plot [smooth, tension=1] coordinates {(3/3,2/12) (-1/4,-6/12) (1/3,-14/12) };
			\fill [black] (0,-0.26) circle (0.06);
			\draw (0,-0.26+1/16) node [above,red] { \scriptsize $\alpha\texttt{+}\gamma\text{-}4$};
			\fill [black] (0,-0.99) circle (0.06);
			\draw (0,-0.99-2/16) node [below,red] { \scriptsize $\alpha\texttt{+}\gamma\text{-}4$};
			\fill [black] (-3/4+0.14/3,0.14/2) circle (0.06);
			\draw(3/4-0.14/3+1/6,0.14/2) node [above,red] { \scriptsize $1$};
			\fill [black] (3/4-0.14/3,0.14/2) circle (0.06);
			\draw (-3/4+0.14/3-1/16,0.14/2) node [above,red] { $\frac{2\text{-}\gamma}{2}$};
		\end{scope}
		\begin{scope}[shift={(-2/2,-1)}, yscale=-1]
			\draw [thick,cborange] plot [smooth, tension=1] coordinates {(-3/3,2/12) (1/4,-6/12) (-1/3,-14/12) };
		\draw [thick,cborange] plot [smooth, tension=1] coordinates {(3/3,2/12) (-1/4,-6/12) (1/3,-14/12) };
		\fill [black] (0,-0.26) circle (0.06);
		\draw (0,-0.26+2/16) node [below,red] { \scriptsize $\beta\texttt{+}\gamma\text{-}4$};
		\fill [black] (0,-0.99) circle (0.06);
		\draw (0,-0.99-1/16) node [above,red] { \scriptsize $\beta\texttt{+}\gamma\text{-}4$};
		\fill [black] (-3/4+0.14/3,0.14/2) circle (0.06);
		\draw(3/4-0.14/3+1/6,0.14/2) node [below,red] { $\frac{2\text{-}\gamma}{2}$};
		\fill [black] (3/4-0.14/3,0.14/2) circle (0.06);
		\draw (-3/4+0.14/3-1/16,0.14/2) node [below,red]{ \scriptsize $1$}; 
	\end{scope}
	\draw [thick, black, shorten >= -0.25cm, shorten <=-0.25cm] (0.14/3+1/2,0.14/2+1/2)--(-0.14/3 -1/4,-0.14/2-1);
	
	\begin{scope} [shift={(-1/8+1/4,-1/4 )}]
		\fill [black] (0,0) circle (0.06);
				\draw(-4/16,-1/16)  node [above,red] { $\frac{\gamma}{2}$};
	\draw [thick,black] plot  coordinates {(4/12,-8/8) (-1/12,1/4) };
	\begin{scope}[shift={(1/9,-3/9)}]
		\draw [thick,black] plot  coordinates {(8/64,3/64) (-8/64,-3/64) };
	\end{scope}		
	\begin{scope}[shift={(2/9,-6/9)}]
		\draw [thick,black] plot  coordinates {(8/64,3/64) (-8/64,-3/64) };
	\end{scope}
	\end{scope}

		\begin{scope}[shift={(3/2+1/2,0+1/2)}]
			\begin{scope}[shift={(-1/48,0)}]
				\draw [thick,black] plot  coordinates {(4/12,-8/8) (-1/12,1/4) };
				\begin{scope}[shift={(1/9,-3/9)}]
						\draw [thick,black] plot  coordinates {(8/64,3/64) (-8/64,-3/64) };
				\end{scope}		
				\begin{scope}[shift={(2/9,-6/9)}]
					\draw [thick,black] plot  coordinates {(8/64,3/64) (-8/64,-3/64) };
				\end{scope}
			\end{scope}
		\end{scope}
		\begin{scope}[shift={(-3/2-1/4,-1)}, yscale=-1, xscale=-1]
			\begin{scope}[shift={(-1/48,0)}]
			\draw [thick,black] plot  coordinates {(4/12,-8/8) (-1/12,1/4) };
			\begin{scope}[shift={(1/9,-3/9)}]
				\draw [thick,black] plot  coordinates {(8/64,3/64) (-8/64,-3/64) };
			\end{scope}		
			\begin{scope}[shift={(2/9,-6/9)}]
				\draw [thick,black] plot  coordinates {(8/64,3/64) (-8/64,-3/64) };
			\end{scope}
		\end{scope}
	\end{scope}
	\end{tikzpicture}}

\def\FigTwoCaseDix{\begin{tikzpicture}[baseline=13pt]
		\begin{scope}[shift={(3/4+1/2,0+1/2)}]
			\draw [thick,cborange] plot [smooth, tension=1] coordinates {(-3/3,2/12) (1/4,-6/12) (-1/3,-14/12) };
		\draw [thick,cborange] plot [smooth, tension=1] coordinates {(3/3,2/12) (-1/4,-6/12) (1/3,-14/12) };
		\fill [black] (0,-0.26) circle (0.06);
		\draw (0,-0.26+1/16) node [above,red] { \scriptsize $\alpha\texttt{+}\gamma\text{-}4$};
		\fill [black] (0,-0.99) circle (0.06);
		\draw (0,-0.99-2/16) node [below,red] { \scriptsize $\alpha\texttt{+}\gamma\text{-}4$};
		\fill [black] (-3/4+0.14/3,0.14/2) circle (0.06);
		\draw(3/4-0.14/3+1/6,0.14/2) node [above,red] { \scriptsize $1$};
		\fill [black] (3/4-0.14/3,0.14/2) circle (0.06);
		\draw (-3/4+0.14/3-1/16,0.14/2) node [above,red] { $\frac{2\text{-}\gamma}{2}$};
		\end{scope}
		\begin{scope}[shift={(-2/2,-1)}, yscale=-1]
		\draw [thick, cborange] plot coordinates {(-3/4+0.14/3-3/8,0.14/2) (3/4-0.14/3+3/8,0.14/2) };
		\draw (0,0) node [below] {\hbox{\footnotesize$\geqone$}};
		\fill [black] (-3/4+0.14/3,0.14/2) circle (0.06);
		\draw (-3/4+0.14/3-1/16,0.14/2) node [below,red]  {\scriptsize $1$};
		\fill [black] (3/4-0.14/3,0.14/2) circle (0.06);
		\draw(3/4-0.14/3+1/6,0.14/2) node [below,red]{ $\frac{\beta\text{-}\gamma}{4}$};
		\end{scope}
	\draw [thick, black, shorten >= -0.25cm, shorten <=-0.25cm] (0.14/3+1/2,0.14/2+1/2)--(-0.14/3 -1/4,-0.14/2-1);
	
		\begin{scope} [shift={(-1/8+1/4,-1/4 )}]
		\fill [black] (0,0) circle (0.06);
				\draw(-4/16,-1/16)  node [above,red] { $\frac{\gamma}{2}$};
		\draw [thick,black] plot  coordinates {(4/12,-8/8) (-1/12,1/4) };
		\begin{scope}[shift={(1/9,-3/9)}]
			\draw [thick,black] plot  coordinates {(8/64,3/64) (-8/64,-3/64) };
		\end{scope}		
		\begin{scope}[shift={(2/9,-6/9)}]
			\draw [thick,black] plot  coordinates {(8/64,3/64) (-8/64,-3/64) };
		\end{scope}
	\end{scope}
		
		\begin{scope}[shift={(3/2+1/2,0+1/2)}]
					\begin{scope}[shift={(-1/48,0)}]
			\draw [thick,black] plot  coordinates {(4/12,-8/8) (-1/12,1/4) };
			\begin{scope}[shift={(1/9,-3/9)}]
				\draw [thick,black] plot  coordinates {(8/64,3/64) (-8/64,-3/64) };
			\end{scope}		
			\begin{scope}[shift={(2/9,-6/9)}]
				\draw [thick,black] plot  coordinates {(8/64,3/64) (-8/64,-3/64) };
			\end{scope}
		\end{scope}
		\end{scope}
		\begin{scope}[shift={(-3/2-1/4,-1)}, yscale=-1, xscale=-1]
						\begin{scope}[shift={(-1/48,0)}]
				\draw [thick,black] plot  coordinates {(4/12,-8/8) (-1/12,1/4) };
				\begin{scope}[shift={(1/9,-3/9)}]
					\draw [thick,black] plot  coordinates {(8/64,3/64) (-8/64,-3/64) };
				\end{scope}		
				\begin{scope}[shift={(2/9,-6/9)}]
					\draw [thick,black] plot  coordinates {(8/64,3/64) (-8/64,-3/64) };
				\end{scope}
			\end{scope}
		\end{scope}
\end{tikzpicture}}

\def\FigTwoCaseDviii{\begin{tikzpicture}[baseline=13pt]
		\begin{scope}[shift={(3/4+1/2,0+1/2)}]
			\draw [thick,cborange] plot [smooth, tension=1] coordinates {(-3/3,2/12) (1/4,-6/12) (-1/3,-14/12) };
			\draw [thick,cborange] plot [smooth, tension=1] coordinates {(3/3,2/12) (-1/4,-6/12) (1/3,-14/12) };
			\fill [black] (0,-0.26) circle (0.06);
						\draw (0,-0.26+1/16) node [above,red] { \scriptsize $\alpha\texttt{+}\gamma\text{-}4$};
			\fill [black] (0,-0.99) circle (0.06);
					\draw (0,-0.99-2/16) node [below,red] { \scriptsize $\alpha\texttt{+}\gamma\text{-}4$};
			\fill [black] (-3/4+0.14/3,0.14/2) circle (0.06);
				\draw(3/4-0.14/3+1/6,0.14/2) node [above,red] { \scriptsize $1$};
			\fill [black] (3/4-0.14/3,0.14/2) circle (0.06);
				\draw (-3/4+0.14/3-1/16,0.14/2) node [above,red] { $\frac{2\text{-}\gamma}{2}$};
		\end{scope}
		\begin{scope}[shift={(-2/2,-1)}, yscale=-1]
			\draw [thick, cborange] plot coordinates {(-3/4+0.14/3-3/8,0.14/2) (3/4-0.14/3+3/8,0.14/2) };
		\draw (0,-7/8) node [above] {\hbox{\footnotesize$\geqone$}};
		\draw [thick, cbblue] plot coordinates {(0,0.14/2+3/8) (0,0.14/2-8/8) } ;
		\fill [black] (-3/4+0.14/3,0.14/2) circle (0.06);
		\draw (-3/4+0.14/3-1/16,0.14/2) node [below,red] { $\frac{\beta\texttt{+}\gamma}{4}$};
		\fill [black] (3/4-0.14/3,0.14/2) circle (0.06);
		\draw(3/4-0.14/3+1/6,0.14/2) node [below,red] { $\frac{\beta\text{-}\gamma}{4}$};
		\fill [black] (0,0.14/2) circle (0.06);
		\draw(0+4/16+3/16-0/32,0.14/2-1/16) node [above,red] { $\frac{4\text{-}\beta\text{-}\gamma}{12}$};
		\end{scope}
		\draw [thick, black, shorten >= -0.25cm, shorten <=-0.25cm] (0.14/3+1/2,0.14/2+1/2)--(-0.14/3 -1/4,-0.14/2-1);
			\begin{scope} [shift={(-1/8+1/4,-1/4 )}]
			\fill [black] (0,0) circle (0.06);
					\draw(-4/16,-1/16)  node [above,red] { $\frac{\gamma}{2}$};
			\draw [thick,black] plot  coordinates {(4/12,-8/8) (-1/12,1/4) };
			\begin{scope}[shift={(1/9,-3/9)}]
				\draw [thick,black] plot  coordinates {(8/64,3/64) (-8/64,-3/64) };
			\end{scope}		
			\begin{scope}[shift={(2/9,-6/9)}]
				\draw [thick,black] plot  coordinates {(8/64,3/64) (-8/64,-3/64) };
			\end{scope}
		\end{scope}
		\begin{scope}[shift={(3/2+1/2,0+1/2)}]
	\begin{scope}[shift={(-1/48,0)}]
	\draw [thick,black] plot  coordinates {(4/12,-8/8) (-1/12,1/4) };
	\begin{scope}[shift={(1/9,-3/9)}]
		\draw [thick,black] plot  coordinates {(8/64,3/64) (-8/64,-3/64) };
	\end{scope}		
	\begin{scope}[shift={(2/9,-6/9)}]
		\draw [thick,black] plot  coordinates {(8/64,3/64) (-8/64,-3/64) };
	\end{scope}
\end{scope}
		\end{scope}
		\begin{scope}[shift={(-3/2-1/4,-1)}, yscale=-1, xscale=-1]
				\begin{scope}[shift={(-1/48,0)}]
				\draw [thick,black] plot  coordinates {(4/12,-8/8) (-1/12,1/4) };
				\begin{scope}[shift={(1/9,-3/9)}]
					\draw [thick,black] plot  coordinates {(8/64,3/64) (-8/64,-3/64) };
				\end{scope}		
				\begin{scope}[shift={(2/9,-6/9)}]
					\draw [thick,black] plot  coordinates {(8/64,3/64) (-8/64,-3/64) };
				\end{scope}
			\end{scope}
		\end{scope}
\end{tikzpicture}}

\def\FigTwoCaseDvii{\begin{tikzpicture}[baseline=13pt]
		\begin{scope}[shift={(3/4+1/2,0+1/2)}]
		\draw [thick, cborange] plot coordinates {(-3/4+0.14/3-3/8,0.14/2) (3/4-0.14/3+3/8,0.14/2) };
		\draw (0,0) node [below] {\hbox{\footnotesize$\geqone$}};
		\fill [black] (-3/4+0.14/3,0.14/2) circle (0.06);
		\draw (-3/4+0.14/3-1/16,0.14/2) node [above,red] { $\frac{\alpha\text{-}\gamma}{4}$};
		\fill [black] (3/4-0.14/3,0.14/2) circle (0.06);
		\draw(3/4-0.14/3+1/6,0.14/2) node [above,red] {\scriptsize $1$};
		\end{scope}
		\begin{scope}[shift={(-2/2,-1)}, yscale=-1]
		\draw [thick, cborange] plot coordinates {(-3/4+0.14/3-3/8,0.14/2) (3/4-0.14/3+3/8,0.14/2) };
		\draw (0,0) node [below] {\hbox{\footnotesize$\geqone$}};
		\fill [black] (-3/4+0.14/3,0.14/2) circle (0.06);
		\draw (-3/4+0.14/3-1/16,0.14/2) node [below,red]  {\scriptsize $1$};
		\fill [black] (3/4-0.14/3,0.14/2) circle (0.06);
		\draw(3/4-0.14/3+1/6,0.14/2) node [below,red]{ $\frac{\beta\text{-}\gamma}{4}$};
		\end{scope}
		\draw [thick, black, shorten >= -0.25cm, shorten <=-0.25cm] (0.14/3+1/2,0.14/2+1/2)--(-0.14/3 -1/4,-0.14/2-1);
		
			\begin{scope} [shift={(-1/8+1/4,-1/4 )}]
			\fill [black] (0,0) circle (0.06);
					\draw(-4/16,-1/16)  node [above,red] { $\frac{\gamma}{2}$};
			\draw [thick,black] plot  coordinates {(4/12,-8/8) (-1/12,1/4) };
			\begin{scope}[shift={(1/9,-3/9)}]
				\draw [thick,black] plot  coordinates {(8/64,3/64) (-8/64,-3/64) };
			\end{scope}		
			\begin{scope}[shift={(2/9,-6/9)}]
				\draw [thick,black] plot  coordinates {(8/64,3/64) (-8/64,-3/64) };
			\end{scope}
		\end{scope}
		
		\begin{scope}[shift={(3/2+1/2,0+1/2)}]
					\begin{scope}[shift={(-1/48,0)}]
			\draw [thick,black] plot  coordinates {(4/12,-8/8) (-1/12,1/4) };
			\begin{scope}[shift={(1/9,-3/9)}]
				\draw [thick,black] plot  coordinates {(8/64,3/64) (-8/64,-3/64) };
			\end{scope}		
			\begin{scope}[shift={(2/9,-6/9)}]
				\draw [thick,black] plot  coordinates {(8/64,3/64) (-8/64,-3/64) };
			\end{scope}
		\end{scope}
		\end{scope}
		\begin{scope}[shift={(-3/2-1/4,-1)}, yscale=-1, xscale=-1]
			\begin{scope}[shift={(-1/48,0)}]
	\draw [thick,black] plot  coordinates {(4/12,-8/8) (-1/12,1/4) };
	\begin{scope}[shift={(1/9,-3/9)}]
		\draw [thick,black] plot  coordinates {(8/64,3/64) (-8/64,-3/64) };
	\end{scope}		
	\begin{scope}[shift={(2/9,-6/9)}]
		\draw [thick,black] plot  coordinates {(8/64,3/64) (-8/64,-3/64) };
	\end{scope}
\end{scope}
		\end{scope}
\end{tikzpicture}}

\def\FigTwoCaseDvi{\begin{tikzpicture}[baseline=13pt]
		\begin{scope}[shift={(3/4+1/2,0+1/2)}]
	\draw [thick, cborange] plot coordinates {(-3/4+0.14/3-3/8,0.14/2) (3/4-0.14/3+3/8,0.14/2) };
	\draw (0,0) node [below] {\hbox{\footnotesize$\geqone$}};
	\fill [black] (-3/4+0.14/3,0.14/2) circle (0.06);
		\draw (-3/4+0.14/3-1/16,0.14/2) node [above,red] { $\frac{\alpha\text{-}\gamma}{4}$};
	\fill [black] (3/4-0.14/3,0.14/2) circle (0.06);
			\draw(3/4-0.14/3+1/6,0.14/2) node [above,red] {\scriptsize $1$};
		\end{scope}
	\begin{scope}[shift={(-2/2,-1)}, yscale=-1]
					\draw [thick, cborange] plot coordinates {(-3/4+0.14/3-3/8,0.14/2) (3/4-0.14/3+3/8,0.14/2) };
				\draw (0,-7/8) node [above] {\hbox{\footnotesize$\geqone$}};
				\draw [thick, cbblue] plot coordinates {(0,0.14/2+3/8) (0,0.14/2-8/8) } ;
				\fill [black] (-3/4+0.14/3,0.14/2) circle (0.06);
				\draw (-3/4+0.14/3-1/16,0.14/2) node [below,red] { $\frac{\beta\texttt{+}\gamma}{4}$};
				\fill [black] (3/4-0.14/3,0.14/2) circle (0.06);
				\draw(3/4-0.14/3+1/6,0.14/2) node [below,red] { $\frac{\beta\text{-}\gamma}{4}$};
				\fill [black] (0,0.14/2) circle (0.06);
				\draw(0+4/16+3/16-0/32,0.14/2-1/16) node [above,red] { $\frac{4\text{-}\beta\text{-}\gamma}{12}$};
	\end{scope}
\draw [thick, black, shorten >= -0.25cm, shorten <=-0.25cm] (0.14/3+1/2,0.14/2+1/2)--(-0.14/3 -1/4,-0.14/2-1);

	\begin{scope} [shift={(-1/8+1/4,-1/4 )}]
	\fill [black] (0,0) circle (0.06);
		\draw(-4/16,-1/16)  node [above,red] { $\frac{\gamma}{2}$};
	\draw [thick,black] plot  coordinates {(4/12,-8/8) (-1/12,1/4) };
	\begin{scope}[shift={(1/9,-3/9)}]
		\draw [thick,black] plot  coordinates {(8/64,3/64) (-8/64,-3/64) };
	\end{scope}		
	\begin{scope}[shift={(2/9,-6/9)}]
		\draw [thick,black] plot  coordinates {(8/64,3/64) (-8/64,-3/64) };
	\end{scope}
\end{scope}

		\begin{scope}[shift={(3/2+1/2,1/2)}]
						\begin{scope}[shift={(-1/48,0)}]
				\draw [thick,black] plot  coordinates {(4/12,-8/8) (-1/12,1/4) };
				\begin{scope}[shift={(1/9,-3/9)}]
					\draw [thick,black] plot  coordinates {(8/64,3/64) (-8/64,-3/64) };
				\end{scope}		
				\begin{scope}[shift={(2/9,-6/9)}]
					\draw [thick,black] plot  coordinates {(8/64,3/64) (-8/64,-3/64) };
				\end{scope}
			\end{scope}
			\end{scope}
		\begin{scope}[shift={(-3/2-1/4,-1)}, yscale=-1, xscale=-1]
						\begin{scope}[shift={(-1/48,0)}]
				\draw [thick,black] plot  coordinates {(4/12,-8/8) (-1/12,1/4) };
				\begin{scope}[shift={(1/9,-3/9)}]
					\draw [thick,black] plot  coordinates {(8/64,3/64) (-8/64,-3/64) };
				\end{scope}		
				\begin{scope}[shift={(2/9,-6/9)}]
					\draw [thick,black] plot  coordinates {(8/64,3/64) (-8/64,-3/64) };
				\end{scope}
			\end{scope}
		\end{scope}
\end{tikzpicture}}

\def\FigTwoCaseDv{\begin{tikzpicture}[baseline=13pt]
		\begin{scope}[shift={(3/4+1/2,1/2)}]
			\draw [thick, cborange] plot coordinates {(-3/4+0.14/3-3/8,0.14/2) (3/4-0.14/3+3/8,0.14/2) };
			\draw (0,-7/8) node [below] {\hbox{\footnotesize$\geqone$}};
			\draw [thick, cbblue] plot coordinates {(0,0.14/2+3/8) (0,0.14/2-8/8) } ;
			\fill [black] (-3/4+0.14/3,0.14/2) circle (0.06);
			\draw (-3/4+0.14/3-1/16,0.14/2) node [above,red] { $\frac{\alpha\texttt{-}\gamma}{4}$};
			\fill [black] (3/4-0.14/3,0.14/2) circle (0.06);
			\draw(3/4-0.14/3+1/6,0.14/2) node [above,red] { $\frac{\alpha\texttt{+}\gamma}{4}$};
			\fill [black] (0,0.14/2) circle (0.06);
			\draw(0-4/16-3/16+1/64,0.14/2-1/16) node [below,red] { $\frac{4\text{-}\alpha\text{-}\gamma}{12}$};
		\end{scope}
		\begin{scope}[shift={(-2/2,-1)}, yscale=-1]
				\draw [thick, cborange] plot coordinates {(-3/4+0.14/3-3/8,0.14/2) (3/4-0.14/3+3/8,0.14/2) };
			\draw (0,-7/8) node [above] {\hbox{\footnotesize$\geqone$}};
			\draw [thick, cbblue] plot coordinates {(0,0.14/2+3/8) (0,0.14/2-8/8) } ;
			\fill [black] (-3/4+0.14/3,0.14/2) circle (0.06);
			\draw (-3/4+0.14/3-1/16,0.14/2) node [below,red] { $\frac{\beta\texttt{+}\gamma}{4}$};
			\fill [black] (3/4-0.14/3,0.14/2) circle (0.06);
			\draw(3/4-0.14/3+1/6,0.14/2) node [below,red] { $\frac{\beta\text{-}\gamma}{4}$};
			\fill [black] (0,0.14/2) circle (0.06);
			\draw(0+4/16+3/16-0/32,0.14/2-1/16) node [above,red] { $\frac{4\text{-}\beta\text{-}\gamma}{12}$};
				
		\end{scope}

\draw [thick, black, shorten >= -0.25cm, shorten <=-0.25cm] (0.14/3+1/2,0.14/2+1/2)--(-0.14/3 -1/4,-0.14/2-1);

	\begin{scope} [shift={(-1/8+1/4,-1/4 )}]
	\fill [black] (0,0) circle (0.06);
		\draw(-4/16,-1/16)  node [above,red] { $\frac{\gamma}{2}$};
	\draw [thick,black] plot  coordinates {(4/12,-8/8) (-1/12,1/4) };
	\begin{scope}[shift={(1/9,-3/9)}]
		\draw [thick,black] plot  coordinates {(8/64,3/64) (-8/64,-3/64) };
	\end{scope}		
	\begin{scope}[shift={(2/9,-6/9)}]
		\draw [thick,black] plot  coordinates {(8/64,3/64) (-8/64,-3/64) };
	\end{scope}
\end{scope}

		\begin{scope}[shift={(3/2+1/2,1/2)}]
					\begin{scope}[shift={(-1/48,0)}]
			\draw [thick,black] plot  coordinates {(4/12,-8/8) (-1/12,1/4) };
			\begin{scope}[shift={(1/9,-3/9)}]
				\draw [thick,black] plot  coordinates {(8/64,3/64) (-8/64,-3/64) };
			\end{scope}		
			\begin{scope}[shift={(2/9,-6/9)}]
				\draw [thick,black] plot  coordinates {(8/64,3/64) (-8/64,-3/64) };
			\end{scope}
		\end{scope}
		\end{scope}
		\begin{scope}[shift={(-3/2-1/4,-1)}, yscale=-1, xscale=-1]
						\begin{scope}[shift={(-1/48,0)}]
				\draw [thick,black] plot  coordinates {(4/12,-8/8) (-1/12,1/4) };
				\begin{scope}[shift={(1/9,-3/9)}]
					\draw [thick,black] plot  coordinates {(8/64,3/64) (-8/64,-3/64) };
				\end{scope}		
				\begin{scope}[shift={(2/9,-6/9)}]
					\draw [thick,black] plot  coordinates {(8/64,3/64) (-8/64,-3/64) };
				\end{scope}
			\end{scope}
		\end{scope}
\end{tikzpicture}}

\def\FigTwoCaseDxi{\begin{tikzpicture}[baseline=13pt]
		\draw [thick,black] plot  coordinates {(-8/6,-8/4) (1/6,1/4) };
			\begin{scope} [shift={(-10/27,-15/27 )}]
				\draw [thick,black] plot  coordinates {(8/12,-8/8) (-1/6,1/4) };
			\fill [black] (0,0) circle (0.06);
						\draw (-1/16,0) node [left,red] { $\frac{\gamma}{2}$};
			\begin{scope}[shift={(2/9,-3/9)}]
				\draw [thick,black] plot  coordinates {(3/32,2/32) (-3/32,-2/32) };
			\end{scope}
			\begin{scope}[shift={(4/9,-6/9)}]
				\draw [thick,black] plot  coordinates {(3/32,2/32) (-3/32,-2/32) };
			\end{scope}
		\end{scope}
			\begin{scope} [shift={(-20/27,-30/27 )}]
				\draw [thick,black] plot  coordinates {(8/12,-8/8) (-1/6,1/4) };
			\fill [black] (0,0) circle (0.06);
							\draw (-1/16,0) node [left,red] { $\frac{\gamma}{2}$};
			\begin{scope}[shift={(2/9,-3/9)}]
				\draw [thick,black] plot  coordinates {(3/32,2/32) (-3/32,-2/32) };
			\end{scope}
			\begin{scope}[shift={(4/9,-6/9)}]
				\draw [thick,black] plot  coordinates {(3/32,2/32) (-3/32,-2/32) };
			\end{scope}
		\end{scope}
		\draw [thick,cborange] plot  coordinates {(-1/2,0) (4/2,0) };
		\fill [black] (0,0) circle (0.06);
			\draw (-1/16,0) node [above,red] { $\frac{\alpha\text{-}\gamma}{4}$};
		\draw [thick,cbblue] plot  coordinates {(3/4,1/3) (3/4,-6/4) };
		\draw (3/4,-3/2) node [right] {\hbox{\footnotesize$\geqone$}};
		\fill [black] (3/4,0) circle (0.06);
			\draw(3/4-4/16-3/16+1/64,0.14/2-1/16) node [below,red] { $\frac{4\text{-}\alpha\text{-}\gamma}{12}$};

		\begin{scope}[shift={(3/2,0)}]
			\draw [thick,black] plot  coordinates {(8/12,-8/8) (-1/6,1/4) };
			\fill [black] (0,0) circle (0.06);
				\draw (1/16,0) node [above,red] { $\frac{\alpha\texttt{+}\gamma}{4}$};
			\begin{scope}[shift={(2/9,-3/9)}]
				\draw [thick,black] plot  coordinates {(3/32,2/32) (-3/32,-2/32) };
			\end{scope}
			\begin{scope}[shift={(4/9,-6/9)}]
				\draw [thick,black] plot  coordinates {(3/32,2/32) (-3/32,-2/32) };
			\end{scope}
		\end{scope}
		\draw [thick,cbblue] plot  coordinates {(-10/9-5/6,-15/9+5/4) (-10/9+1/6,-15/9-1/4) };
		\fill [black] (-10/9,-15/9) circle (0.06);
			\draw (-10/9,-15/9) node [left,red] { $\frac{2\text{-}\gamma}{6}$};
		\draw (-10/9-5/6,-15/9+5/4+1/6) node [] {\hbox{\footnotesize$\geqone$}};
\end{tikzpicture}}

\def\FigTwoCaseDxii{\begin{tikzpicture}[baseline=13pt]
		\draw [thick,black] plot  coordinates {(-8/6,-8/4) (1/6,1/4) };
		\begin{scope} [shift={(-10/27,-15/27 )}]
			\draw [thick,black] plot  coordinates {(8/12,-8/8) (-1/6,1/4) };
			\fill [black] (0,0) circle (0.06);
			\draw (-1/16,0) node [left,red] { $\frac{\gamma}{2}$};
			\begin{scope}[shift={(2/9,-3/9)}]
				\draw [thick,black] plot  coordinates {(3/32,2/32) (-3/32,-2/32) };
			\end{scope}
			\begin{scope}[shift={(4/9,-6/9)}]
				\draw [thick,black] plot  coordinates {(3/32,2/32) (-3/32,-2/32) };
			\end{scope}
		\end{scope}
		\begin{scope} [shift={(-20/27,-30/27 )}]
			\draw [thick,black] plot  coordinates {(8/12,-8/8) (-1/6,1/4) };
			\fill [black] (0,0) circle (0.06);
			\draw (-1/16,0) node [left,red] { $\frac{\gamma}{2}$};
			\begin{scope}[shift={(2/9,-3/9)}]
				\draw [thick,black] plot  coordinates {(3/32,2/32) (-3/32,-2/32) };
			\end{scope}
			\begin{scope}[shift={(4/9,-6/9)}]
				\draw [thick,black] plot  coordinates {(3/32,2/32) (-3/32,-2/32) };
			\end{scope}
		\end{scope}
		\draw [thick,cborange] plot  coordinates {(-1/2,0) (4/2,0) };
		\fill [black] (0,0) circle (0.06);
		\draw (-1/16,0) node [above,red] { $\frac{\alpha\text{-}\gamma}{4}$};
		\draw (3/4,0) node [below] {\hbox{\footnotesize$\geqone$}};

		\begin{scope}[shift={(3/2,0)}]
			\draw [thick,black] plot  coordinates {(8/12,-8/8) (-1/6,1/4) };
			\fill [black] (0,0) circle (0.06);
			\draw (1/16,0) node [above,red] { $\frac{\alpha\texttt{+}\gamma}{4}$};
			\begin{scope}[shift={(2/9,-3/9)}]
				\draw [thick,black] plot  coordinates {(3/32,2/32) (-3/32,-2/32) };
			\end{scope}
			\begin{scope}[shift={(4/9,-6/9)}]
				\draw [thick,black] plot  coordinates {(3/32,2/32) (-3/32,-2/32) };
			\end{scope}
		\end{scope}
		\draw [thick,cbblue] plot  coordinates {(-10/9-5/6,-15/9+5/4) (-10/9+1/6,-15/9-1/4) };
		\fill [black] (-10/9,-15/9) circle (0.06);
		\draw (-10/9,-15/9) node [left,red] { $\frac{2\text{-}\gamma}{6}$};
		\draw (-10/9-5/6,-15/9+5/4+1/6) node [] {\hbox{\footnotesize$\geqone$}};
\end{tikzpicture}}

\def\FigTwoCaseDxiii{\begin{tikzpicture}[baseline=13pt]
		\draw [thick,black] plot  coordinates {(-16/12,-16/8) (1/6,1/4) };
\begin{scope} [shift={(-10/27,-15/27 )}]
	\draw [thick,black] plot  coordinates {(8/12,-8/8) (-1/6,1/4) };
	\fill [black] (0,0) circle (0.06);
	\draw (-1/16,0) node [left,red] { $\frac{\gamma}{2}$};
	\begin{scope}[shift={(2/9,-3/9)}]
		\draw [thick,black] plot  coordinates {(3/32,2/32) (-3/32,-2/32) };
	\end{scope}
	\begin{scope}[shift={(4/9,-6/9)}]
		\draw [thick,black] plot  coordinates {(3/32,2/32) (-3/32,-2/32) };
	\end{scope}
\end{scope}
\begin{scope} [shift={(-20/27,-30/27 )}]
	\draw [thick,black] plot  coordinates {(8/12,-8/8) (-1/6,1/4) };
	\fill [black] (0,0) circle (0.06);
	\draw (-1/16,0) node [left,red] { $\frac{\gamma}{2}$};
	\begin{scope}[shift={(2/9,-3/9)}]
		\draw [thick,black] plot  coordinates {(3/32,2/32) (-3/32,-2/32) };
	\end{scope}
	\begin{scope}[shift={(4/9,-6/9)}]
		\draw [thick,black] plot  coordinates {(3/32,2/32) (-3/32,-2/32) };
	\end{scope}
\end{scope}
		\begin{scope}[shift={(3/4,0)}]
			\draw [thick,cborange] plot [smooth, tension=1] coordinates {(-3/3,2/12) (1/4,-6/12) (-1/3,-14/12) };
			\draw [thick,cborange] plot [smooth, tension=1] coordinates {(3/3,2/12) (-1/4,-6/12) (1/3,-14/12) };
			\fill [black] (0,-0.26) circle (0.06);
								\draw (0,-0.26+1/16)  node [above,red] {\scriptsize $\alpha\texttt{+}\gamma\text{-}4$};
			\fill [black] (0,-0.99) circle (0.06);
											\draw (0,-0.99-2/16)  node [below,red] {\scriptsize $\alpha\texttt{+}\gamma\text{-}4$};
			\fill [black] (-3/4+0.14/3,0.14/2) circle (0.06);
						\draw (-3/4+0.14/3-1/16,0.14/2)  node [above,red] { $\frac{2\text{-}\gamma}{2}$};
			\fill [black] (3/4-0.14/3,0.14/2) circle (0.06);
									\draw (3/4-0.14/3+1/16,0.14/2)  node [above,red] {\scriptsize $1$};
		\end{scope}
		\begin{scope}[shift={(3/2,0)}]
			\draw [thick,black] plot  coordinates {(8/12,-8/8) (-1/6,1/4) };
			\begin{scope}[shift={(2/9,-3/9)}]
				\draw [thick,black] plot  coordinates {(3/32,2/32) (-3/32,-2/32) };
			\end{scope}
			\begin{scope}[shift={(4/9,-6/9)}]
				\draw [thick,black] plot  coordinates {(3/32,2/32) (-3/32,-2/32) };
			\end{scope}
		\end{scope}
		\draw [thick,cbblue] plot  coordinates {(-10/9-5/6,-15/9+5/4) (-10/9+1/6,-15/9-1/4) };
		\fill [black] (-10/9,-15/9) circle (0.06);
					\draw (-10/9,-15/9) node [left,red] { $\frac{2\text{-}\gamma}{6}$};
		\draw (-10/9-5/6,-15/9+5/4+1/6) node [] {\hbox{\footnotesize$\geqone$}};
\end{tikzpicture}}

\def\FigTwoCaseDxiv{\begin{tikzpicture}[baseline=13pt]
		\draw [thick,black] plot  coordinates {(-14/12,-14/8) (1/6,1/4) };
		\begin{scope} [shift={(-4/9,-6/9 )}, xscale=-1, yscale=-1]
			\draw [thick,black] plot  coordinates {(8/12,-8/8) (-1/6,1/4) };
			\fill [black] (0,0) circle (0.06);
			\draw (0,0+1/16)  node [below,red] { $\frac{\gamma}{2}$};
			\begin{scope}[shift={(2/9,-3/9)}]
				\draw [thick,black] plot  coordinates {(3/32,2/32) (-3/32,-2/32) };
			\end{scope}
			\begin{scope}[shift={(4/9,-6/9)}]
				\draw [thick,black] plot  coordinates {(3/32,2/32) (-3/32,-2/32) };
			\end{scope}
		\end{scope}
		\begin{scope} [shift={(-7/9,-10.5/9 )}, xscale=-1,yscale=-1 ]
			\draw [thick,black] plot  coordinates {(8/12,-8/8) (-1/6,1/4) };
			\fill [black] (0,0) circle (0.06);
				\draw (0,0+1/16)  node [below,red] { $\frac{\gamma}{2}$};
			\begin{scope}[shift={(2/9,-3/9)}]
				\draw [thick,black] plot  coordinates {(3/32,2/32) (-3/32,-2/32) };
			\end{scope}
			\begin{scope}[shift={(4/9,-6/9)}]
				\draw [thick,black] plot  coordinates {(3/32,2/32) (-3/32,-2/32) };
			\end{scope}
		\end{scope}
		\draw [thick,cborange] plot  coordinates {(-1/4,1/6) (6/4,-6/6) };
		\draw [thick,cborange] plot  coordinates {(5/2+1/4,1/6) (5/2-6/4,-6/6) };
		\fill [black] (0,0) circle (0.06);
		\draw (0-1/16,0+1/16)  node [above,red] { $\frac{\delta \text{-}\gamma}{4}$};
		\fill [black] (5/4,-5/6) circle (0.06);
				\draw (5/4,-5/6-2/16)  node [below,red] { $\frac{\alpha\texttt{+}\gamma \text{-} 2\delta}{4}$};
		\begin{scope}[shift={(3/5,-2/5)}]
			\draw [thick,cbblue] plot coordinates {(1/12,1/3) (-7/24,-7/6) };
			\fill [black] (0,0) circle (0.06);
						\draw (0+9/16,0+1/16)  node [above,red] { $\frac{4\text{-}3\delta \texttt{+} \gamma}{12}$};
			\draw (-6.5/24,-6.5/6) node [below] {\hbox{\footnotesize$\geqone$}};
		\end{scope}
		\begin{scope}[shift={(5/2-3/5,-2/5)}]
			\draw [thick,cbblue] plot coordinates {(-1/12,1/3) (7/24,-7/6) };
			\fill [black] (0,0) circle (0.06);
			\draw (0+8/16,0+1/16)  node [below,red] { $\frac{4\text{-}\alpha \text{-} \delta}{12}$};
			\draw (6.5/24,-6.5/6) node [below] {\hbox{\footnotesize$\geqone$}};
		\end{scope}
		\begin{scope}[shift={(5/2,0)}]
			\draw [thick,black] plot  coordinates {(8/12,-8/8) (-1/6,1/4) };
			\fill [black] (0,0) circle (0.06);
				\draw (0+1/16,0+1/16)  node [above,red] { $\frac{\alpha\texttt{+}\delta}{4}$};
			\begin{scope}[shift={(2/9,-3/9)}]
				\draw [thick,black] plot  coordinates {(3/32,2/32) (-3/32,-2/32) };
			\end{scope}
			\begin{scope}[shift={(4/9,-6/9)}]
				\draw [thick,black] plot  coordinates {(3/32,2/32) (-3/32,-2/32) };
			\end{scope}
		\end{scope}
\end{tikzpicture}}

\def\FigTwoCaseDxv{\begin{tikzpicture}[baseline=13pt]
		\draw [thick,black] plot  coordinates {(-14/12,-14/8) (1/6,1/4) };
		\begin{scope} [shift={(-4/9,-6/9 )}, xscale=-1, yscale=-1]
			\draw [thick,black] plot  coordinates {(8/12,-8/8) (-1/6,1/4) };
			\fill [black] (0,0) circle (0.06);
			\draw (0,0+1/16)  node [below,red] { $\frac{\gamma}{2}$};
			\begin{scope}[shift={(2/9,-3/9)}]
				\draw [thick,black] plot  coordinates {(3/32,2/32) (-3/32,-2/32) };
			\end{scope}
			\begin{scope}[shift={(4/9,-6/9)}]
				\draw [thick,black] plot  coordinates {(3/32,2/32) (-3/32,-2/32) };
			\end{scope}
		\end{scope}
		\begin{scope} [shift={(-7/9,-10.5/9 )}, xscale=-1,yscale=-1 ]
			\draw [thick,black] plot  coordinates {(8/12,-8/8) (-1/6,1/4) };
			\fill [black] (0,0) circle (0.06);
			\draw (0,0+1/16)  node [below,red] { $\frac{\gamma}{2}$};
			\begin{scope}[shift={(2/9,-3/9)}]
				\draw [thick,black] plot  coordinates {(3/32,2/32) (-3/32,-2/32) };
			\end{scope}
			\begin{scope}[shift={(4/9,-6/9)}]
				\draw [thick,black] plot  coordinates {(3/32,2/32) (-3/32,-2/32) };
			\end{scope}
		\end{scope}
		\draw [thick,cborange] plot  coordinates {(-1/4,1/6) (6/4,-6/6) };
		\draw [thick,cborange] plot  coordinates {(5/2+1/4,1/6) (5/2-6/4,-6/6) };
		\fill [black] (0,0) circle (0.06);
				\draw (0-1/16,0+1/16)  node [above,red] { $\frac{\delta \text{-}\gamma}{4}$};
		\fill [black] (5/4,-5/6) circle (0.06);
						\draw (5/4,-5/6-2/16)  node [below,red] { $\frac{\alpha\texttt{+}\gamma \text{-} 2\delta}{4}$};
		\begin{scope}[shift={(3/5,-2/5)}]
			\draw [thick,cbblue] plot coordinates {(1/12,1/3) (-7/24,-7/6) };
			\fill [black] (0,0) circle (0.06);
						\draw (0+9/16,0+1/16)  node [above,red] { $\frac{4\text{-}3\delta \texttt{+} \gamma}{12}$};
			\draw (-6.5/24,-6.5/6) node [below] {\hbox{\footnotesize$\geqone$}};
		\end{scope}
		\draw (1+4/4+1/12,-1+4/6) node [below] {\hbox{\footnotesize$\geqone$}};
		\begin{scope}[shift={(5/2,0)}]
			\draw [thick,black] plot  coordinates {(8/12,-8/8) (-1/6,1/4) };
			\fill [black] (0,0) circle (0.06);
							\draw (0+1/16,0+1/16)  node [above,red] { \scriptsize $1$};
			\begin{scope}[shift={(2/9,-3/9)}]
				\draw [thick,black] plot  coordinates {(3/32,2/32) (-3/32,-2/32) };
			\end{scope}
			\begin{scope}[shift={(4/9,-6/9)}]
				\draw [thick,black] plot  coordinates {(3/32,2/32) (-3/32,-2/32) };
			\end{scope}
		\end{scope}
\end{tikzpicture}}

\def\FigTwoCaseDxvi{\begin{tikzpicture}[baseline=13pt]
		\draw [thick,black] plot  coordinates {(-1-14/12,-14/8) (-1+1/6,1/4) };
		\begin{scope} [shift={(-4/9-1,-6/9 )}, xscale=-1, yscale=-1]
			\draw [thick,black] plot  coordinates {(8/12,-8/8) (-1/6,1/4) };
			\fill [black] (0,0) circle (0.06);
						\draw (0,0+1/16)  node [below,red] { $\frac{\gamma}{2}$};
			\begin{scope}[shift={(2/9,-3/9)}]
				\draw [thick,black] plot  coordinates {(3/32,2/32) (-3/32,-2/32) };
			\end{scope}
			\begin{scope}[shift={(4/9,-6/9)}]
				\draw [thick,black] plot  coordinates {(3/32,2/32) (-3/32,-2/32) };
			\end{scope}
		\end{scope}
		\begin{scope} [shift={(-7/9-1,-10.5/9 )}, xscale=-1, yscale=-1]
			\draw [thick,black] plot  coordinates {(8/12,-8/8) (-1/6,1/4) };
			\fill [black] (0,0) circle (0.06);
						\draw (0,0+1/16)  node [below,red] { $\frac{\gamma}{2}$};
			\begin{scope}[shift={(2/9,-3/9)}]
				\draw [thick,black] plot  coordinates {(3/32,2/32) (-3/32,-2/32) };
			\end{scope}
			\begin{scope}[shift={(4/9,-6/9)}]
				\draw [thick,black] plot  coordinates {(3/32,2/32) (-3/32,-2/32) };
			\end{scope}
		\end{scope}
		\draw [thick,cborange] plot [smooth, tension=1] coordinates {(4/12,1/4) (-6/12,-1/4) (-16/12,1/4) };
		\begin{scope}[shift={(3/4,0)}]
			\draw [thick,cborange] plot [smooth, tension=1] coordinates {(-3/3,2/12) (1/4,-6/12) (-1/3,-14/12) };
			\draw [thick,cborange] plot [smooth, tension=1] coordinates {(3/3,2/12) (-1/4,-6/12) (1/3,-14/12) };
			\fill [black] (0,-0.26) circle (0.06);
				\draw(0,-0.26+2/16)  node [above,red] {\scriptsize $\alpha \texttt{+} \delta \text{-}4$};
			\fill [black] (0,-0.99) circle (0.06);
							\draw(0,-0.99-2/16)  node [below,red] {\scriptsize $\alpha \texttt{+} \delta \text{-}4$};
			\fill [black] (3/4-0.14/3,0.14/2) circle (0.06);
			\draw(3/4-0.14/3+1/16,0.14/2)  node [above,red] {\scriptsize $1$};
		\end{scope}
		\fill [black] (-5/6-0.6/3,1/4-0.6/2) circle (0.06);
					\draw(-5/6-0.6/3-1/16,1/4-0.6/2+1/16)  node [above,red] {$\frac{\delta \text{-}\gamma}{4}$};
		\fill [black] (0.14,0.036) circle (0.06);
			\draw(0.14-1/16,0.036-0/16)  node [above,red] {$\frac{4\text{-}3\delta \texttt{+}\gamma}{4}$};
		\begin{scope}[shift={(3/2,0)}]
			\draw [thick,black] plot  coordinates {(8/12,-8/8) (-1/6,1/4) };
			\begin{scope}[shift={(2/9,-3/9)}]
				\draw [thick,black] plot  coordinates {(3/32,2/32) (-3/32,-2/32) };
			\end{scope}
			\begin{scope}[shift={(4/9,-6/9)}]
				\draw [thick,black] plot  coordinates {(3/32,2/32) (-3/32,-2/32) };
			\end{scope}
		\end{scope}
		\draw [thick,cbblue] plot  coordinates {(-6/12,-15/9+7/4) (-6/12,-15/9+1/4) };
		\fill [black] (-6/12,-.25) circle (0.06);
		\draw (-6/6,-10/6) node [right] {\hbox{\footnotesize$\geqone$}};
\end{tikzpicture}}

\def\FigTwoCaseDxvii{\begin{tikzpicture}[baseline=13pt]
		\draw [thick,black] plot  coordinates {(-14/12,-14/8) (1/6,1/4) };
		\begin{scope} [shift={(-4/9,-6/9 )}, xscale=-1, yscale=-1]
			\draw [thick,black] plot  coordinates {(8/12,-8/8) (-1/6,1/4) };
			\fill [black] (0,0) circle (0.06);
									\draw (0,0+1/16)  node [below,red] { $\frac{\gamma}{2}$};
			\begin{scope}[shift={(2/9,-3/9)}]
				\draw [thick,black] plot  coordinates {(3/32,2/32) (-3/32,-2/32) };
			\end{scope}
			\begin{scope}[shift={(4/9,-6/9)}]
				\draw [thick,black] plot  coordinates {(3/32,2/32) (-3/32,-2/32) };
			\end{scope}
		\end{scope}
		\begin{scope} [shift={(-7/9,-10.5/9 )},xscale=-1, yscale=-1]
			\draw [thick,black] plot  coordinates {(8/12,-8/8) (-1/6,1/4) };
			\fill [black] (0,0) circle (0.06);
									\draw (0,0+1/16)  node [below,red] { $\frac{\gamma}{2}$};
			\begin{scope}[shift={(2/9,-3/9)}]
				\draw [thick,black] plot  coordinates {(3/32,2/32) (-3/32,-2/32) };
			\end{scope}
			\begin{scope}[shift={(4/9,-6/9)}]
				\draw [thick,black] plot  coordinates {(3/32,2/32) (-3/32,-2/32) };
			\end{scope}
		\end{scope}
		\begin{scope}[shift={(3/4,0)}]
			\draw [thick,cborange] plot [smooth, tension=1] coordinates {(-3/3,2/12) (1,-8/12) (-1/2,-20/12) };
			\draw [thick,cborange] plot [smooth, tension=.8] coordinates {(1/2+3/3,2/12) (5/12,-8/12) (11/12,-16/12) };
			\fill [black] (0.69,-0.36) circle (0.06);
					\draw (0.69-2/16,-0.36)  node [above,red] { $\frac{3\alpha\text{-}8\texttt{+}\gamma}{8}$};
			\fill [black] (0.66,-1.16) circle (0.06);
			\draw (0.69,-1.16-2/16)  node [below,red] { $\frac{3\alpha\text{-}8\texttt{+}\gamma}{8}$};
			\draw (-1/2+1/4,-20/12+1/8) node [below] {\hbox{\footnotesize$\geqone$}};
		\end{scope}
		\fill [black] (0.195/3,0.195/2) circle (0.06);
		\draw (0.195/3+3/16,0.195/2)  node [below,red] { $\frac{2\text{-}\gamma}{6}$};
		\begin{scope}[shift={(2,0)}]
			\draw [thick,black] plot  coordinates {(8/12,-8/8) (-1/6,1/4) };
			\fill [black] (-0.035/3,0.035/2) circle (0.06);
				\draw (-0.035/3+1/16,0.035/2)  node [above,red] { \scriptsize $1$};
			\begin{scope}[shift={(2/9,-3/9)}]
				\draw [thick,black] plot  coordinates {(3/32,2/32) (-3/32,-2/32) };
			\end{scope}
			\begin{scope}[shift={(4/9,-6/9)}]
				\draw [thick,black] plot  coordinates {(3/32,2/32) (-3/32,-2/32) };
			\end{scope}
		\end{scope}
\end{tikzpicture}}

\def\FigTwoCaseDxviii{\begin{tikzpicture}[baseline=15pt]
		\draw [thick,black] plot  coordinates {(-14/12,-14/8) (1/6,1/4) };
		\begin{scope} [shift={(-4/9,-6/9 )}, xscale=-1, yscale=-1]
			\draw [thick,black] plot  coordinates {(8/12,-8/8) (-1/6,1/4) };
			\fill [black] (0,0) circle (0.06);
												\draw (0,0+1/16)  node [below,red] { $\frac{\gamma}{2}$};
			\begin{scope}[shift={(2/9,-3/9)}]
				\draw [thick,black] plot  coordinates {(3/32,2/32) (-3/32,-2/32) };
			\end{scope}
			\begin{scope}[shift={(4/9,-6/9)}]
				\draw [thick,black] plot  coordinates {(3/32,2/32) (-3/32,-2/32) };
			\end{scope}
		\end{scope}
		\begin{scope} [shift={(-7/9,-10.5/9 )},xscale=-1, yscale=-1]
			\draw [thick,black] plot  coordinates {(8/12,-8/8) (-1/6,1/4) };
			\fill [black] (0,0) circle (0.06);
												\draw (0,0+1/16)  node [below,red] { $\frac{\gamma}{2}$};
			\begin{scope}[shift={(2/9,-3/9)}]
				\draw [thick,black] plot  coordinates {(3/32,2/32) (-3/32,-2/32) };
			\end{scope}
			\begin{scope}[shift={(4/9,-6/9)}]
				\draw [thick,black] plot  coordinates {(3/32,2/32) (-3/32,-2/32) };
			\end{scope}
		\end{scope}
		\draw [thick,cborange] plot  coordinates {(-1/2,0) (5/2,0) };
		\fill [black] (0,0) circle (0.06);
		\draw (0+3/16,0+0/16)  node [below,red] { $\frac{\alpha \text{-}\gamma}{8}$};
		\begin{scope}[shift={(4/2,0)}]
			\draw [thick,black] plot  coordinates {(8/12,-8/8) (-1/6,1/4) };
			\fill [black] (0,0) circle (0.06);
					\draw (0-3/16,0+0/16)  node [below,red] { $\frac{3\alpha \texttt{+}\gamma}{8}$};
			\begin{scope}[shift={(2/9,-3/9)}]
				\draw [thick,black] plot  coordinates {(3/32,2/32) (-3/32,-2/32) };
			\end{scope}
			\begin{scope}[shift={(4/9,-6/9)}]
				\draw [thick,black] plot  coordinates {(3/32,2/32) (-3/32,-2/32) };
			\end{scope}
		\end{scope}
		\begin{scope}[shift={(2/3,0)}]
			\draw [thick,cbblue] plot  coordinates {(0,1/3) (0,-6/4) };
			\fill [black] (0,0) circle (0.06);
			\draw (0-3/16,0+1/16)  node [above,red] { $\frac{8\text{-}3\alpha \text{-}\gamma}{24}$};
			\draw (0,-6/4) node [left] {\hbox{\footnotesize$\geqone$\kern-1pt}};
		\end{scope}
		\begin{scope}[shift={(4/3,0)}]
			\draw [thick,cbblue] plot  coordinates {(0,1/3) (0,-6/4) };
			\fill [black] (0,0) circle (0.06);
			\draw (0+3/16,0+1/16)  node [above,red] { $\frac{8\text{-}3\alpha \text{-}\gamma}{24}$};
			\draw (0,-6/4) node [right] {\hbox{\footnotesize$\geqone$}};
		\end{scope}
\end{tikzpicture}}

\def\FigTwoCaseDxix{\begin{tikzpicture}[baseline=15pt]
		\draw [thick,black] plot  coordinates {(-14/12,-14/8) (1/6,1/4) };
\begin{scope} [shift={(-4/9,-6/9 )}, xscale=-1, yscale=-1]
	\draw [thick,black] plot  coordinates {(8/12,-8/8) (-1/6,1/4) };
	\fill [black] (0,0) circle (0.06);
	\draw (0-1/16,0+0/16)  node [right,red] { $\frac{\gamma}{2}$};
	\begin{scope}[shift={(2/9,-3/9)}]
		\draw [thick,black] plot  coordinates {(3/32,2/32) (-3/32,-2/32) };
	\end{scope}
	\begin{scope}[shift={(4/9,-6/9)}]
		\draw [thick,black] plot  coordinates {(3/32,2/32) (-3/32,-2/32) };
	\end{scope}
\end{scope}
\begin{scope} [shift={(-7/9,-10.5/9 )},xscale=-1, yscale=-1]
	\draw [thick,black] plot  coordinates {(8/12,-8/8) (-1/6,1/4) };
	\fill [black] (0,0) circle (0.06);
	\draw (0,0+1/16)  node [below,red] { $\frac{\gamma}{2}$};
	\begin{scope}[shift={(2/9,-3/9)}]
		\draw [thick,black] plot  coordinates {(3/32,2/32) (-3/32,-2/32) };
	\end{scope}
	\begin{scope}[shift={(4/9,-6/9)}]
		\draw [thick,black] plot  coordinates {(3/32,2/32) (-3/32,-2/32) };
	\end{scope}
\end{scope}
\draw [thick,cborange] plot  coordinates {(-1/2,0) (5/2,0) };
\fill [black] (0,0) circle (0.06);
\draw (0-3/16,0+0/16)  node [above,red] { $\frac{\alpha \text{-}\gamma}{8}$};
\begin{scope}[shift={(4/2,0)}]
	\draw [thick,black] plot  coordinates {(8/12,-8/8) (-1/6,1/4) };
	\fill [black] (0,0) circle (0.06);
	\draw (0+3/16,0+0/16)  node [above,red] { $\frac{3\alpha \texttt{+}\gamma}{8}$};
	\begin{scope}[shift={(2/9,-3/9)}]
		\draw [thick,black] plot  coordinates {(3/32,2/32) (-3/32,-2/32) };
	\end{scope}
	\begin{scope}[shift={(4/9,-6/9)}]
		\draw [thick,black] plot  coordinates {(3/32,2/32) (-3/32,-2/32) };
	\end{scope}
\end{scope}
		\begin{scope}[shift={(1/2,0)}]
			\draw [thick,cbgreen] plot  coordinates {(0,1/3) (0,-8/4) };
			\fill [black] (0,0) circle (0.06);
				\draw (0+3/16,0+0/16)  node [above,red] { $\frac{2\delta \text{-}\alpha \text{-}\gamma}{8}$};
		\end{scope}
		\begin{scope}[shift={(1/2,-5/6)}]
			\draw [thick,cbblue] plot  coordinates {(-1/3,0) (4/3,0)};
			\fill [black] (0,0) circle (0.06);
							\draw (0+12/16,0-1/16)  node [above,red] { $\frac{4 \text{-}5\delta  \texttt{+}2\gamma\texttt{+}\alpha}{12}$};
			\draw (2/2-24/16,0) node [below] {\hbox{\footnotesize$\geqone$\kern-1pt}};
		\end{scope}
		\begin{scope}[shift={(1/2,-5/3)}]
			\draw [thick,cbblue] plot  coordinates {(-1/3,0) (4/3,0)};
			\fill [black] (0,0) circle (0.06);
				\draw (0+12/16,0-1/16)  node [above,red] { $\frac{4 \text{-}5\delta  \texttt{+}2\gamma\texttt{+}\alpha}{12}$};
			\draw (2/2,0) node [below] {\hbox{\footnotesize$\geqone$\kern-1pt}};
		\end{scope}
\end{tikzpicture}}

\def\FigTwoCaseDxx{\begin{tikzpicture}[baseline=15pt]
			\draw [thick,black] plot  coordinates {(-14/12,-14/8) (1/6,1/4) };
		\begin{scope} [shift={(-4/9,-6/9 )}, xscale=-1, yscale=-1]
			\draw [thick,black] plot  coordinates {(8/12,-8/8) (-1/6,1/4) };
			\fill [black] (0,0) circle (0.06);
			\draw (0,0+1/16)  node [below,red] { $\frac{\gamma}{2}$};
			\begin{scope}[shift={(2/9,-3/9)}]
				\draw [thick,black] plot  coordinates {(3/32,2/32) (-3/32,-2/32) };
			\end{scope}
			\begin{scope}[shift={(4/9,-6/9)}]
				\draw [thick,black] plot  coordinates {(3/32,2/32) (-3/32,-2/32) };
			\end{scope}
		\end{scope}
		\begin{scope} [shift={(-7/9,-10.5/9 )},xscale=-1, yscale=-1]
			\draw [thick,black] plot  coordinates {(8/12,-8/8) (-1/6,1/4) };
			\fill [black] (0,0) circle (0.06);
			\draw (0,0+1/16)  node [below,red] { $\frac{\gamma}{2}$};
			\begin{scope}[shift={(2/9,-3/9)}]
				\draw [thick,black] plot  coordinates {(3/32,2/32) (-3/32,-2/32) };
			\end{scope}
			\begin{scope}[shift={(4/9,-6/9)}]
				\draw [thick,black] plot  coordinates {(3/32,2/32) (-3/32,-2/32) };
			\end{scope}
		\end{scope}
		\draw [thick,cborange] plot  coordinates {(-1/2,0) (5/2,0) };
		\fill [black] (0,0) circle (0.06);
		\draw (0-3/16,0+0/16)  node [above,red] { $\frac{\alpha \text{-}\gamma}{8}$};
		\begin{scope}[shift={(4/2,0)}]
			\draw [thick,black] plot  coordinates {(8/12,-8/8) (-1/6,1/4) };
			\fill [black] (0,0) circle (0.06);
			\draw (0+3/16,0+0/16)  node [above,red] { $\frac{3\alpha \texttt{+}\gamma}{8}$};
			\begin{scope}[shift={(2/9,-3/9)}]
				\draw [thick,black] plot  coordinates {(3/32,2/32) (-3/32,-2/32) };
			\end{scope}
			\begin{scope}[shift={(4/9,-6/9)}]
				\draw [thick,black] plot  coordinates {(3/32,2/32) (-3/32,-2/32) };
			\end{scope}
		\end{scope}
		\begin{scope}[shift={(1,0)}]
			\draw [thick,cbblue] plot  coordinates {(0,1/3) (0,-6/4) };
			\fill [black] (0,0) circle (0.06);
				\draw (0+8/16,0+0/16)  node [below,red] { $\frac{8\text{-}3\alpha \text{-}\gamma}{40}$};
			\draw (0,-6/4) node [right] {\hbox{\footnotesize$\geqtwo$}};
		\end{scope}
\end{tikzpicture}}

\def\FigTwoCaseDxxi{\begin{tikzpicture}[baseline=15pt]
		\draw [thick,black] plot  coordinates {(-14/12,-14/8) (1/6,1/4) };
	\begin{scope} [shift={(-4/9,-6/9 )}, xscale=-1, yscale=-1]
		\draw [thick,black] plot  coordinates {(8/12,-8/8) (-1/6,1/4) };
		\fill [black] (0,0) circle (0.06);
		\draw (0,0+1/16)  node [below,red] { $\frac{\gamma}{2}$};
		\begin{scope}[shift={(2/9,-3/9)}]
			\draw [thick,black] plot  coordinates {(3/32,2/32) (-3/32,-2/32) };
		\end{scope}
		\begin{scope}[shift={(4/9,-6/9)}]
			\draw [thick,black] plot  coordinates {(3/32,2/32) (-3/32,-2/32) };
		\end{scope}
	\end{scope}
	\begin{scope} [shift={(-7/9,-10.5/9 )},xscale=-1, yscale=-1]
		\draw [thick,black] plot  coordinates {(8/12,-8/8) (-1/6,1/4) };
		\fill [black] (0,0) circle (0.06);
		\draw (0,0+1/16)  node [below,red] { $\frac{\gamma}{2}$};
		\begin{scope}[shift={(2/9,-3/9)}]
			\draw [thick,black] plot  coordinates {(3/32,2/32) (-3/32,-2/32) };
		\end{scope}
		\begin{scope}[shift={(4/9,-6/9)}]
			\draw [thick,black] plot  coordinates {(3/32,2/32) (-3/32,-2/32) };
		\end{scope}
	\end{scope}
	\draw [thick,cborange] plot  coordinates {(-1/2,0) (5/2,0) };
	\fill [black] (0,0) circle (0.06);
	\draw (0-3/16,0+0/16)  node [above,red] { $\frac{\alpha \text{-}\gamma}{8}$};
	\begin{scope}[shift={(4/2,0)}]
		\draw [thick,black] plot  coordinates {(8/12,-8/8) (-1/6,1/4) };
		\fill [black] (0,0) circle (0.06);
		\draw (0+3/16,0+0/16)  node [above,red] { $\frac{3\alpha \texttt{+}\gamma}{8}$};
		\begin{scope}[shift={(2/9,-3/9)}]
			\draw [thick,black] plot  coordinates {(3/32,2/32) (-3/32,-2/32) };
		\end{scope}
		\begin{scope}[shift={(4/9,-6/9)}]
			\draw [thick,black] plot  coordinates {(3/32,2/32) (-3/32,-2/32) };
		\end{scope}
	\end{scope}
		\draw (1,0) node [below] {\hbox{\footnotesize$\geqtwo$}};
\end{tikzpicture}}

\def\FigTwoCaseDxxii{\begin{tikzpicture}[baseline=15pt]
		\begin{scope} [shift={(-3/4,0 )}, xscale=-1, yscale=1]
			\draw [thick,black] plot  coordinates {(8/12,-8/8) (-1/6,1/4) };
			\fill [black] (0,0) circle (0.06);
						\draw (0+1/16,0+1/16)  node [above,red] { $\frac{\alpha}{2}$};
			\begin{scope}[shift={(2/9,-3/9)}]
				\draw [thick,black] plot  coordinates {(3/32,2/32) (-3/32,-2/32) };
			\end{scope}
			\begin{scope}[shift={(4/9,-6/9)}]
				\draw [thick,black] plot  coordinates {(3/32,2/32) (-3/32,-2/32) };
			\end{scope}
		\end{scope}
		\begin{scope} [shift={(0,0 )}, xscale=-1, yscale=1]
		\draw [thick,black] plot  coordinates {(8/12,-8/8) (-1/6,1/4) };
		\fill [black] (0,0) circle (0.06);
		\draw (0+1/16,0+1/16)  node [above,red] { $\frac{\alpha}{2}$};
		\begin{scope}[shift={(2/9,-3/9)}]
			\draw [thick,black] plot  coordinates {(3/32,2/32) (-3/32,-2/32) };
		\end{scope}
		\begin{scope}[shift={(4/9,-6/9)}]
			\draw [thick,black] plot  coordinates {(3/32,2/32) (-3/32,-2/32) };
		\end{scope}
	\end{scope}
		\draw [thick,black] plot  coordinates {(-3/2+1/4,0) (5/2,0) };
		\begin{scope}[shift={(4/2,0)}]
			\draw [thick,black] plot  coordinates {(8/12,-8/8) (-1/6,1/4) };
			\fill [black] (0,0) circle (0.06);
			\draw (0+1/16,0+1/16)  node [above,red] { $\frac{\alpha}{2}$};
			\begin{scope}[shift={(2/9,-3/9)}]
				\draw [thick,black] plot  coordinates {(3/32,2/32) (-3/32,-2/32) };
			\end{scope}
			\begin{scope}[shift={(4/9,-6/9)}]
				\draw [thick,black] plot  coordinates {(3/32,2/32) (-3/32,-2/32) };
			\end{scope}
		\end{scope}
		\begin{scope}[shift={(2/3,0)}]
			\draw [thick,cbblue] plot  coordinates {(0,1/3) (0,-6/4) };
			\fill [black] (0,0) circle (0.06);
			\draw (0-5/16,0+0/16)  node [below,red] { $\frac{2\text{-}\alpha}{6}$};
			\draw (0,-6/4) node [left] {\hbox{\footnotesize$\geqone$\kern-1pt}};
		\end{scope}
		\begin{scope}[shift={(4/3,0)}]
			\draw [thick,cbblue] plot  coordinates {(0,1/3) (0,-6/4) };
			\fill [black] (0,0) circle (0.06);
						\draw (0-5/16,0+0/16)  node [below,red] { $\frac{2\text{-}\alpha}{6}$};
			\draw (0,-6/4) node [right] {\hbox{\footnotesize$\geqone$}};
		\end{scope}
\end{tikzpicture}}

\def\FigTwoCaseDxxiii{\begin{tikzpicture}[baseline=15pt]
	\begin{scope} [shift={(-3/4,0 )}, xscale=-1, yscale=1]
		\draw [thick,black] plot  coordinates {(8/12,-8/8) (-1/6,1/4) };
		\fill [black] (0,0) circle (0.06);
		\draw (0+1/16,0+1/16)  node [above,red] { $\frac{\alpha}{2}$};
		\begin{scope}[shift={(2/9,-3/9)}]
			\draw [thick,black] plot  coordinates {(3/32,2/32) (-3/32,-2/32) };
		\end{scope}
		\begin{scope}[shift={(4/9,-6/9)}]
			\draw [thick,black] plot  coordinates {(3/32,2/32) (-3/32,-2/32) };
		\end{scope}
	\end{scope}
	\begin{scope} [shift={(0,0 )}, xscale=-1, yscale=1]
		\draw [thick,black] plot  coordinates {(8/12,-8/8) (-1/6,1/4) };
		\fill [black] (0,0) circle (0.06);
		\draw (0+1/16,0+1/16)  node [above,red] { $\frac{\alpha}{2}$};
		\begin{scope}[shift={(2/9,-3/9)}]
			\draw [thick,black] plot  coordinates {(3/32,2/32) (-3/32,-2/32) };
		\end{scope}
		\begin{scope}[shift={(4/9,-6/9)}]
			\draw [thick,black] plot  coordinates {(3/32,2/32) (-3/32,-2/32) };
		\end{scope}
	\end{scope}
	\draw [thick,black] plot  coordinates {(-3/2+1/4,0) (5/2,0) };
	\begin{scope}[shift={(4/2,0)}]
		\draw [thick,black] plot  coordinates {(8/12,-8/8) (-1/6,1/4) };
		\fill [black] (0,0) circle (0.06);
		\draw (0+1/16,0+1/16)  node [above,red] { $\frac{\alpha}{2}$};
		\begin{scope}[shift={(2/9,-3/9)}]
			\draw [thick,black] plot  coordinates {(3/32,2/32) (-3/32,-2/32) };
		\end{scope}
		\begin{scope}[shift={(4/9,-6/9)}]
			\draw [thick,black] plot  coordinates {(3/32,2/32) (-3/32,-2/32) };
		\end{scope}
	\end{scope}
\begin{scope}[shift={(1/2,0)}]
	\draw [thick,cbgreen] plot  coordinates {(0,1/3) (0,-8/4) };
	\fill [black] (0,0) circle (0.06);
		\draw (0+3/16,0+1/16)  node [above,red] { $\frac{\delta\text{-}\alpha}{4}$};
\end{scope}
\begin{scope}[shift={(1/2,-5/6)}]
	\draw [thick,cbblue] plot  coordinates {(-1/3,0) (4/3,0)};
	\fill [black] (0,0) circle (0.06);
								\draw (0+10/16,0-1/16)  node [above,red] { $\frac{4 \text{-}5\delta  \texttt{+}3\alpha}{12}$};
	\draw (2/2-24/16,0) node [below] {\hbox{\footnotesize$\geqone$\kern-1pt}};
\end{scope}
\begin{scope}[shift={(1/2,-5/3)}]
	\draw [thick,cbblue] plot  coordinates {(-1/3,0) (4/3,0)};
	\fill [black] (0,0) circle (0.06);
								\draw (0+10/16,0-1/16)  node [above,red] { $\frac{4 \text{-}5\delta  \texttt{+}3\alpha}{12}$};
	\draw (2/2,0) node [below] {\hbox{\footnotesize$\geqone$\kern-1pt}};
\end{scope}
\end{tikzpicture}}

\def\FigTwoCaseDxxiv{\begin{tikzpicture}[baseline=15pt]
	\begin{scope} [shift={(-3/4,0 )}, xscale=-1, yscale=1]
		\draw [thick,black] plot  coordinates {(8/12,-8/8) (-1/6,1/4) };
		\fill [black] (0,0) circle (0.06);
		\draw (0+1/16,0+1/16)  node [above,red] { $\frac{\alpha}{2}$};
		\begin{scope}[shift={(2/9,-3/9)}]
			\draw [thick,black] plot  coordinates {(3/32,2/32) (-3/32,-2/32) };
		\end{scope}
		\begin{scope}[shift={(4/9,-6/9)}]
			\draw [thick,black] plot  coordinates {(3/32,2/32) (-3/32,-2/32) };
		\end{scope}
	\end{scope}
	\begin{scope} [shift={(0,0 )}, xscale=-1, yscale=1]
		\draw [thick,black] plot  coordinates {(8/12,-8/8) (-1/6,1/4) };
		\fill [black] (0,0) circle (0.06);
		\draw (0+1/16,0+1/16)  node [above,red] { $\frac{\alpha}{2}$};
		\begin{scope}[shift={(2/9,-3/9)}]
			\draw [thick,black] plot  coordinates {(3/32,2/32) (-3/32,-2/32) };
		\end{scope}
		\begin{scope}[shift={(4/9,-6/9)}]
			\draw [thick,black] plot  coordinates {(3/32,2/32) (-3/32,-2/32) };
		\end{scope}
	\end{scope}
	\draw [thick,black] plot  coordinates {(-3/2+1/4,0) (5/2,0) };
	\begin{scope}[shift={(4/2,0)}]
		\draw [thick,black] plot  coordinates {(8/12,-8/8) (-1/6,1/4) };
		\fill [black] (0,0) circle (0.06);
		\draw (0+1/16,0+1/16)  node [above,red] { $\frac{\alpha}{2}$};
		\begin{scope}[shift={(2/9,-3/9)}]
			\draw [thick,black] plot  coordinates {(3/32,2/32) (-3/32,-2/32) };
		\end{scope}
		\begin{scope}[shift={(4/9,-6/9)}]
			\draw [thick,black] plot  coordinates {(3/32,2/32) (-3/32,-2/32) };
		\end{scope}
	\end{scope}
		\begin{scope}[shift={(1,0)}]
			\draw [thick,cbblue] plot  coordinates {(0,1/3) (0,-6/4) };
			\fill [black] (0,0) circle (0.06);
					\draw (0-5/16,0+1/16)  node [below,red] { $\frac{2\text{-}\alpha}{10}$};
			\draw (0,-6/4) node [right] {\hbox{\footnotesize$\geqtwo$}};
		\end{scope}
\end{tikzpicture}}

\def\includefigD#1{
\makebox(130pt,87pt){#1}}

\def\FigTwoCaseDone{
\begin{tabular}{|c|c|c|}
\hline
{\large\strut} (D1) & (D2) & (D3) \\
\hline
${\Large\strut}
\alpha=\beta=\gamma=2$ 
& $\alpha>\beta=\gamma=2$ 
& $\beta>\gamma=2$ 
\\
\hline
\includefigD{\FigTwoCaseDi} & 
\includefigD{\FigTwoCaseDii} & 
\includefigD{\FigTwoCaseDiii}
 \\
\hline\hline
{\large\strut} (D4) & (D5) & (D6) \\
\hline
${\Large\strut}
\gamma >2$ 
& $\beta>\gamma \ \land \ \alpha+\gamma<4$ 
& $\beta>\gamma \ \land \ \beta+\gamma<4=\alpha+\gamma$ 
\\
\hline
\includefigD{\FigTwoCaseDiv} &
\includefigD{\FigTwoCaseDv}&
\includefigD{\FigTwoCaseDvi}
\\
\hline\hline
{\large\strut} (D7) & (D8) & (D9) \\
\hline
${\Large\strut}
\beta>\gamma \ \land \ \beta+\gamma=4=\alpha+\gamma$
& $\beta>\gamma \ \land \ \beta+\gamma<4<\alpha+\gamma$
& $2>\gamma \ \land \ \beta+\gamma=4<\alpha+\gamma$ 
\\
\hline
\includefigD{\FigTwoCaseDvii} &
\includefigD{\FigTwoCaseDviii}&
\includefigD{\FigTwoCaseDix}
\\
\hline\hline
{\large\strut} (D10) & (D11) & (D12) \\
\hline
${\Large\strut}
2>\gamma<\beta \ \land \ \beta+\gamma>4$ 
& $\beta=\gamma=\delta<\alpha \ \land \ \alpha+\gamma<4$ 
& $\beta=\gamma=\delta<\alpha \ \land \ \alpha+\gamma=4$ 
\\
\hline
\includefigD{\FigTwoCaseDx} &
\includefigD{\FigTwoCaseDxi} &
\includefigD{\FigTwoCaseDxii}
\\
\hline 

\end{tabular}}

\def\FigTwoCaseDtwo{
\begin{tabular}{|c|c|c|}
	\hline
	{\large\strut} (D13) & (D14) & (D15) \\
	\hline
	${\Large\strut}
	\beta=\gamma=\delta<2 \ \land \ \alpha+\gamma>4$ 
	& $ \gamma <\delta<\frac{\alpha+\gamma}{2} \ \land \ \alpha +\delta<4
	 $ 
	& $ \gamma <\delta<\frac{\alpha+\gamma}{2} \ \land \ \alpha +\delta=4
	$ 
	\\
	\hline
	\includefigD{\FigTwoCaseDxiii} & 
	\includefigD{\FigTwoCaseDxiv} & 
	\includefigD{\FigTwoCaseDxv}
	\\
	\hline\hline
	{\large\strut} (D16) & (D17) & (D18) \\
	\hline
	${\Large\strut}
	\gamma<\delta<\frac{4+\gamma}{3} \ \land \ \alpha+\delta>4
	$ 
	& $\delta\geq \frac{4+\gamma}{3} \ \land \ \alpha>\frac{8-\gamma}{3}>2  
	$ 
	& $  3\alpha+\gamma < 8 \ \land \  \delta = \frac{\alpha+\gamma}{2}>\gamma$
	\\
	\hline
	\includefigD{\FigTwoCaseDxvi} & 
	\includefigD{\FigTwoCaseDxvii} & 
	\includefigD{\FigTwoCaseDxviii}
	\\
	\hline\hline
		{\large\strut} (D19) & (D20) & (D21) \\
	\hline
	{\Large\strut}
	$\gamma<\alpha \ \land \   \frac{\alpha+\gamma}{2}<\delta<\frac{\alpha+2\gamma+4}{5} $ 
	& $\gamma< \frac{3\alpha+\gamma}{4} < 2 \ \land \  \delta\geq \frac{\alpha+2\gamma+4}{5} $
	& $  3\alpha+\gamma = 8 \ \land \   \delta\geq\frac{\alpha+\gamma}{2}>\gamma $
	\\
	\hline
	\includefigD{\FigTwoCaseDxix} & 
	\includefigD{\FigTwoCaseDxx} & 
	\includefigD{\FigTwoCaseDxxi}
	\\
	\hline\hline
		{\large\strut} (D22) & (D23) & (D24) \\
	\hline
	{\Large\strut}
	$\alpha=\gamma=\delta<2 $ 
	& $2>\alpha=\gamma<\delta<\frac{\alpha+2\gamma +4}{5}$
	&$2>\alpha=\gamma \ \land \ \delta\geq\frac{\alpha+2\gamma +4}{5}$
	\\
	\hline
	\includefigD{\FigTwoCaseDxxii} & 
	\includefigD{\FigTwoCaseDxxiii} & 
	\includefigD{\FigTwoCaseDxxiv}
	\\
	\hline
\end{tabular}}

\def\FigTwoCaseEFGi{\begin{tikzpicture}[baseline=12pt]
\draw (7/4,0/4) node {};
\draw [thick,black] plot  coordinates {(13/8,-13/8) (-1/4,1/4) };
\begin{scope}[shift={(4/3,-4/3)}]
	\draw [thick,cbblue] plot  coordinates {(8/12,8/8) (-1/6,-1/4) };
	\fill [black] (0,0) circle (0.06);
			\draw (0+1/16,0+1/32)  node [right,red] { $\frac{1}{3}$};
	\draw (2/4,4/4) node [above] {\hbox{\footnotesize$\geqone$}};
\end{scope}
\begin{scope}[xscale=-1]
\draw [thick,black] plot  coordinates {(13/8,-13/8) (-1/4,1/4) };
\begin{scope}[shift={(1/3,-1/3)}]
	\draw [thick,black] plot  coordinates {(-1/12,-1/12) (1/12,1/12) };
\end{scope}
\begin{scope}[shift={(2/3,-2/3)}]
	\draw [thick,black] plot  coordinates {(-1/12,-1/12) (1/12,1/12) };
\end{scope}
\begin{scope}[shift={(3/3,-3/3)}]
	\draw [thick,black] plot  coordinates {(-1/12,-1/12) (1/12,1/12) };
\end{scope}

\begin{scope}[shift={(4/3,-4/3)}]
	\draw [thick,cbblue] plot  coordinates {(8/12,8/8) (-1/6,-1/4) };
	\fill [black] (0,0) circle (0.06);
		\draw (0+1/16,0+1/32)  node [left,red] { $\frac{1}{3}$};
	\draw (2/4,4/4) node [above] {\hbox{\footnotesize$\geqone$}};
\end{scope}
\end{scope}

\begin{scope}[shift={(0,0)}]
\fill [white] (0,0) circle (0.07);
\draw [thick] (0,0) circle (0.07);
		\draw (0+0/16,0+1/16)  node [above,red] { $\frac{\epsilon}{2}$};
\end{scope}

\begin{scope}[shift={(1/3,-1/3)}]
\draw [thick,black] plot  coordinates {(-1/12,-1/12) (1/12,1/12) };
\end{scope}
\begin{scope}[shift={(2/3,-2/3)}]
\draw [thick,black] plot  coordinates {(-1/12,-1/12) (1/12,1/12) };
\end{scope}
\begin{scope}[shift={(3/3,-3/3)}]
\draw [thick,black] plot  coordinates {(-1/12,-1/12) (1/12,1/12) };
\end{scope}
\end{tikzpicture}}

\def\FigTwoCaseEFGii{\begin{tikzpicture}[baseline=12pt]
		\draw (7/4,0/4) node {};
\draw [thick,black] plot  coordinates {(3/3,-3/3) (-1/4,1/4) };
\begin{scope}[shift={(-1/3,+1/3)}]
\draw [thick,cborange] plot  coordinates {(2/3,-1) (8/3,-1) };
\fill [black] (+3/3,-1) circle (0.06);
		\draw (+3/3,-1)  node [below,red] { $\frac{\alpha}{4}$};
\begin{scope}[shift={(-1/3,0)}]
\draw [thick,cbblue] plot  coordinates {(2,-4/3) (2,1/8) };
\fill [black] (2,-1) circle (0.06);
		\draw (2+1/16,-1+6/16)  node [left,red] { $\frac{4\text{-}\alpha}{12}$};
\draw (2,0) node [right] {\hbox{\footnotesize$\geqone$}};
\end{scope}
\end{scope}

\begin{scope} [shift={(2,-2/3 )}, xscale=-1, yscale=1]
	\draw [thick,black] plot  coordinates {(0,-8/8) (0,1/4) };
	\fill [black] (0,0) circle (0.06);
			\draw (0-3/16,0)  node [above,red] { $\frac{\alpha}{4}$};
	\begin{scope}[shift={(0,-3/9)}]
		\draw [thick,black] plot  coordinates {(3.6/32,0) (-3.6/32,0) };
	\end{scope}
	\begin{scope}[shift={(0,-6/9)}]
		\draw [thick,black] plot  coordinates {(3.6/32,0) (-3.6/32,0) };
	\end{scope}
\end{scope}

		\begin{scope}[xscale=-1]
			\draw [thick,black] plot  coordinates {(13/8,-13/8) (-1/4,1/4) };
			\begin{scope}[shift={(1/3,-1/3)}]
				\draw [thick,black] plot  coordinates {(-1/12,-1/12) (1/12,1/12) };
			\end{scope}
			\begin{scope}[shift={(2/3,-2/3)}]
				\draw [thick,black] plot  coordinates {(-1/12,-1/12) (1/12,1/12) };
			\end{scope}
			\begin{scope}[shift={(3/3,-3/3)}]
				\draw [thick,black] plot  coordinates {(-1/12,-1/12) (1/12,1/12) };
			\end{scope}
			
			\begin{scope}[shift={(4/3,-4/3)}]
				\draw [thick,cbblue] plot  coordinates {(8/12,8/8) (-1/6,-1/4) };
				\fill [black] (0,0) circle (0.06);
					\draw (0+1/16,0+1/32)  node [left,red] { $\frac{1}{3}$};
				\draw (2/4,4/4) node [above] {\hbox{\footnotesize$\geqone$}};
			\end{scope}
		\end{scope}
		
		\begin{scope}[shift={(0,0)}]
			\fill [white] (0,0) circle (0.07);
			\draw [thick] (0,0) circle (0.07);
					\draw (0+0/16,0+1/16)  node [above,red] { $\frac{\epsilon}{2}$};
		\end{scope}
		
		\begin{scope}[shift={(1/3,-1/3)}]
			\draw [thick,black] plot  coordinates {(-1/12,-1/12) (1/12,1/12) };
		\end{scope}
\end{tikzpicture}}

\def\FigTwoCaseEFGiii{\begin{tikzpicture}[baseline=12pt]
		\draw (7/4,0/4) node {};
	\draw [thick,black] plot  coordinates {(3/3,-3/3) (-1/4,1/4) };
	\begin{scope}[shift={(-1/3,+1/3)}]
		\draw [thick,cborange] plot  coordinates {(2/3,-1) (8/3,-1) };
		\fill [black] (+3/3,-1) circle (0.06);
		\draw (+3/3,-1)  node [below,red] { $\frac{\alpha}{4}$};
		\begin{scope}[shift={(-1/3,0)}]
			\draw [thick,cbblue] plot  coordinates {(2,-4/3) (2,1/8) };
			\fill [black] (2,-1) circle (0.06);
			\draw (2+1/16,-1+6/16)  node [left,red] { $\frac{4\text{-}\alpha}{12}$};
			\draw (2,0) node [right] {\hbox{\footnotesize$\geqone$}};
		\end{scope}
	\end{scope}
		
	\begin{scope} [shift={(2,-2/3 )}, xscale=-1, yscale=1]
		\draw [thick,black] plot  coordinates {(0,-8/8) (0,1/4) };
		\fill [black] (0,0) circle (0.06);
			\draw(0-3/16,0)  node [above,red] { $\frac{\alpha}{4}$};
		\begin{scope}[shift={(0,-3/9)}]
			\draw [thick,black] plot  coordinates {(3.6/32,0) (-3.6/32,0) };
		\end{scope}
		\begin{scope}[shift={(0,-6/9)}]
			\draw [thick,black] plot  coordinates {(3.6/32,0) (-3.6/32,0) };
		\end{scope}
	\end{scope}
		
		\begin{scope}[xscale=-1]
					\draw [thick,black] plot  coordinates {(3/3,-3/3) (-1/4,1/4) };
				\begin{scope}[shift={(-1/3,+1/3)}]
					\draw [thick,cborange] plot  coordinates {(2/3,-1) (8/3,-1) };
					\fill [black] (+3/3,-1) circle (0.06);
					\draw (+3/3,-1)  node [below,red] { $\frac{\beta}{4}$};
					\begin{scope}[shift={(-1/3,0)}]
						\draw [thick,cbblue] plot  coordinates {(2,-4/3) (2,1/8) };
						\fill [black] (2,-1) circle (0.06);
						\draw (2+1/16,-1+6/16)  node [right,red] { $\frac{4\text{-}\beta}{12}$};
						\draw (2,0) node [left] {\hbox{\footnotesize$\geqone$}};
					\end{scope}
				\end{scope}
				
				\begin{scope} [shift={(2,-2/3 )}, xscale=-1, yscale=1]
					\draw [thick,black] plot  coordinates {(0,-8/8) (0,1/4) };
					\fill [black] (0,0) circle (0.06);
					\draw(0-3/16,0)  node [above,red] { $\frac{\beta}{4}$};
					\begin{scope}[shift={(0,-3/9)}]
						\draw [thick,black] plot  coordinates {(3.6/32,0) (-3.6/32,0) };
					\end{scope}
					\begin{scope}[shift={(0,-6/9)}]
						\draw [thick,black] plot  coordinates {(3.6/32,0) (-3.6/32,0) };
					\end{scope}
				\end{scope}
		\end{scope}
		
		\begin{scope}[shift={(0,0)}]
			\fill [white] (0,0) circle (0.07);
			\draw [thick] (0,0) circle (0.07);
					\draw (0+0/16,0+1/16)  node [above,red] { $\frac{\epsilon}{2}$};
		\end{scope}
		
		\begin{scope}[shift={(1/3,-1/3)}]
			\draw [thick,black] plot  coordinates {(-1/12,-1/12) (1/12,1/12) };
		\end{scope}
		\begin{scope}[xscale=-1]
			\begin{scope}[shift={(1/3,-1/3)}]
				\draw [thick,black] plot  coordinates {(-1/12,-1/12) (1/12,1/12) };
			\end{scope}
			\end{scope}
\end{tikzpicture}}

\def\FigTwoCaseEFGiv{\begin{tikzpicture}[baseline=12pt]
		\draw (7/4,0/4) node {};
		\draw [thick,black] plot  coordinates {(3/3,-3/3) (-1/4,1/4) };
		\begin{scope}[shift={(-1/3,+1/3)}]
			\draw [thick,cborange] plot  coordinates {(2/3,-1) (8/3,-1) };
			\fill [black] (+3/3,-1) circle (0.06);
				\draw (+3/3,-1)   node [below,red] {\scriptsize $1$};
			\begin{scope}[shift={(-1/3,0)}]
				\draw (2,-1) node [above] {\hbox{\footnotesize$\geqone$}};
			\end{scope}
		\end{scope}
		
\begin{scope} [shift={(2,-2/3 )}, xscale=-1, yscale=1]
	\draw [thick,black] plot  coordinates {(0,-8/8) (0,1/4) };
	\fill [black] (0,0) circle (0.06);
				\draw (0-3/16,0)  node [above,red] {\scriptsize $1$};
	\begin{scope}[shift={(0,-3/9)}]
		\draw [thick,black] plot  coordinates {(3.6/32,0) (-3.6/32,0) };
	\end{scope}
	\begin{scope}[shift={(0,-6/9)}]
		\draw [thick,black] plot  coordinates {(3.6/32,0) (-3.6/32,0) };
	\end{scope}
\end{scope}
		
		\begin{scope}[xscale=-1]
			\draw [thick,black] plot  coordinates {(13/8,-13/8) (-1/4,1/4) };
			\begin{scope}[shift={(1/3,-1/3)}]
				\draw [thick,black] plot  coordinates {(-1/12,-1/12) (1/12,1/12) };
			\end{scope}
			\begin{scope}[shift={(2/3,-2/3)}]
				\draw [thick,black] plot  coordinates {(-1/12,-1/12) (1/12,1/12) };
			\end{scope}
			\begin{scope}[shift={(3/3,-3/3)}]
				\draw [thick,black] plot  coordinates {(-1/12,-1/12) (1/12,1/12) };
			\end{scope}
			
			\begin{scope}[shift={(4/3,-4/3)}]
				\draw [thick,cbblue] plot  coordinates {(8/12,8/8) (-1/6,-1/4) };
				\fill [black] (0,0) circle (0.06);
					\draw (0+1/16,0+1/32)  node [left,red] { $\frac{1}{3}$};
				\draw (2/4,4/4) node [above] {\hbox{\footnotesize$\geqone$}};
			\end{scope}
		\end{scope}
		
		\begin{scope}[shift={(0,0)}]
			\fill [white] (0,0) circle (0.07);
			\draw [thick] (0,0) circle (0.07);
					\draw (0+0/16,0+1/16)  node [above,red] { $\frac{\epsilon}{2}$};
		\end{scope}
		
		\begin{scope}[shift={(1/3,-1/3)}]
			\draw [thick,black] plot  coordinates {(-1/12,-1/12) (1/12,1/12) };
		\end{scope}
\end{tikzpicture}}

\def\FigTwoCaseEFGv{\begin{tikzpicture}[baseline=12pt]
		\draw (7/4,0/4) node {};
		\draw [thick,black] plot  coordinates {(3/3,-3/3) (-1/4,1/4) };
		\begin{scope}[shift={(-1/3,+1/3)}]
			\draw [thick,cborange] plot  coordinates {(2/3,-1) (8/3,-1) };
			\fill [black] (+3/3,-1) circle (0.06);
			\draw (+3/3,-1)   node [below,red] {\scriptsize $1$};
			\begin{scope}[shift={(-1/3,0)}]
				\draw (2,-1) node [above] {\hbox{\footnotesize$\geqone$}};
			\end{scope}
		\end{scope}
		
		\begin{scope} [shift={(2,-2/3 )}, xscale=-1, yscale=1]
			\draw [thick,black] plot  coordinates {(0,-8/8) (0,1/4) };
			\fill [black] (0,0) circle (0.06);
			\draw (0-3/16,0)  node [above,red] {\scriptsize $1$};
			\begin{scope}[shift={(0,-3/9)}]
				\draw [thick,black] plot  coordinates {(3.6/32,0) (-3.6/32,0) };
			\end{scope}
			\begin{scope}[shift={(0,-6/9)}]
				\draw [thick,black] plot  coordinates {(3.6/32,0) (-3.6/32,0) };
			\end{scope}
		\end{scope}
		
	\begin{scope}[xscale=-1]
		\draw [thick,black] plot  coordinates {(3/3,-3/3) (-1/4,1/4) };
		\begin{scope}[shift={(-1/3,+1/3)}]
			\draw [thick,cborange] plot  coordinates {(2/3,-1) (8/3,-1) };
			\fill [black] (+3/3,-1) circle (0.06);
			\draw (+3/3,-1)  node [below,red] { $\frac{\beta}{4}$};
			\begin{scope}[shift={(-1/3,0)}]
				\draw [thick,cbblue] plot  coordinates {(2,-4/3) (2,1/8) };
				\fill [black] (2,-1) circle (0.06);
				\draw (2+1/16,-1+6/16)  node [right,red] { $\frac{4\text{-}\beta}{12}$};
				\draw (2,0) node [right] {\hbox{\footnotesize$\geqone$}};
			\end{scope}
		\end{scope}
		
		\begin{scope} [shift={(2,-2/3 )}, xscale=-1, yscale=1]
			\draw [thick,black] plot  coordinates {(0,-8/8) (0,1/4) };
			\fill [black] (0,0) circle (0.06);
			\draw(0-3/16,0)  node [above,red] { $\frac{\beta}{4}$};
			\begin{scope}[shift={(0,-3/9)}]
				\draw [thick,black] plot  coordinates {(3.6/32,0) (-3.6/32,0) };
			\end{scope}
			\begin{scope}[shift={(0,-6/9)}]
				\draw [thick,black] plot  coordinates {(3.6/32,0) (-3.6/32,0) };
			\end{scope}
		\end{scope}
	\end{scope}
		
		\begin{scope}[shift={(0,0)}]
			\fill [white] (0,0) circle (0.07);
			\draw [thick] (0,0) circle (0.07);
					\draw (0+0/16,0+1/16)  node [above,red] { $\frac{\epsilon}{2}$};
		\end{scope}
		
				\begin{scope}[shift={(1/3,-1/3)}]
			\draw [thick,black] plot  coordinates {(-1/12,-1/12) (1/12,1/12) };
		\end{scope}
		\begin{scope}[xscale=-1]
			\begin{scope}[shift={(1/3,-1/3)}]
				\draw [thick,black] plot  coordinates {(-1/12,-1/12) (1/12,1/12) };
			\end{scope}
		\end{scope}
		
\end{tikzpicture}}

\def\FigTwoCaseEFGvi{\begin{tikzpicture}[baseline=12pt]
		\draw (7/4,0/4) node {};
	\draw [thick,black] plot  coordinates {(3/3,-3/3) (-1/4,1/4) };
	\begin{scope}[shift={(-1/3,+1/3)}]
		\draw [thick,cborange] plot  coordinates {(2/3,-1) (8/3,-1) };
		\fill [black] (+3/3,-1) circle (0.06);
		\draw (+3/3,-1)   node [below,red] {\scriptsize $1$};
		\begin{scope}[shift={(-1/3,0)}]
			\draw (2,-1) node [above] {\hbox{\footnotesize$\geqone$}};
		\end{scope}
	\end{scope}
	
	\begin{scope} [shift={(2,-2/3 )}, xscale=-1, yscale=1]
		\draw [thick,black] plot  coordinates {(0,-8/8) (0,1/4) };
		\fill [black] (0,0) circle (0.06);
		\draw (0-3/16,0)  node [above,red] {\scriptsize $1$};
		\begin{scope}[shift={(0,-3/9)}]
			\draw [thick,black] plot  coordinates {(3.6/32,0) (-3.6/32,0) };
		\end{scope}
		\begin{scope}[shift={(0,-6/9)}]
			\draw [thick,black] plot  coordinates {(3.6/32,0) (-3.6/32,0) };
		\end{scope}
	\end{scope}
		
		\begin{scope}[xscale=-1]
			\draw [thick,black] plot  coordinates {(3/3,-3/3) (-1/4,1/4) };
		\begin{scope}[shift={(-1/3,+1/3)}]
			\draw [thick,cborange] plot  coordinates {(2/3,-1) (8/3,-1) };
			\fill [black] (+3/3,-1) circle (0.06);
			\draw (+3/3,-1)   node [below,red] {\scriptsize $1$};
			\begin{scope}[shift={(-1/3,0)}]
				\draw (2,-1) node [above] {\hbox{\footnotesize$\geqone$}};
			\end{scope}
		\end{scope}
		
		\begin{scope} [shift={(2,-2/3 )}, xscale=-1, yscale=1]
			\draw [thick,black] plot  coordinates {(0,-8/8) (0,1/4) };
			\fill [black] (0,0) circle (0.06);
			\draw (0-3/16,0)  node [above,red] {\scriptsize $1$};
			\begin{scope}[shift={(0,-3/9)}]
				\draw [thick,black] plot  coordinates {(3.6/32,0) (-3.6/32,0) };
			\end{scope}
			\begin{scope}[shift={(0,-6/9)}]
				\draw [thick,black] plot  coordinates {(3.6/32,0) (-3.6/32,0) };
			\end{scope}
		\end{scope}
		\end{scope}
		
		\begin{scope}[shift={(0,0)}]
			\fill [white] (0,0) circle (0.07);
			\draw [thick] (0,0) circle (0.07);
					\draw (0+0/16,0+1/16)  node [above,red] { $\frac{\epsilon}{2}$};
		\end{scope}
		
		\begin{scope}[shift={(1/3,-1/3)}]
	\draw [thick,black] plot  coordinates {(-1/12,-1/12) (1/12,1/12) };
\end{scope}
\begin{scope}[xscale=-1]
	\begin{scope}[shift={(1/3,-1/3)}]
		\draw [thick,black] plot  coordinates {(-1/12,-1/12) (1/12,1/12) };
	\end{scope}
\end{scope}
\end{tikzpicture}}

\def\FigTwoCaseEFGvii{\begin{tikzpicture}[baseline=12pt]
		\draw (7/4,0/4) node {};
		\draw [thick,black] plot  coordinates {(3/3,-3/3) (-1/4,1/4) };
		\begin{scope}[shift={(-1/3,+1/3)}]
			\draw [thick,cborange] plot [smooth, tension=1] coordinates {(2/3,-9/8) (2,-4/8) (4/3,1/8) };
		\draw [thick,cborange] plot [smooth, tension=1] coordinates {(7/2-2/3,-9/8) (7/2-2,-4/8) (7/2-4/3,1/8) };
		\fill [black] (7/4,-0.075) circle (0.06);
		\draw (7/4,-0.075+1/16)  node [above,red] {\scriptsize $\alpha \text{-} 4$};
		\fill [black] (7/4,-0.72) circle (0.06);
			\draw (7/4,-0.72-1/16)  node [below,red] {\scriptsize $\alpha \text{-} 4$};
		\fill [black] (1.025,-1.025) circle (0.06);
		\draw(1.025,-1.025)  node [below,red] {\scriptsize $1$};
		\end{scope}
	
				\begin{scope} [shift={(2+1/32,-2/3 )}, xscale=-1, yscale=1]
			\draw [thick,black] plot  coordinates {(0,-8/8) (0,1/4) };
			\fill [black] (0,0) circle (0.06);
				\draw(0-3/16,0)  node [above,red] {\scriptsize $1$};
			\begin{scope}[shift={(0,-3/9)}]
				\draw [thick,black] plot  coordinates {(3.6/32,0) (-3.6/32,0) };
			\end{scope}
			\begin{scope}[shift={(0,-6/9)}]
				\draw [thick,black] plot  coordinates {(3.6/32,0) (-3.6/32,0) };
			\end{scope}
		\end{scope}
		
		\begin{scope}[xscale=-1]
			\draw [thick,black] plot  coordinates {(13/8,-13/8) (-1/4,1/4) };
			\begin{scope}[shift={(1/3,-1/3)}]
				\draw [thick,black] plot  coordinates {(-1/12,-1/12) (1/12,1/12) };
			\end{scope}
			\begin{scope}[shift={(2/3,-2/3)}]
				\draw [thick,black] plot  coordinates {(-1/12,-1/12) (1/12,1/12) };
			\end{scope}
			\begin{scope}[shift={(3/3,-3/3)}]
				\draw [thick,black] plot  coordinates {(-1/12,-1/12) (1/12,1/12) };
			\end{scope}
			
			\begin{scope}[shift={(4/3,-4/3)}]
				\draw [thick,cbblue] plot  coordinates {(8/12,8/8) (-1/6,-1/4) };
				\fill [black] (0,0) circle (0.06);
				\draw (0+1/16,0+1/32)  node [left,red] { $\frac{1}{3}$};
				\draw (2/4,4/4) node [above] {\hbox{\footnotesize$\geqone$}};
			\end{scope}
		\end{scope}
		
		\begin{scope}[shift={(0,0)}]
			\fill [white] (0,0) circle (0.07);
			\draw [thick] (0,0) circle (0.07);
					\draw (0+0/16,0+1/16)  node [above,red] { $\frac{\epsilon}{2}$};
		\end{scope}
		
		\begin{scope}[shift={(1/3,-1/3)}]
			\draw [thick,black] plot  coordinates {(-1/12,-1/12) (1/12,1/12) };
		\end{scope}
\end{tikzpicture}}

\def\FigTwoCaseEFGviii{\begin{tikzpicture}[baseline=12pt]
		\draw (7/4,0/4) node {};
		\draw [thick,black] plot  coordinates {(3/3,-3/3) (-1/4,1/4) };
		\begin{scope}[shift={(-1/3,+1/3)}]
			\draw [thick,cborange] plot [smooth, tension=1] coordinates {(2/3,-9/8) (2,-4/8) (4/3,1/8) };
			\draw [thick,cborange] plot [smooth, tension=1] coordinates {(7/2-2/3,-9/8) (7/2-2,-4/8) (7/2-4/3,1/8) };
			\fill [black] (7/4,-0.075) circle (0.06);
			\draw (7/4,-0.075+1/16)  node [above,red] {\scriptsize $\alpha \text{-} 4$};
			\fill [black] (7/4,-0.72) circle (0.06);
			\draw (7/4,-0.72-1/16)  node [below,red] {\scriptsize $\alpha \text{-} 4$};
			\fill [black] (1.025,-1.025) circle (0.06);
			\draw(1.025,-1.025)  node [below,red] {\scriptsize $1$};
		\end{scope}
		
		\begin{scope} [shift={(2+1/32,-2/3 )}, xscale=-1, yscale=1]
			\draw [thick,black] plot  coordinates {(0,-8/8) (0,1/4) };
			\fill [black] (0,0) circle (0.06);
			\draw(0-3/16,0)  node [above,red] {\scriptsize $1$};
			\begin{scope}[shift={(0,-3/9)}]
				\draw [thick,black] plot  coordinates {(3.6/32,0) (-3.6/32,0) };
			\end{scope}
			\begin{scope}[shift={(0,-6/9)}]
				\draw [thick,black] plot  coordinates {(3.6/32,0) (-3.6/32,0) };
			\end{scope}
		\end{scope}
		
		\begin{scope}[xscale=-1]
			\draw [thick,black] plot  coordinates {(3/3,-3/3) (-1/4,1/4) };
			\begin{scope}[shift={(-1/3,+1/3)}]
				\draw [thick,cborange] plot  coordinates {(2/3,-1) (8/3,-1) };
				\fill [black] (+3/3,-1) circle (0.06);
				\draw (+3/3,-1)  node [below,red] { $\frac{\beta}{4}$};
				\begin{scope}[shift={(-1/3,0)}]
					\draw [thick,cbblue] plot  coordinates {(2,-4/3) (2,1/8) };
					\fill [black] (2,-1) circle (0.06);
					\draw (2+1/16,-1+6/16)  node [right,red] { $\frac{4\text{-}\beta}{12}$};
					\draw (2,0) node [left] {\hbox{\footnotesize$\geqone$}};
				\end{scope}
			\end{scope}
			
			\begin{scope} [shift={(2,-2/3 )}, xscale=-1, yscale=1]
				\draw [thick,black] plot  coordinates {(0,-8/8) (0,1/4) };
				\fill [black] (0,0) circle (0.06);
				\draw(0-3/16,0)  node [above,red] { $\frac{\beta}{4}$};
				\begin{scope}[shift={(0,-3/9)}]
					\draw [thick,black] plot  coordinates {(3.6/32,0) (-3.6/32,0) };
				\end{scope}
				\begin{scope}[shift={(0,-6/9)}]
					\draw [thick,black] plot  coordinates {(3.6/32,0) (-3.6/32,0) };
				\end{scope}
			\end{scope}
		\end{scope}
		
		\begin{scope}[shift={(0,0)}]
			\fill [white] (0,0) circle (0.07);
			\draw [thick] (0,0) circle (0.07);
					\draw (0+0/16,0+1/16)  node [above,red] { $\frac{\epsilon}{2}$};
		\end{scope}
				\begin{scope}[shift={(1/3,-1/3)}]
			\draw [thick,black] plot  coordinates {(-1/12,-1/12) (1/12,1/12) };
		\end{scope}
		\begin{scope}[xscale=-1]
			\begin{scope}[shift={(1/3,-1/3)}]
				\draw [thick,black] plot  coordinates {(-1/12,-1/12) (1/12,1/12) };
			\end{scope}
		\end{scope}
\end{tikzpicture}}

\def\FigTwoCaseEFGix{\begin{tikzpicture}[baseline=12pt]
		\draw (7/4,0/4) node {};
		\draw [thick,black] plot  coordinates {(3/3,-3/3) (-1/4,1/4) };
		\begin{scope}[shift={(-1/3,+1/3)}]
			\draw [thick,cborange] plot [smooth, tension=1] coordinates {(2/3,-9/8) (2,-4/8) (4/3,1/8) };
			\draw [thick,cborange] plot [smooth, tension=1] coordinates {(7/2-2/3,-9/8) (7/2-2,-4/8) (7/2-4/3,1/8) };
			\fill [black] (7/4,-0.075) circle (0.06);
			\draw (7/4,-0.075+1/16)  node [above,red] {\scriptsize $\alpha \text{-} 4$};
			\fill [black] (7/4,-0.72) circle (0.06);
			\draw (7/4,-0.72-1/16)  node [below,red] {\scriptsize $\alpha \text{-} 4$};
			\fill [black] (1.025,-1.025) circle (0.06);
			\draw(1.025,-1.025)  node [below,red] {\scriptsize $1$};
		\end{scope}
		
		\begin{scope} [shift={(2+1/32,-2/3 )}, xscale=-1, yscale=1]
			\draw [thick,black] plot  coordinates {(0,-8/8) (0,1/4) };
			\fill [black] (0,0) circle (0.06);
			\draw(0-3/16,0)  node [above,red] {\scriptsize $1$};
			\begin{scope}[shift={(0,-3/9)}]
				\draw [thick,black] plot  coordinates {(3.6/32,0) (-3.6/32,0) };
			\end{scope}
			\begin{scope}[shift={(0,-6/9)}]
				\draw [thick,black] plot  coordinates {(3.6/32,0) (-3.6/32,0) };
			\end{scope}
		\end{scope}
		
	\begin{scope}[xscale=-1]
		\draw [thick,black] plot  coordinates {(3/3,-3/3) (-1/4,1/4) };
		\begin{scope}[shift={(-1/3,+1/3)}]
			\draw [thick,cborange] plot  coordinates {(2/3,-1) (8/3,-1) };
			\fill [black] (+3/3,-1) circle (0.06);
			\draw (+3/3,-1)   node [below,red] {\scriptsize $1$};
			\begin{scope}[shift={(-1/3,0)}]
				\draw (2,-1) node [above] {\hbox{\footnotesize$\geqone$}};
			\end{scope}
		\end{scope}
		
		\begin{scope} [shift={(2,-2/3 )}, xscale=-1, yscale=1]
			\draw [thick,black] plot  coordinates {(0,-8/8) (0,1/4) };
			\fill [black] (0,0) circle (0.06);
			\draw (0-3/16,0)  node [above,red] {\scriptsize $1$};
			\begin{scope}[shift={(0,-3/9)}]
				\draw [thick,black] plot  coordinates {(3.6/32,0) (-3.6/32,0) };
			\end{scope}
			\begin{scope}[shift={(0,-6/9)}]
				\draw [thick,black] plot  coordinates {(3.6/32,0) (-3.6/32,0) };
			\end{scope}
		\end{scope}
	\end{scope}
		
		\begin{scope}[shift={(0,0)}]
			\fill [white] (0,0) circle (0.07);
			\draw [thick] (0,0) circle (0.07);
					\draw (0+0/16,0+1/16)  node [above,red] { $\frac{\epsilon}{2}$};
		\end{scope}
				\begin{scope}[shift={(1/3,-1/3)}]
			\draw [thick,black] plot  coordinates {(-1/12,-1/12) (1/12,1/12) };
		\end{scope}
		\begin{scope}[xscale=-1]
			\begin{scope}[shift={(1/3,-1/3)}]
				\draw [thick,black] plot  coordinates {(-1/12,-1/12) (1/12,1/12) };
			\end{scope}
		\end{scope}
\end{tikzpicture}}

\def\FigTwoCaseEFGx{\begin{tikzpicture}[baseline=12pt]
		\draw (7/4,0/4) node {};
		\draw [thick,black] plot  coordinates {(3/3,-3/3) (-1/4,1/4) };
		\begin{scope}[shift={(-1/3,+1/3)}]
			\draw [thick,cborange] plot [smooth, tension=1] coordinates {(2/3,-9/8) (2,-4/8) (4/3,1/8) };
			\draw [thick,cborange] plot [smooth, tension=1] coordinates {(7/2-2/3,-9/8) (7/2-2,-4/8) (7/2-4/3,1/8) };
			\fill [black] (7/4,-0.075) circle (0.06);
			\draw (7/4,-0.075+1/16)  node [above,red] {\scriptsize $\alpha \text{-} 4$};
			\fill [black] (7/4,-0.72) circle (0.06);
			\draw (7/4,-0.72-1/16)  node [below,red] {\scriptsize $\alpha \text{-} 4$};
			\fill [black] (1.025,-1.025) circle (0.06);
			\draw(1.025,-1.025)  node [below,red] {\scriptsize $1$};
		\end{scope}
		
		\begin{scope} [shift={(2+1/32,-2/3 )}, xscale=-1, yscale=1]
			\draw [thick,black] plot  coordinates {(0,-8/8) (0,1/4) };
			\fill [black] (0,0) circle (0.06);
			\draw(0-3/16,0)  node [above,red] {\scriptsize $1$};
			\begin{scope}[shift={(0,-3/9)}]
				\draw [thick,black] plot  coordinates {(3.6/32,0) (-3.6/32,0) };
			\end{scope}
			\begin{scope}[shift={(0,-6/9)}]
				\draw [thick,black] plot  coordinates {(3.6/32,0) (-3.6/32,0) };
			\end{scope}
		\end{scope}
		
		\begin{scope}[xscale=-1]
			\draw [thick,black] plot  coordinates {(3/3,-3/3) (-1/4,1/4) };
			\begin{scope}[shift={(-1/3,+1/3)}]
				\draw [thick,cborange] plot [smooth, tension=1] coordinates {(2/3,-9/8) (2,-4/8) (4/3,1/8) };
				\draw [thick,cborange] plot [smooth, tension=1] coordinates {(7/2-2/3,-9/8) (7/2-2,-4/8) (7/2-4/3,1/8) };
				\fill [black] (7/4,-0.075) circle (0.06);
				\draw (7/4,-0.075+1/16)  node [above,red] {\scriptsize $\beta \text{-} 4$};
				\fill [black] (7/4,-0.72) circle (0.06);
				\draw (7/4,-0.72-1/16)  node [below,red] {\scriptsize $\beta \text{-} 4$};
				\fill [black] (1.025,-1.025) circle (0.06);
				\draw(1.025,-1.025)  node [below,red] {\scriptsize $1$};
			\end{scope}
			
			\begin{scope} [shift={(2+1/32,-2/3 )}, xscale=-1, yscale=1]
				\draw [thick,black] plot  coordinates {(0,-8/8) (0,1/4) };
				\fill [black] (0,0) circle (0.06);
				\draw(0-3/16,0)  node [above,red] {\scriptsize $1$};
				\begin{scope}[shift={(0,-3/9)}]
					\draw [thick,black] plot  coordinates {(3.6/32,0) (-3.6/32,0) };
				\end{scope}
				\begin{scope}[shift={(0,-6/9)}]
					\draw [thick,black] plot  coordinates {(3.6/32,0) (-3.6/32,0) };
				\end{scope}
			\end{scope}
		\end{scope}
		
		\begin{scope}[shift={(0,0)}]
			\fill [white] (0,0) circle (0.07);
			\draw [thick] (0,0) circle (0.07);
					\draw (0+0/16,0+1/16)  node [above,red] { $\frac{\epsilon}{2}$};
		\end{scope}
				\begin{scope}[shift={(1/3,-1/3)}]
			\draw [thick,black] plot  coordinates {(-1/12,-1/12) (1/12,1/12) };
		\end{scope}
		\begin{scope}[xscale=-1]
			\begin{scope}[shift={(1/3,-1/3)}]
				\draw [thick,black] plot  coordinates {(-1/12,-1/12) (1/12,1/12) };
			\end{scope}
		\end{scope}
\end{tikzpicture}}

\def\includefigEFG#1{
\makebox(140pt,70pt){#1}}

\def\FigTwoCaseE{
\begin{tabular}{|c|}
\hline{\large\strut} (E1) \\
\hline
\includefigEFG{\FigTwoCaseEFGi} \\
\hline\end{tabular}}

\def\FigTwoCaseF{
\begin{tabular}{|c|c|c|}
\hline{\large\strut} (F1) & (F2) & (F3) \\
\hline{\Large\strut}
$\alpha<4$ & $\alpha=4$ & $\alpha>4$\\
\hline
\includefigEFG{\FigTwoCaseEFGii} &
\includefigEFG{\FigTwoCaseEFGiv} & 
\includefigEFG{\FigTwoCaseEFGvii}\\
\hline\end{tabular}}

\def\FigTwoCaseG{
\begin{tabular}{|c|c|c|}
\hline{\large\strut} (G1) & (G2) & (G3) \\
\hline{\Large\strut}
$4>\alpha\ge\beta$ & $4=\alpha>\beta$ & $\alpha=4=\beta$\\
\hline
\includefigEFG{\FigTwoCaseEFGiii} &
\includefigEFG{\FigTwoCaseEFGv} & 
\includefigEFG{\FigTwoCaseEFGvi}\\
\hline\hline
{\large\strut} (G4) & (G5) & (G6) \\ 
\hline
{\Large\strut} $\alpha>4>\beta$ & $\alpha>4=\beta$ & $\alpha\ge\beta>4$ \\
\hline
\includefigEFG{\FigTwoCaseEFGviii} & 
\includefigEFG{\FigTwoCaseEFGix} &
\includefigEFG{\FigTwoCaseEFGx} \\
\hline 
\end{tabular}}

\def\FigTwoUnmarkedi{\begin{tikzpicture}[baseline=13pt]
		\draw [thick,black] plot  coordinates {(-1/2,0) (4/2,0) };
		\draw (-2/4,0) node [left] {\hbox{\footnotesize$\geqtwo$}};
\end{tikzpicture}}

\def\FigTwoUnmarkedii{\begin{tikzpicture}[baseline=13pt]
				\draw [thick, black] (0,0) arc (-45:225:0.4);
				\draw [thick, black] plot coordinates {(0,0) (-1,-1)};
				\begin{scope}[shift={(-0.707106781*0.8,0)},xscale=-1]
									\draw [thick, black] plot coordinates {(0,0) (-1,-1)};
				\end{scope}
		\fill [black] (-0.707106781*0.4,-0.4*0.707106781) circle (0.06);
		\draw (-0.707106781*0.4+2/16,-0.4*0.707106781-1/16) node [right,red] { $\frac{6\alpha \texttt{+}2\gamma \text{-}16 }{3}$};
		\draw (-3/4,-1/2) node [left] {\hbox{\footnotesize$\geqone$}};
\end{tikzpicture}}

\def\FigTwoUnmarkediii{\begin{tikzpicture}[baseline=13pt]
		\draw [thick, black] (0,0) arc (-45:225:0.4);
		\draw [thick, black] plot coordinates {(0,0) (-1,-1)};
		\begin{scope}[shift={(-0.707106781*0.4+1/2,-0.4*0.707106781-1/4)},xscale=-1]
					\draw[scale=0.5, domain=-1:1, smooth, variable=\x, thick, black] plot ({\x}, {\x*\x/2});
		\end{scope}
		
		\begin{scope}[shift={(-0.707106781*0.8,0)},xscale=-1]
			\draw [thick, black] plot coordinates {(0,0) (-0.707106781*0.4,-0.707106781*0.4)};
		\end{scope}
		\fill [black] (-0.707106781*0.4,-0.4*0.707106781) circle (0.06);
				\draw (-0.707106781*0.4-1/16,-0.4*0.707106781) node [left,red] {\scriptsize $2\beta \text{-}4$};
		\begin{scope}[shift={(1-0.8*0.707106781,0)},xscale=-1]
	\draw [thick, black] (0,0) arc (-45:225:0.4);
\draw [thick, black] plot coordinates {(0,0) (-1,-1)};
	\begin{scope}[shift={(-0.707106781*0.8,0)},xscale=-1]
	\draw [thick, black] plot coordinates {(0,0) (-0.707106781*0.4,-0.707106781*0.4)};
\end{scope}
		\fill [black] (-0.707106781*0.4,-0.4*0.707106781) circle (0.06);
			\draw (-0.707106781*0.4-1/16,-0.4*0.707106781) node [right,red] {\scriptsize $2\alpha \text{-}4$};

		\end{scope}

\end{tikzpicture}}

\def\FigTwoUnmarkediv{\begin{tikzpicture}[baseline=13pt]
		
							\draw[scale=0.4, domain=0:3*pi, smooth, variable=\x, thick, black] plot ({\x}, {cos(\x r)});
							\draw[scale=0.4, domain=0:3*pi, smooth, variable=\x, thick, black] plot ({\x}, {-cos(\x r)});
							\fill [black] (0.4*1/2*pi,0) circle (0.06);
										\draw (0.4*1/2*pi,4/16) node [above,red] {\scriptsize $2\gamma \text{-}4$};
							\fill [black] (0.4*3/2*pi,0) circle (0.06);
										\draw (0.4*3/2*pi,4/16) node [above,red] {\scriptsize $2\beta \text{-}4$};
							\fill [black] (0.4*5/2*pi,0) circle (0.06);
										\draw (0.4*5/2*pi,4/16) node [above,red] {\scriptsize $2\alpha \text{-}4$};
\end{tikzpicture}}

\def\FigTwoUnmarkedv{\begin{tikzpicture}[baseline=13pt]
				\draw [thick,black] plot  coordinates {(1/4,1/4) (-1,-1) };
				\draw [thick,black] plot  coordinates {(-1/4,1/4) (1,-1) };
						\fill [black] (0,0) circle (0.06);
							\draw (0,+2/16) node [above,red] {$\frac{4\text{-}5 \delta \texttt{+}\alpha \texttt{+}\beta \texttt{+}\gamma\texttt{+}3\epsilon}{6}$};
				\draw (-1,-1) node [left] {\hbox{\footnotesize$\geqone$}};
					\draw (1,-1) node [right] {\hbox{\footnotesize$\geqone$}};
\end{tikzpicture}}

\def\FigTwoUnmarkedvi{\begin{tikzpicture}[baseline=13pt]
		\draw [thick,black] plot  coordinates {(1/4,1/4) (-1,-1) };
			\draw [thick,black] plot  coordinates {(-1/4,1/4) (0,0) };
\begin{scope}[yscale=-1, shift={(1,1)}]
	\draw [thick, black] (0,0) arc (-45:225:0.4);
	\draw [thick, black] plot coordinates {(0,0) (-1,-1)};
	\begin{scope}[shift={(-0.707106781*0.8,0)},xscale=-1]
		\draw [thick, black] plot coordinates {(0,0) (-1/2,-1/2)};
	\end{scope}
	\fill [black] (-0.707106781*0.4,-0.4*0.707106781) circle (0.06);
	\draw (-0.707106781*0.4+1/16,-0.4*0.707106781) node [right,red] {\scriptsize $2\alpha\texttt{+} 2 \delta \text{-}8$};
\end{scope}
			
		\fill [black] (0,0) circle (0.06);
									\draw (0,+2/16) node [above,red] {$\frac{8\text{-}6 \delta  \texttt{+}\beta \texttt{+}\gamma\texttt{+}3\epsilon}{6}$};
		
		\draw (-1,-1) node [left] {\hbox{\footnotesize$\geqone$}};
\end{tikzpicture}}

\def\FigTwoUnmarkedvii{\begin{tikzpicture}[baseline=13pt]
		\draw [thick,black] plot  coordinates {(-1/4,1/4) (0,0) };
\begin{scope}[yscale=-1, shift={(1,1)}]
	\draw [thick, black] (0,0) arc (-45:225:0.4);
	\draw [thick, black] plot coordinates {(0,0) (-1,-1)};
	\begin{scope}[shift={(-0.707106781*0.8,0)},xscale=-1]
		\draw [thick, black] plot coordinates {(0,0) (-1/2,-1/2)};
	\end{scope}
	\fill [black] (-0.707106781*0.4,-0.4*0.707106781) circle (0.06);
	\draw (-0.707106781*0.4+1/16,-0.4*0.707106781) node [right,red] {\scriptsize $2\alpha\texttt{+} 2 \gamma \text{-}8$};
\end{scope}
		\begin{scope}[xscale=-1]
				\draw [thick,black] plot  coordinates {(-1/4,1/4) (0,0) };
\begin{scope}[yscale=-1, shift={(1,1)}]
	\draw [thick, black] (0,0) arc (-45:225:0.4);
	\draw [thick, black] plot coordinates {(0,0) (-1,-1)};
	\begin{scope}[shift={(-0.707106781*0.8,0)},xscale=-1]
		\draw [thick, black] plot coordinates {(0,0) (-1/2,-1/2)};
	\end{scope}
	\fill [black] (-0.707106781*0.4,-0.4*0.707106781) circle (0.06);
	\draw (-0.707106781*0.4+1/16,-0.4*0.707106781) node [left,red] {\scriptsize $2\beta\texttt{+} 2 \gamma \text{-}8$};
\end{scope} 
		\end{scope}
		
		\fill [black] (0,0) circle (0.06);
			\draw (0,0+1/16) node [above,red] {$\frac{4\text{-}2\gamma\texttt{+}\epsilon}{4}$};
\end{tikzpicture}}

\def\includefigUnmarked#1{
	\makebox(130pt,85pt){#1}}
\def\FigTwoUnmarked{
\begin{tabular}{|c|c|c|}
	\hline
	{\large\strut} (I) & (II) & (III) \\
	\hline
	{\Large\strut}
$3\alpha +\gamma\leq 8 \ \land  \  \delta \geq \frac{4+\alpha+2\gamma}{5}$ 
	&  \vphantom{\huge\strut} \makecell{ $3 \alpha+\gamma >8 \ \ \land  \ $ \\ $3\delta-\gamma>4\geq \beta+\gamma$ } 
	& $\gamma=2<\beta$ \\
	\hline
	\includefigUnmarked{\FigTwoUnmarkedi} & 
	\includefigUnmarked{\FigTwoUnmarkedii} & 
	\includefigUnmarked{\FigTwoUnmarkediii} \\
	\hline\hline
	{\large\strut} (IV) & (V) & (VI) \\
	\hline
	{\Large\strut}
	 $\gamma>2$
	& $\alpha+\delta\leq 4 \ \land \  \delta<\frac{4+\alpha+2\gamma}{5}$ 
	& \vphantom{\huge\strut} \makecell{$3\delta-\gamma<4\geq \beta+\gamma$ \\   $ \land\ \ \delta+\alpha>4$} \\
	\hline
	\includefigUnmarked{\FigTwoUnmarkediv} &
	\includefigUnmarked{\FigTwoUnmarkedv} & 
	\includefigUnmarked{\FigTwoUnmarkedvi}   \\
	\hline 
	\noalign{\global\arrayrulewidth=2.0pt}
	\arrayrulecolor{white}\hline
	\noalign{\global\arrayrulewidth=0.4pt}
	\arrayrulecolor{black} \cline{1-1}
	{\large\strut} (VII) & \multicolumn{1}{c}{} & \multicolumn{1}{c}{}  \\
	\cline{1-1}
	{\Large\strut}

	$\gamma<2 \ \land \  \beta+\gamma>4$
	& \multicolumn{1}{c}{}
	& \multicolumn{1}{c}{}  \\
	\cline{1-1}
	\includefigUnmarked{\FigTwoUnmarkedvii}   &
	\multicolumn{1}{c}{} & 
	\multicolumn{1}{c}{} \\
	\cline{1-1}

	\end{tabular}}

\def\bigtimes{\mathop{\raise-2pt\hbox{\huge$\times$}}}

\newbox\circbulletbox
\setbox\circbulletbox=\hbox{\,{\large$\circledcirc$}\kern-8.7pt\raise .6pt\hbox{$\bullet$}\,}

\let\ge\geqslant
\let\leq\leqslant
\let\geq\geqslant

\def\circVbig{\hbox{\text{\it\r{V}}}}
\def\circVscript{\hbox{\scriptsize\text{\it\r{V}}}}
\def\circVscriptscript{\mbox{\tiny\text{\it\r{V}}}}

\def\circVlimits_#1^#2{{\mathchoice%
{\circVbig{}^{\kern2pt #2}_{\kern-2pt #1}}%
{\circVbig{}^{\kern2pt #2}_{\kern-2pt #1}}%
{\scriptstyle\circVscript{}^{\kern1.7pt #2}_{\kern-1pt #1}}%
{\scriptscriptstyle\circVscriptscript{}^{\kern1.5pt #2}_{\kern-1pt #1}}%
}}
\def\circVr_#1{\circVlimits_#1^r}
\def\circVs_#1{\circVlimits_#1^s}

\def\circWbig{\hbox{\text{\it\r{W}}}}
\def\circWscript{\hbox{\scriptsize\text{\it\r{W}}}}
\def\circWscriptscript{\mbox{\tiny\text{\it\r{W}}}}

\def\circWlimits_#1^#2{{\mathchoice%
{\circWbig{}^{\kern2pt #2}_{\kern-2pt #1}}%
{\circWbig{}^{\kern2pt #2}_{\kern-2pt #1}}%
{\scriptstyle\circWscript{}^{\kern1.7pt #2}_{\kern-1pt #1}}%
{\scriptscriptstyle\circWscriptscript{}^{\kern1.5pt #2}_{\kern-1pt #1}}%
}}

\def\OM{\mathchoice
{\rlap{\kern3.2pt$\overline{\phantom{L}}$}M}
{\rlap{\kern3.2pt$\overline{\phantom{L}}$}M}
{\rlap{\kern2.4pt$\scriptstyle\overline{\phantom{L}}$}M}
{\rlap{\kern1.8pt$\scriptscriptstyle\overline{\phantom{L}}$}M}}

\def\mycirc{{\kern1pt\circ\kern2pt}}

\let\bbar\widehat

\def\Gal{\mathop{\rm Gal}\nolimits}

\def\Spec{\mathop{\rm Spec}\nolimits}

\let\phi\varphi
\let\theta\vartheta
\let\epsilon\varepsilon
\let\setminus\smallsetminus

\newcommand{\BP}{{\mathbb{P}}}
\newcommand{\BQ}{{\mathbb{Q}}}
\newcommand{\BR}{{\mathbb{R}}}

\newcommand{\Fm}{{\mathfrak{m}}}

\newcommand{\CC}{{\cal C}}

\newcommand{\CP}{{\cal P}}

\newcommand{\CY}{{\cal Y}}

\newcommand{\op}{\operatorname}

\newcommand{\Cst}{C^{\op{st}}}

\newcommand{\wbar}{{\overline{w}}}
\newbox\mybox
\def\arrover#1{\mathrel{
\setbox\mybox=\hbox spread 1.4em
{\hfil$\scriptstyle#1$\hfil}
\vbox{\offinterlineskip\copy\mybox
\hbox to\wd\mybox{\rightarrowfill}}}}

\def\larrover#1{\mathrel{
\setbox\mybox=\hbox spread 1.4em
{\hfil$\scriptstyle#1\vphantom{g}$\hfil}
\vbox{\offinterlineskip\copy\mybox
\hbox to\wd\mybox{\leftarrowfill}}}}

\def\ontoover#1{\mathrel{
\setbox\mybox=\hbox spread 1.4em
{\hfil$\scriptstyle#1\vphantom{g}$\hfil}
\vbox{\offinterlineskip\copy\mybox
\hbox to\wd\mybox{\rightarrowfill\hskip-2.8mm
$\rightarrow$}}}}
\def\leftontoover#1{\mathrel{
\setbox\mybox=\hbox spread 1.4em
{\hfil$\scriptstyle#1\vphantom{g}$\hfil}
\vbox{\offinterlineskip\copy\mybox
\hbox to\wd\mybox{$\leftarrow$\hskip-2.8mm
\leftarrowfill}}}}

\let\onto\twoheadrightarrow

\def\isoto{\mathrel{
\setbox\mybox=\hbox spread 0.9em
{\hfil$\scriptstyle\sim$\hfil}
\vbox{\offinterlineskip\copy\mybox
\hbox to\wd\mybox{\rightarrowfill}}}}

\def\Bigskip{\bigskip\bigskip}

\newtheorem{Thm}{Theorem}[section]
\newtheorem{Prop}[Thm]{Proposition}

\newtheorem{Conj}[Thm]{Conjecture}
\newtheorem{Def}[Thm]{Definition}

\newtheorem{Rem}[Thm]{Remark}
\newtheorem{Ex}[Thm]{Example}

\newtheorem{thmx}{Theorem}

\newtheorem{conjx}[thmx]{Conjecture}

\numberwithin{Thm}{subsection}
\def\UseTheoremCounterForNextEquation{\setcounter{equation}{\value{Thm}}\addtocounter{Thm}{1}}



\def\qed{{\hskip0pt\unskip\unskip\nobreak\hfil\penalty50
\hskip1em\hbox{}\nobreak\hfil
{$\square$}
\parfillskip=0pt\finalhyphendemerits=0
\par}\medskip}

\newenvironment{Proof}
{\noindent{\bf Proof.}\ }
{\qed}

{\noindent{\bf Proof of #1.}\ }
{\qed}

\begin{document}

\title{\strut
	\vskip-80pt
	The Classification of the Stable Marked Reduction \\
	 of Genus $2$ Curves in Residue Characteristic $2$
}
\author{
	\begin{minipage}{.3\hsize}
		Tim Gehrunger\\[12pt]
		\small Department of Mathematics \\
		ETH Z\"urich\\
		8092 Z\"urich\\
		Switzerland \\
		tim.gehrunger@math.ethz.ch\\[9pt]
	\end{minipage}
}
\date{\today}

\maketitle

\Bigskip

\begin{abstract}
		Consider a hyperelliptic curve of genus $2$ over a field $K$ of characteristic zero. After extending $K$ we can view it as a marked curve with its $6$ Weierstrass points. We classify the structure of the potentially stable reduction of such curves for a valuation of residue characteristic~$2$. We implement this classification into a computer algebra system and compute it for 
a list of curves defined over $\BQ$ with conductor at most $2^{20}$.

\end{abstract}

{\renewcommand{\thefootnote}{}
}

\newpage
\tableofcontents
\newpage

\section{Introduction}
\label{Intro}

{\bf 1.1 Motivation and strategy:}
Let $K$ be a valued field of characteristic $0$ and residue characteristic $2$. Moreover, let $C$ be a hyperelliptic curve over~$K$, that is, a curve defined by a Weierstrass equation $z^2=F(x)$ for some polynomial~$F$. After extending $K$ if necessary, the ramification points of the canonical double covering $\pi\colon C\onto\bar C \cong \mathbb{P}^1_K$ are defined over $K$ and there exists a stable marked model $\CC$ of $C$ with these points marked. From this model, one can read off the stable (unmarked) model of $C$, whose special fiber encodes many arithmetic properties of $C$ such as the local $L$-factor. 

In \cite{GP24}, Richard Pink and the author of this article provide an algorithm that describes $\CC$ explicitly. For genus $2$, we gave a classification of the possible combinatorial structure of the special fiber $C_0$ of $\mathcal{C}$.
In this article, we enrich the classification by the thicknesses of the double points of $C_0$. Moreover, we reprove part of the classification using elementary arguments, making it independent of the criterion of Liu in \cite{Liucriterionarticle}, which was previously used. Finally, we implement the classification into SageMath \cite{sagemath} and compute it for a list of curves defined over $\mathbb{Q}$ with conductor at most $2^{20}$.

In the case of residue characteristic $\neq 2$, the notion of the cluster picture as described in \cite{DokchitserDokchitser} is used to encode the combinatorics of the root configuration of $F(x)$. In that same article, several invariants of $C$ such as the conductor exponent, valuation of the minimal discriminant and the Tamagawa number are expressed in terms of the cluster picture and the action of $\Gal(L/K)$ on it, where $L/K$ is a finite Galois extension over which $C$ admits stable reduction. We expect that the classification presented in this article will eventually lead to similar results. However, the combinatorial complexity is much larger; as there are 54 possibilities for the dual graph of the special fiber $C_0$ of the stable marked model in residue characteristic $2$, whereas in residue characteristic $\neq 2$ there are only the $7$ possibilities listed in \cite[Fig.1]{GP21}.

\medskip
\textbf{1.2 Overview:} 
We now explain the content of this article in greater detail. First, we reduce the general case throughout to the case that $K$ is algebraically closed to avoid the recurring need for field extensions and the resulting changes in notation. Let $R$ denote its valuation ring with maximal ideal $\Fm$ and $k=R/\Fm$ its residue field. Let $v$ be the valuation on $K$ corresponding to $R$ normalized to $v(2)=1$. 


\textbf{Classification:}
Let $C$ be a genus $2$ curve and write $\bar \CC$ for the stable marked model of $\bar C$ with the $6$ branch points of $\pi$ marked. We denote the special fiber of $\bar \CC$ by $\bar C_0$. Then $\bar C_0$ has at most three even double points and at most one odd double point, and the combinatorial structure is invariant under symmetries interchanging the former transitively. Let $\alpha\ge\beta\ge\gamma\ge0$ denote the (normalized) thicknesses of the even double points and $\epsilon\ge0$ that of the odd double point, where we interpret $0$ as the thickness of a double point that does not exist. See Section \ref{ModelsofCurvesSubSection} for the definition of the normalized thickness. 
As $K$ is algebraically closed, there exists an equation of the form 
\UseTheoremCounterForNextEquation
\begin{equation}\label{IntroGenus2Equation} 
	z^2\ =\ F(x)\ =\ ax+bx^2+cx^3+dx^4+ex^5,
\end{equation}
where $v(F)=0$ and $v(a)=\alpha+2\epsilon$ and $v(e)=\beta$. We choose a square root $\sqrt{bd}$ of $bd$ and  set $\delta:=v\bigl(c-2\sqrt{bd}\kern2pt\bigr)$.
As established in \cite[Section 5.2]{G24}, the reduction type of $C$ depends only on $\alpha, \beta, \gamma, \epsilon$ and $\delta$. There are $54$ different reduction types, denoted (A1)-(A3), (B1)-(B11), (C1)-(C6), (D1)-(D24), (E1), (F1)-(F3) and (G1)-(G6).

Using the thickness bound from \cite[Prop. 5.1.1]{G24}, we compute the thicknesses of the double points of $C_0$ and express these in terms of $\alpha, \beta, \gamma, \epsilon$  and $\delta$. Moreover, we use this result to distinguish between the cases (B9) and (B10), (D19) and (D20), and (D23) and (D24), which in \cite{GP24} required invoking Liu's criterion from \cite{Liucriterionarticle}. 


Computing $\delta$ using its definition requires an equation of the form $\eqref{IntroGenus2Equation}$. When working with a concrete curve defined over a field that is not algebraically closed, producing such a normalized Weierstrass equation requires explicitly determining Weierstrass points over a (ramified) extension of the base field, which is inconvenient in practice. 
 Let $J_8(F)$ be the eighth Igusa invariant of $F$ as defined by Igusa in \cite{Igusa1960}. 
Define $\delta':=\frac{v(J_8(F))}{8}+2$. It turns out that for the purposes of the classification, one may replace $\delta$ by $\delta'$ (even though they are not always equal).

Why should it be more convenient to compute $\delta'$, when its definition also uses an equation of the form \eqref{IntroGenus2Equation}?
Because one can use that the Igusa invariants are weighted projective invariants. In particular, let $K_1\subset K$ be a discretely valued subfield and assume we are given a Weierstrass equation $y^2=\tilde F$ with $\tilde F \in  K_1[x]$. In our implementation in \cite{WorksheetsG25} we compute $\delta'$ by first computing $J_8(\tilde F)$ and then using that we know how $J_8$ transforms under the action of $\operatorname{PGL}_2(K)$, which in particular does not require working with an extension of $K_1$. 
In summary, we prove: 
\begin{thmx}[See Theorem \ref{Genus2SummaryTheorem} and Theorem \ref{DeltaPrimeTheorem}]\label{MainTheoremIntroduction}
%
%
There are 54 cases for the combinatorial structure of $C_0$. The space of parameters 
$(\alpha, \beta, \gamma, \delta', \epsilon)$ 
decomposes into half-open polyhedral regions\footnote{Here we mean regions in Euclidean space described by a finite set of linear equalities and strict linear inequalities.} $P_i\subset \BR^5_{\geq 0}$ associated to these 54 reduction types such that a curve $C$ is of reduction type $i$ if and only if the corresponding parameters $(\alpha, \beta, \gamma, \delta', \epsilon)$ are contained in $P_i$. The thicknesses of the double points of $C_0$ depend only on $(\alpha, \beta, \gamma, \delta', \epsilon)$. 
The same result holds true when replacing $\delta'$ with $\delta$. 
We present the 54 cases in  Figures \ref{FigTwoCaseA}-\ref{FigTwoCaseG}.
%
\end{thmx}
%
%

\textbf{Implementation and curves with small conductor:} Using this reformulation  of the genus $2$ classification in terms of $\delta'$, we implement the classification into SageMath. Given a Weierstrass equation of a curve defined over $\BQ_2$, it returns $(\alpha, \beta, \gamma, \delta', \epsilon)$ from Theorem \ref{MainTheoremIntroduction} and the type of reduction of $C_0$. For $\alpha, \beta, \gamma$ and $\epsilon$, we use the existing implementation of the cluster picture package \cite{best_vanbommel_cluster_pictures_github}, accompanying \cite{userguide}.

We then use our implementation to compute the classification for a list of genus $2$ curves with conductor at most $2^{20}$ defined over $\BQ$. It turns out that 53 out of the 54 cases of the classification are realized for such curves, missing only (D21). The case (D19), whose arithmetic conditions on the thicknesses of $C_0$ and $\delta$ are similar to those of (D21) is the rarest of the other 53 cases and the curves realizing it have by far the largest mean conductor exponent out of all the cases. 
This motivates the following conjecture:
\begin{conjx}[= Conjecture \ref{RealisationConjecture}]
	All 54 reduction types of the genus $2$ classification are realized over $\mathbb{Q}$.	
	\end{conjx} 

Finally, we plot the behavior of the conductor exponent and the valuation of the minimal discriminant for some one-dimensional families. 

\medskip
{\bf 1.3 Structure of the paper:}
Section \ref{SectionSetUp} contains preparatory material: In Subsections \ref{ModelsofCurvesSubSection} and \ref{SubSectionHyperellipticCurves} we review basic facts about semistable and stable marked models of hyperelliptic curves over~$R$. 
In the final Subsection \ref{ThicknessBoundSubSection} we state the thickness bound that is later used to refine the classification.

In Section \ref{ClassificationSection} we turn to the genus $2$ classification. In Subsection \ref{Genus2SubSection}, we present the refined classification including the thicknesses of the double points of $C_0$. In Subsection \ref{SectionComputingDelta}, we define $\delta'$ and prove that replacing $\delta$ by $\delta'$ does not change the classification.

In Subsection \ref{SubSectionImplementation}, we explain and present the implementation of the classification into SageMath. Finally, in Subsection \ref{SubSectionCurvesWithSmallConductor}, we summarize the results for the curves with conductor at most $2^{20}$.

\medskip
{\bf 1.4 Relation with other work:}
This article builds on \cite{GP24}, which was written by Richard Pink and the author of this article and where the classification that is extended in this article was established.

Historically, semistable reductions of hyperelliptic curves have mainly been studied and constructed in residue characteristic $\not =$ 2, see for example the construction of Bosch in \cite{Bosch1980}. A more recent approach including a classification for semistable genus $2$ curves is the article \cite{DokchitserDokchitser} by Dokchitser, Dokchitser, Maistret, and Morgan, which explicitly constructs the minimal regular model of $C$ and deduces the stable model by iteratively contracting unstable components of the special fiber. The authors then use the classification to compute many arithmetic invariants. 
%
%

Moreover, in  \cite{Liucriterionarticle}, Liu gives criteria for the type of the stable reduction of the unmarked curve $C$ in terms of Igusa invariants. This result is first proved in the setting of $\op{char}(k)\neq 2$ and carries over to the wild case by a moduli argument.  

In \cite{NoteOnOrdinary}, Dokchitser and Morgan show that the cluster picture controls the local arithmetic for hyperelliptic curves  with good ordinary reduction in the case of residue characteristic $2$. For genus $2$, this corresponds to case $(D1)$, which is the only case in our classification where the cluster picture carries the same information as the reduction type of $C$.

{\bf 1.5 Acknowledgments:}
First and foremost, I want to thank Richard Pink and Robert Nowak for many helpful and interesting conversations regarding the content of this paper. I am very grateful to Raymond van Bommel for his help with the cluster picture package and for Andrew Booker and Andrew Sutherland for granting me access to their data set of genus $2$ curves before the release of their own article. 
I also thank Johannes Schmitt and Jeremy Feusi for many valuable comments on earlier versions of this paper. I am very grateful to the anonymous reviewer for providing many valuable comments, which greatly improved the quality of this work.

\section{Set-up}\label{SectionSetUp}

This article uses the notation and conventions of \cite{GP24} and \cite{G24}, which we will briefly review in this section.  
Throughout this article let $R_1$ be a complete discrete valuation ring with quotient
field $K_1$. At several places we will need to replace $R_1$ by its integral closure in a finite
extension of $K_1$. As in \cite{GP24}, we find it more convenient to
work over an algebraic closure instead. So we fix an algebraic closure $K$ of $K_1$ and let $R$
denote the integral closure of $R_1$ in $K$. Since $R_1$ is complete, the valuation on $K_1$ extends
to a unique valuation with values in $\BQ$ on $K$, whose associated valuation ring is $R.$
Let $v$ denote a corresponding valuation on $K$, and let $\Fm$ be the maximal ideal and $k := R/\Fm$ the
residue field of $R$.
By \cite[Thm 2.1.5]{GP24} we can work over
the single field $K$ and avoid cumbersome changes of notation. 

\medskip

From now on we assume that $K$ has characteristic $0$ and $k$ has characteristic $2$. We normalize the valuation $v$ on $K$ in such a way that $v(2)=1$. For every integer $n\ge1$ we fix an $n$-th root $2^{1/n}\in K$ in a compatible way, such that for all $n,m\ge1$ we have $(2^{1/mn})^m=2^{1/n}$. For any rational number $\alpha = m/n$ we then set $2^\alpha := (2^{1/n})^m$. This defines a group homomorphism $\BQ \to K^\times$, which by the normalization of $v$ satisfies $v(2^\alpha)=\alpha$.

	\subsection{Models of curves}\label{ModelsofCurvesSubSection}
	
	Let $Y$ be a connected smooth proper algebraic curve over $K$. By a \emph{model of~$Y$} we mean a flat and finitely presented curve $\CY$ over $R$ with generic fiber~$Y$. We call such a model \emph{semistable} if the special fiber $Y_0$ is smooth except possibly for finitely many ordinary double points. Every double point $p\in Y_0$ then possesses an \'etale neighborhood in $\CY$ which is \'etale over $\Spec R[x,y]/(xy-a)$ for some nonzero $a\in \Fm$, such that $p$ corresponds to the point $x=y=0$. 	Here the valuation $v(a)$ depends only on the local ring of $\CY$ at~$p$, for instance by Liu \cite[\S10.3.2 Cor.\,3.22]{LiuAlgGeo2002}. We call $v(a)$ the \emph{(normalized) thickness} of $p$. 
	
	\begin{Rem}\rm 
		In the literature, the thickness (henceforth called \emph{classical thickness}) of a double point is usually defined only for models over a discrete valuation ring $R_1$ and uses a valuation  $v_1$ normalized to $v_1(\pi)=1$, where $\pi$ is a uniformizer of $R_1$, see for example 	\cite[\S10.3.2 Def.\,3.23]{LiuAlgGeo2002}.  Hence the classical thickness of a double point is always an integer, but is not stable under replacing $R_1$ by a finite extension of itself. 
		 If $Y$ and $\CY$ are defined over $\op{Frac}{R_1}\subset K$ and $R_1\subset R$ respectively and the restriction of $v$ to $\op{Frac}R_1$ is equivalent to $v_1$, the classical thickness of a double point is  $v_1(2)$ times the normalized thickness. 
		 
		 The special fiber of the stable model together with the classical thicknesses of the double points determines the special fiber of the minimal regular model (see \cite[\S10.3.2 Cor.\,3.25]{LiuAlgGeo2002}). 
In particular, for a genus $2$ curve over a discrete valuation ring which admits semistable reduction, our results can be used to determine the special fiber of the minimal regular model. Finding an extension over which a curve admits semistable reduction can be done using the methods established in \cite{GP24}. Note that in many cases, such an extension can also be found using the mclf-package for SageMath \cite{wewers_ruth_mclf_github}. We plan to address this in detail in a future article.
%
%
	\end{Rem}
	
	\medskip
	Any model is an integral separated scheme. Thus for any two models $\CY$ and $\CY'$ over~$R$, the identity morphism on $Y$ extends to at most one morphism $\CY\to\CY'$. If this morphism exists, we say that $\CY$ \emph{dominates}~$\CY'$. This defines a partial order on the collection of all models of $Y$ up to isomorphism. By blowing up one model one can construct many other models that dominate it. 
	
	\medskip
	Now consider an integer $n\ge0$ and distinct $K$-rational points $P_1,\ldots,P_n \in Y(K)$. This turns $Y$ into a \emph{smooth semistable marked curve} $(Y,P_1,\ldots,P_n)$ over~$K$. If $\CY$ is a semistable model of $Y$ such that these points extend to pairwise disjoint sections $\CP_1,\ldots,\CP_n \in \CY(R)$ which avoid all double points of the special fiber, we call $(\CY,\CP_1,\ldots,\CP_n)$ a \emph{semistable model} of $(Y,P_1,\ldots,P_n)$ over~$R$. 
	
	\medskip
	From now on we assume that $Y$ is of genus $g$ and that $2g+n\ge3$. Then the group of automorphisms of $Y$ which preserve the given points is finite, and $(Y,P_1,\ldots,P_n)$ is a smooth \emph{stable marked curve}. A \emph{stable model} of $(Y,P_1,\ldots,P_n)$ over $R$ is a semistable model such that the group of automorphisms of the closed fiber which preserve the given sections is finite as well. A semistable model is stable if and only if its closed fiber possesses no irreducible component that is isomorphic to $\BP^1_k$ and contains at most two double or marked points. 
	
	\begin{Prop}\label{AbsStabMod}
		\begin{enumerate}
			\item[(a)] A stable model $(\CY,\CP_1,\ldots,\CP_n)$ of $(Y,P_1,\ldots,P_n)$ exists.
			\item[(b)] This model is unique up to unique isomorphism.
			\item[(c)] For every semistable model $(\CY',\CP_1',\ldots,\CP_n')$ the model $\CY'$ dominates $\CY$.
		\end{enumerate}
	\end{Prop}
	\begin{Proof}
		See Liu \cite[2.19-21]{Liu2006} or Temkin \cite[1.2-5]{Temkin2010} or Cuzub \cite[Th.\,3.4]{Cuzub2018}.
	\end{Proof}
	
	\subsection{Hyperelliptic curves}\label{SubSectionHyperellipticCurves}

	Now let $C$ be a \emph{hyperelliptic curve} of genus $2$ over~$K$. Thus $C$ is a connected smooth proper algebraic curve which comes with a double covering $\pi\colon C \onto \bar C$ of a rational curve $\bar C\cong\BP^1_K$. 
	
	%

	
%
	\medskip
	As $K$ has characteristic $0$, the covering $\pi$ is only tamely ramified, and by the Hurwitz formula it is ramified at precisely $6$ closed points, namely, at the Weierstrass points of~$C$. Let $P_1,\ldots,P_{6} \in C(K)$ denote these points and $\bar P_1,\ldots,\bar P_{6} \in \bar C(K)$ their images under~$\pi$. Here both $(C,P_1,\ldots,P_{6})$ and $(\bar C,\bar P_1,\ldots,\bar P_{6})$ are stable marked curves. 
	
	\medskip
	Let $(\CC, \CP_1, \dots, \CP_{6})$ and $(\bar \CC, \bar \CP_1, \dots, \bar \CP_{6})$
	be the stable  models of  $(C, P_1, \dots, P_{6})$   and $(\bar C, \bar P_1, \dots, \bar P_{6})$, respectively. 
	 By \cite[\S 10.3.3 Prop. 3.48]{LiuAlgGeo2002} and the uniqueness of the stable model, the quotient $\widehat{\CC}$ of $\CC$ by the extension of the hyperelliptic involution is a semistable marked model of $(\bar C, \bar P_1, \dots, \bar P_{6})$ which dominates $\bar \CC$.

	In \cite{GP24}, Richard Pink and the author of this article developed an algorithm to compute the model $\CC$ for arbitrary genus, which starts with $\bar \CC$ and also computes $\widehat \CC$.  
	
	When the residue characteristic of $K$ is not 2, the morphism $\widehat \CC \onto \bar \CC$ is an isomorphism and the combinatorial structure of $C_0$ only depends on $\bar C_0$, see for example \cite{GP21}. In the article \cite{DokchitserDokchitser} by Dokchitser, Dokchitser, Maistret, and Morgan,  the special fiber $\bar C_0$ is described in their notion of cluster pictures and several arithmetic properties of $C$ are computed in terms of it.

	
	\medskip
	For reference we collect the schemes that we have introduced in the following diagram \ref{AllCPDiagram}. Recall that we have natural morphisms $\CC\onto\bbar\CC\onto\bar\CC$ that are compatible with the given sections. We let $(C_0, p_1, \dots, p_{6})$ and $(\bbar C_0,\bbar p_1, \dots,\bbar p_{6})$ and $(\bar C_0,\bar p_1, \dots,\bar p_{6})$ denote the special fibers of $(\CC,\CP_1,\dots,\CP_{6})$ and $(\bbar\CC,\bbar\CP_1,\dots,\bbar\CP_{6})$ and $(\bar\CC,\bar\CP_1,\dots,\bar\CP_{6})$, respectively. Further recall that $\bbar \CC$ and $\bar \CC$ are models of $\bar C \cong \mathbb{P}^1_K$.
	\UseTheoremCounterForNextEquation

\begin{equation}\label{AllCPDiagram}
	\vcenter{\xymatrix{
		&	\ C\  \ar@{^{ (}->}[r] \ar@{->>}[d]_-\pi & 
			\ \CC\ \ar@{->>}[d]_-{} & 
			\ C_0\ \ar@{_{ (}->}[l] \ar@{->>}[d]_-{}\\
	\mathbb{P}^1_K	\ar@{=}[r]^-{\sim} &	\ \bbar C\ \ar@{^{ (}->}[r] \ar@{=}[d] & 
			\ \bbar\CC\ \ar@{->>}[d] & 
			\ \bbar C_0\ \ar@{->>}[d] \ar@{_{ (}->}[l]  \\
	\mathbb{P}^1_K	\ar@{=}[r]^-{\sim}&	\ \bar C\ \ar@{^{ (}->}[r] \ar@{->>}[d] &
			\ \bar\CC\ \ar@{->>}[d] & 
			\ \bar C_0\ \ar@{_{ (}->}[l] \ar@{->>}[d] \\
		&	\ \Spec K\ \ar@{^{ (}->}[r] & 
			\ \Spec R\  & 
			\ \Spec k\ \ar@{_{ (}->}[l] \\}}
\end{equation}

	\medskip
	To describe the relation between the special fibers $\bbar C_0$ and~$\bar C_0$, we use the terminology concerning the type 
	of an irreducible component of $\bbar C_0$ from \cite[Def 2.3.2]{GP24}. In particular, an irreducible component of $\widehat C_0$ is called 
	\begin{itemize}
		\item of type (a) if it maps isomorphically to an irreducible component of $\bar C_0$;
		\item of type (b) if it lies between irreducible components of type (a);
		\item of type (c) if it is not of type (a) or (b) and is not a leaf;
		\item of type (d) if it is not of type (a) or (b) and is a leaf.
	\end{itemize} 

	We also divide the double points of $\bar C_0$ into two classes. For this recall that $\bar C_0$ is marked with $6$ distinct points in the smooth locus. As the complement of a double point consists of two connected components, this divides the $6$ marked points into two groups.
	\begin{Def}\label{EvenOddDef}
		A double point $\bar p$ of $\bar C_0$ is called \emph{even} if each connected component of $\bar C_0 \setminus \{\bar p\}$ contains an even number of the points $\bar p_1,\ldots,\bar p_{6}$. Otherwise, it is called \emph{odd}. 
	\end{Def}
	
	\subsection{Thickness bound}\label{ThicknessBoundSubSection}
	
Let $X$ be an irreducible component of $\bar C_0$. Moreover, let $x$ be a global coordinate of $\bar \CC$ along $X$ and consider a hyperelliptic equation $z^2=F(x)$ of $C$ with  $F\in R[x^{\pm1}]$ and $v(F)=0$. In \cite[Thm.-Def. 2.4.1]{G24}, the square defect of $X$, denoted by $\wbar(X)$, was defined as  $$\wbar(X):=\min\{2, \, \sup\{v(F-H^2) \ | \ H\in R[x^{\pm1}] \}\}.$$ 
This depends only on $X$ and  is independent of the choice of $x$ and  $F$. 
	
Now consider a closed smooth unmarked point $\bar p\in X$.  Suppose that there exists an irreducible component of $\widehat C_0$ above $\bar p$. Fix such a component $\bbar T$ of type (d). Let $T$ be its inverse image in $C_0$ and write $g(T)$ for the genus of $T$. Furthermore, let $\CY$ be the semistable model of $\bar C$ obtained from blowing down all components of type (c) and (d) above $\bar p$ other than $\bbar T$. Then $\CY$ has a unique double point $\bar q$ above $\bar p$ of thickness $\epsilon$. By \cite[Prop. 5.1.1]{G24}, we have
\UseTheoremCounterForNextEquation
	\begin{equation}
\epsilon\leq \frac{2-\wbar(X)}{2g(T)+1} \label{ThicknessBoundInequality}
	\end{equation}
	 Moreover, equality holds if and only if $\bbar T$ is the only component of type (d) above $\bar p$.

\section{A classification for genus 2 curves}\label{ClassificationSection}
\subsection{Classification with thickness}
\label{Genus2SubSection}

In \cite{GP24}, a classification of the stable marked reduction for genus $g=2$ curves was given. Using  inequality \eqref{ThicknessBoundInequality} from \cite[Prop. 5.1.1]{G24}, we refine this classification by also providing the thicknesses of double points. In all cases but (A2), (B9), (D19) and (D23), this is a straightforward computation using the data already provided in the worksheets \cite{WorksheetsGP24} of \cite{GP24}. The detailed computations for these more difficult cases can be found in \cite{WorksheetsG25}.

Moreover, the verification of the classification in \cite{GP24} partly relied on the stable reduction criteria for the unmarked curve $C$ from Liu \cite[Thm.\,1]{Liucriterionarticle} in terms of Igusa invariants. This was precisely needed for distinguishing between (B9) and (B10), (D19) and (D20), and (D23) and (D24). In the worksheets, we provide an elementary way to derive the stable marked reduction in these cases, making the proof of the classification independent of Liu's work. 
\medskip 

 As a reminder, note that in the case of $g=2$ the morphism $C\onto\bar C$ has $6$ branch points and we start with the associated stable marked model $(\bar\CC,\bar\CP_1,\ldots,\bar\CP_6)$ of~$\bar C$. The seven possibilities for the combinatorial structure of its closed fiber $(\bar C_0,\bar p_1,\ldots,\bar p_6)$ are shown in Figure~\ref{FigTwoAll}. Filled circles signify even double points and empty circles odd double points. The arrows indicate the ways that one type can degenerate into another. 

\begin{figure}[h] \centering 
	\FigTwoAll
	\caption{The possibilities for $(\bar C_0,\bar p_1,\ldots,\bar p_6)$ in genus~$2$.}\label{FigTwoAll}
\end{figure}

We observe that $\bar C_0$ has at most three even double points and at most one odd double point, and that the combinatorial structure is invariant under symmetries interchanging the former transitively. Let $\alpha\ge\beta\ge\gamma\ge0$ denote the thicknesses of the even double points and $\epsilon\ge0$ that of the odd double point, where we interpret $0$ as the thickness of a double point that does not exist. 
Moreover, 
each even double point is connected to a unique leaf component with exactly two marked points. After possibly interchanging the marked points we can assume that $\bar\CP_1$ meets this component for the double point of thickness $\alpha$ if $\alpha>0$. We also assume that $\bar\CP_2$ has maximal distance from~$\bar\CP_1$, that is, that the sum of the thicknesses of the double points between them is maximal. Then $\bar\CP_2$ must meet the leaf component containing the double point of thickness $\beta$ if $\beta>0$. We identify $\bar C$  with $\BP^1_K$ in such a way that $\bar P_1$ is identified with~$0$ and $\bar P_2$ with~$\infty$. Then $C$ is defined by an equation of the form 
\UseTheoremCounterForNextEquation
\begin{equation}\label{Genus2Equation}
	z^2\ =\ F(x)\ =\ ax+bx^2+cx^3+dx^4+ex^5
\end{equation}
with $F\in K[x]$ separable of degree~$5$.
Rescaling $x$ and $z$ by factors in $K^\times$\!, we can now arrange to have $v(F)=0$ and 
$$\left\{\kern-4pt\begin{array}{rl}
	v(a) \kern-4pt&=\, \alpha+2\epsilon \\[3pt]
	v(e) \kern-4pt&=\, \beta \\[3pt]
	v(\operatorname{disc}(F)) \kern-4pt&=\, 2\alpha + 2\gamma + 6 \epsilon \\[3pt]
	v(b) \kern-4pt&\geq \, \epsilon \\[3pt]
	v(b) \kern-4pt&= \, \epsilon\ \ \hbox{if $\alpha>0$}
\end{array}\kern-4pt\right\}
\text{ or equivalently }
\left\{\kern-3pt\begin{array}{rl}
	\alpha \kern-4pt&=\, v(a)-2\epsilon \\[3pt]
	\beta \kern-4pt&=\, v(e) \\[3pt]
	\gamma \kern-4pt&=\, \frac{1}{2}v(\operatorname{disc}(F))-\alpha-3\epsilon \\[3pt]
	\epsilon \kern-4pt&= \, \min\{v(b),\frac{1}{2}v(\operatorname{disc}(F))-v(a)\}
\end{array}\kern-4pt\right\},
$$
where $\operatorname{disc}(F)$ denotes the discriminant of~$F$.
With the equation in this form we choose a square root $\sqrt{bd}$ of $bd$ and  set $\delta:=v\bigl(c-2\sqrt{bd}\kern2pt\bigr)$.
It turns out that the combinatorial structure of  $(C_0,p_1,\ldots,p_6)$  depends only on the values of $\alpha, \beta, \gamma, \delta$ and $\epsilon$, which are all $\geq 0$. 
The computations carried out in \cite{WorksheetsGP24} for each of the cases (A)-(G) show that we always have $\delta \geq \min \{2, \gamma\}$, with equality if $\gamma<\min \{\beta, 2\}$. 


\medskip

All in all, the seven cases 
from Figure \ref{FigTwoAll} divide into 54 subcases for the combinatorial structure of $(C_0,p_1,\ldots,p_6)$.
In each subcase we draw irreducible components of type
(a)  in black, those of type (b) in orange, those of type (c) in green,
and those of type (d) in blue. We also label any irreducible component of genus $g' > 0$ by $g=g'$, while all irreducible components of genus $0$ remain unlabeled. The thickness of a double point is written in red. In summary, we prove:

\begin{Thm}\label{Genus2SummaryTheorem}
		There are 54 cases for the combinatorial structure of $C_0$. The space of parameters 
		$(\alpha, \beta, \gamma, \delta, \epsilon)$ 
		decomposes into half-open polyhedral regions\footnote{Here we mean regions in Euclidean space described by a finite set of linear equalities and strict linear inequalities.} $P_i\subset \BR^5_{\geq 0}$ associated to these 54 reduction types such that a curve $C$ is of reduction type $i$ if and only if the corresponding parameters $(\alpha, \beta, \gamma, \delta, \epsilon)$ are contained in $P_i$. The thicknesses of the double points of $C_0$ depend only on $(\alpha, \beta, \gamma, \delta, \epsilon)$. We present the 54 cases in  Figures \ref{FigTwoCaseA}-\ref{FigTwoCaseG}.
\end{Thm}

\medskip


\begin{figure}[H] \centering 
	\FigTwoCaseA
	\caption{The possibilities for $(C_0,p_1,\ldots,p_6)$ in the case (A).}
	\label{FigTwoCaseA}
\end{figure}

\medskip
{\bf Case (A):} This is the case of ``equidistant geometry'' of Lehr-Matignon \cite{LehrMatignon06}, that is, where $\bar C_0$ is smooth. Hence we have $\alpha=\beta=\gamma=\epsilon=0$. In  \cite[\texttt{Case\_(A).mw}]{WorksheetsGP24}, it is shown that the combinatorial structure of $C_0$ depends only on~$\delta$ and is as written in Figure \ref{FigTwoCaseA}. In case (A1) and (A3),  inequality \eqref{ThicknessBoundInequality} is an equality and directly yields that the thicknesses of the double points of $\widehat C_0$ are $\tfrac{2}{3}$ in case (A1) and $\tfrac{2}{5}$ in case (A3). As there is a single double point of $C_0$ over each double point of $\widehat C_0$  this implies that the thicknesses of the double points of $C_0$ are $\tfrac{1}{3}$ in case (A1) and $\tfrac{1}{5}$ in case (A3).
The thicknesses for case (A2) are computed in \cite[\texttt{CaseDistinctions.mw}]{WorksheetsG25}.

\newpage
\medskip
{\bf Case (B):} Here $\bar C_0$ has exactly one even double point of thickness $\alpha>0$, and we have $\beta=\gamma=\epsilon=0$. The combinatorial structure of $C_0$ depends only on $\alpha$ and $\delta$ as established in  \cite[\texttt{Case\_(B).mw}]{WorksheetsGP24}, with the restriction that there we could not distinguish between (B9) and (B10). In \cite{GP24}, these were  distinguished using Liu's criterion \cite{Liucriterionarticle}. In \cite[\texttt{CaseDistinctions.mw}]{WorksheetsG25}, we present an independent way to distinguish the cases (B9) and (B10). As a corollary, we find the thicknesses of the double points in case (B9).  

	\begin{figure}[h!] \centering 
	\FigTwoCaseB
	\caption{The possibilities for $(C_0,p_1,\ldots,p_6)$ in the case (B).}\label{FigTwoCaseB}
\end{figure}

The 11 possible subcases for $C_0$ are sketched in Figure \ref{FigTwoCaseB}. 
In all of the cases but the already treated (B9), the thicknesses of the double points follow directly from the computation in the worksheets \cite[\texttt{Case\_(B).mw}] {WorksheetsGP24} together with  inequality \eqref{ThicknessBoundInequality} from \cite[Prop. 5.1.1]{G24}. We demonstrate this in the most complicated case (B6).

 Here, we obtain from \cite[\texttt{Case\_(B).mw}]{WorksheetsGP24} that $\widehat C_0$ has exactly six irreducible components $X_1, \dots, X_6$, such that $X_2$ intersects $X_1, X_3$ and $X_4$ and $X_5$ intersects $X_4$ and $X_6$. Furthermore, $X_1$ contains $4$ marked points and $X_6$ contains $2$ marked points. 
In the worksheet, it was shown that $X_2$ intersects $X_1$ in a double point of thickness $\frac{\delta}{2}$ and that $\wbar(X_2)=\tfrac{3\delta}{2}$. Furthermore, it was shown that the thickness of the double point in which $X_2$ intersects $X_4$ is $\frac{4-3\delta}{3}$ and the double point in which $X_4$ intersects $X_5$ is $\alpha+\delta-4$ and the double point in which $X_5$ intersects $X_6$ is $2$.                                                                                                                                                                                         
For the remaining node in which $X_2$ intersects  $X_3$, the bound \eqref{ThicknessBoundInequality} is an equality in this case and yields thickness $\frac{4-3\delta}{6}.$
Using that the preimage of a double point of $\widehat C_0$ consists either of two double points of the same thickness or one double point of half that thickness yields the thicknesses of the double points of $C_0$.

	\medskip
	{\bf Case (C):} Here $\bar C_0$ has two even double points of respective thicknesses $\alpha\ge\beta>0$, and we have $\gamma=\delta=\epsilon=0$. As established in  \cite[\texttt{Case\_(C).mw}]{WorksheetsGP24},  the combinatorial structure of $C_0$ depends only on $\alpha$ and $\beta$. As in case (B6) above, the thicknesses can be derived from the computations in the worksheets together with  inequality \eqref{ThicknessBoundInequality} from \cite[Prop. 5.1.1]{G24}. Alternatively, applying \cite[Proposition 5.2.4]{G24} to both double points directly yields the same result.
%
%
	The 6 possible subcases are sketched in Figure \ref{FigTwoCaseC}.
	\newpage
	\begin{figure}[H] \centering 
		\FigTwoCaseC
		\caption{The possibilities  for $(C_0,p_1,\ldots,p_6)$ in the case (C).}\label{FigTwoCaseC}
	\end{figure}

	
	{\bf Case (D):} 
	Here $\bar C_0$ has one irreducible component without marked points in the middle, which is connected by double points of thicknesses $\alpha\ge\beta\ge\gamma>0$ to leaf components with two marked points each. 
		The combinatorial structure of $C_0$ depends only on $\alpha, \beta, \gamma$ and $\delta$ as established in  \cite[\texttt{Case\_(D).mw}]{WorksheetsGP24}, with the restriction that we could not distinguish (D19) from (D20) and (D23) from (D24). In \cite{GP24}, these were  distinguished using Liu's criterion \cite{Liucriterionarticle}. In \cite[\texttt{CaseDistinctions.mw}]{WorksheetsG25}, we present an independent way to distinguish these cases. As a corollary, we find the thicknesses of the double points in case (D19) and (D23). 	The 24 possible subcases are sketched in Figures \ref{FigTwoCaseDone} and \ref{FigTwoCaseDtwo}.  As in case  (B6) above, the remaining thicknesses can be derived from the computations in the worksheets together with  inequality \eqref{ThicknessBoundInequality} from \cite[Prop. 5.1.1]{G24}.

			\begin{figure}[H] \centering 
		\FigTwoCaseDone
		\caption{The first 12 possibilities for $(C_0,p_1,\ldots,p_6)$ in the case (D).}\label{FigTwoCaseDone}
	\end{figure}


	\newpage
	
	\begin{figure}[H] \centering 
		\FigTwoCaseDtwo
		\caption{The remaining 12 possibilities for $(C_0,p_1,\ldots,p_6)$ in the case (D).}\label{FigTwoCaseDtwo}
	\end{figure}

	%

	\medskip
	{\bf Cases (E--G):} Here $\bar C_0$ has an odd double point~$\bar p$ of thickness $\epsilon>0$, and we always have $\gamma=\delta=0$.  As established in  \cite[\texttt{Cases\_(EFG).mw}]{WorksheetsGP24},  the combinatorial structure of $C_0$ depends only on $\alpha$ and $\beta$. As in case (B6) above, the thicknesses can be derived from the computations in the worksheets together with  inequality \eqref{ThicknessBoundInequality} from \cite[Prop. 5.1.1]{G24}.
	It turns out that the situation on each side of~$\bar p$ is the same as for the reduction of a curve of genus~$1$, and that the two sides are independent of each other.
	\medskip
	{\bf Case (E):} Here we have $\alpha=\beta=0$, and there is a single subcase only, which is sketched in Figure~\ref{FigTwoCaseE}.
	
	\begin{figure}[H] \centering 
		\FigTwoCaseE
		\caption{The single possibility for $(C_0,p_1,\ldots,p_6)$ in the case (E).}\label{FigTwoCaseE}
	\end{figure}
	
	\medskip
	{\bf Case (F):} Here we have $\alpha>\beta=0$. The 3 possible subcases are sketched in Figure~\ref{FigTwoCaseF}.
	
	\begin{figure}[H] \centering 
		\FigTwoCaseF
		\caption{The possibilities for $(C_0,p_1,\ldots,p_6)$ in the case (F).}\label{FigTwoCaseF}
	\end{figure}
	
	\medskip
	{\bf Case (G):} Here we have $\alpha\ge\beta>0$. The 6 possible subcases are sketched in Figure~\ref{FigTwoCaseG}.
	
	\begin{figure}[H] \centering 
		\FigTwoCaseG
		\caption{The possibilities for $(C_0,p_1,\ldots,p_6)$ in the case (G).}\label{FigTwoCaseG}
	\end{figure}
	
	From the above results, one can also determine the closed fiber $\Cst_0$ of the stable reduction of the unmarked curve~$C$ by consecutively contracting non-singular rational curves containing fewer than $3$ double points. The list of possible cases and their names are taken from Liu \cite[Th.\,1]{Liucriterionarticle}. There the case distinctions were given in terms of Igusa invariants, which are complicated polynomials in the coefficients of~$F$. Our results yield relatively simple conditions in terms of the numbers $\alpha,\beta,\gamma,\delta$ alone. The seven cases for $\Cst_0$ are sketched in Figure \ref{FigTwoCaseUnmarked}. 
	\begin{figure}[H] \centering 
		\FigTwoUnmarked
		\caption{The possibilities for the stable reduction of the unmarked curve.}\label{FigTwoCaseUnmarked}
	\end{figure}

	\subsection{Computing $\delta'$}\label{SectionComputingDelta}
	Let $F$ and $x$ be as in Section \ref{Genus2SubSection}. Let $J_8$ be the eighth Igusa invariant of $F$ as defined by Igusa in \cite{Igusa1960}. 
	Define $\delta':=\frac{v(J_8)}{8}+2$. Note that we earlier chose a square root of $bd$, which we denote by $\sqrt{bd}$. Now choose square roots of $b,d$ that we denote by $b^{\frac1 2},\ d^{\frac12}$ such that $\sqrt{bd}=b^{\frac12}d^{\frac12}$. 
	\begin{Thm}\label{DeltaPrimeTheorem}
%
		The \emph{same} half-open polyhedral regions from Theorem \ref{Genus2SummaryTheorem} describe the type of $C_0$ in terms of $(\alpha, \beta, \gamma, \delta', \epsilon)$ and the thicknesses of $C_0$ are given by the same formulas, where $\delta$ is replaced by $\delta'$.
	\end{Thm}
	\begin{Proof}
		We have  
\begin{align}
2^{16}	J_8 &= 3 (c-2b^{\frac12}d^{\frac12})^{8}+48 b^{\frac12}\, d^{\frac12}\, (c-2b^{\frac12}d^{\frac12})^{7} \nonumber\\
	&\quad +\left(-208 a e+304 d b\right) (c-2b^{\frac12}d^{\frac12})^{6} \nonumber\\
	&\quad +\left(960 b^{\frac{3}{2}} d^{\frac{3}{2}}-2496 b^{\frac12}\, d^{\frac12}\, a e+96 a \,d^{2}+96 b^{2} e\right) (c-2b^{\frac12}d^{\frac12})^{5} \nonumber\\
	&\quad +\left(1536 d^{2} b^{2}-10464 d b a e+960 d^{\frac{5}{2}} b^{\frac12}\, a+960 b^{\frac{5}{2}} d^{\frac12}\, e-1632 a^{2} e^{2}\right) (c-2b^{\frac12}d^{\frac12})^{4} \nonumber\\
	&\quad +\left(1024 b^{\frac{5}{2}} d^{\frac{5}{2}}-17152 b^{\frac{3}{2}} d^{\frac{3}{2}} a e+3072 d^{3} b a+3072 d \,b^{3} e\right.\nonumber\\
	&\qquad \left.-13056 b^{\frac12}d^{\frac12}\, a^{2} e^{2}-1152 d^{2} a^{2} e-1152 b^{2} a \,e^{2}\right) (c-2b^{\frac12}d^{\frac12})^{3} \nonumber\\
	&\quad +\left(-6144 d^{2} b^{2} a e+3072 d^{\frac{7}{2}} b^{\frac{3}{2}} a+3072 b^{\frac{7}{2}} d^{\frac{3}{2}} e-29184 d b \,a^{2} e^{2}\right.\nonumber\\
	&\qquad \left.-6912 d^{\frac{5}{2}} b^{\frac12}\, a^{2} e-6912 b^{\frac{5}{2}} d^{\frac12}\, a \,e^{2}+14080 a^{3} e^{3}+512 d^{4} a^{2}+512 b^{4} e^{2}\right) (c-2b^{\frac12}d^{\frac12})^{2} \nonumber\\
	&\quad +\left(6144 b^{\frac{5}{2}} d^{\frac{5}{2}} a e-12288 b^{\frac{3}{2}} d^{\frac{3}{2}} a^{2} e^{2}-8192 d^{3} b \,a^{2} e-8192 d \,b^{3} a \,e^{2}\right.\nonumber\\
	&\qquad \left.+56320 b^{\frac12}\, d^{\frac12}\, a^{3} e^{3}+2048 d^{\frac{9}{2}} b^{\frac12}\, a^{2}+2048 b^{\frac{9}{2}} d^{\frac12}\, e^{2}\right.\nonumber\\
	&\qquad -25600 d^{2} a^{3} e^{2}-25600 b^{2} a^{2} e^{3}\Big) (c-2b^{\frac12}d^{\frac12}) \nonumber\\
	&\quad -83200 a^{4} e^{4}+5120 b^{4} a \,e^{3}+5120 d^{4} a^{3} e+1792 d^{2} b^{2} a^{2} e^{2} \nonumber\\
	&\quad +125440 d b \,a^{3} e^{3}+2048 d^{\frac{7}{2}} b^{\frac{3}{2}} a^{2} e+2048 b^{\frac{7}{2}} d^{\frac{3}{2}} a \,e^{2} \nonumber\\ \nonumber
	&\quad -51200 d^{\frac{5}{2}} b^{\frac12}\, a^{3} e^{2}-51200 b^{\frac{5}{2}} d^{\frac12}\, a^{2} e^{3}.
\end{align}

We proceed by comparing the valuations of $2^{16}J_8$ and $(c-2b^{\frac12}d^{\frac12})=(c-2\sqrt{bd})$ in each of the 54 cases and show that replacing $\delta$ by $\delta'$ always yields the same case of the classification in Section \ref{Genus2SubSection} and that $\delta=\delta'$ whenever $\delta$ appears in a formula for the thickness of a double point. These are straightforward computations, which  can be found in \cite{WorksheetsG25}.
	\end{Proof}

%
\section{Examples}
\subsection{Implementation}\label{SubSectionImplementation}
In \cite[\texttt{genus2\_classification.sage}]{WorksheetsG25}, we provide an implementation of the genus 2 classification of Theorem \ref{MainTheoremIntroduction} into SageMath. 

 Here, we choose to compute $\delta'$ instead of $\delta$ as this is computationally cheaper. Moreover, we use the Cluster picture package \cite{best_vanbommel_cluster_pictures_github} for the computation of $\alpha, \beta, \gamma, \epsilon$. The main functions to call are:
\begin{itemize}
	\item \texttt{AlphaBetaGammaDeltaEpsilonComputation(F)} returns a list of the form $[\alpha, \beta, \gamma, \delta', \epsilon]$ for the curve defined by $y^2=F(x)$.
	\item \texttt{Genus2ClassificationLabel([alpha, beta, gamma, deltaprime, epsilon])} returns the case of the classification from   Section \ref{Genus2SubSection}.
	\item \texttt{ABGDEWithLabel(F)} returns a list of the form  $[\alpha, \beta, \gamma, \delta', \epsilon,\text{`label'}]$, where the label is the case of the classification from Section \ref{ClassificationSection} for the curve defined by $y^2=F(x)$.
	\item \texttt{StableMarkedModelLabel(F)} returns the case of the classification  from   Section \ref{Genus2SubSection} for the curve defined by $y^2=F(x)$.
	\item \texttt{StableModelLabel(F)} returns the case of the stable unmarked reduction from Figure \ref{FigTwoCaseUnmarked} for the curve defined by $y^2=F(x)$. 
	\item \texttt{StableModelLabelAndThicknesses(F)} returns a list of the form [label, thicknesses], where label is the case of the stable unmarked reduction from Figure \ref{FigTwoCaseUnmarked} and the thicknesses of the double points of the stable reduction are listed as they appear in the figure from left to right for the curve defined by $y^2=F(x)$.
\end{itemize}

\begin{Ex}{\rm
\phantom{..}
\begin{lstlisting}
	sage: AlphaBetaGammaDeltaEpsilonComputation(x^5+x^4-4*x^3-10*x^2+12*x)
	[0] : [1, 1/2, 1/2, 1, 0]
	sage: AlphaBetaGammaDeltaEpsilonComputation(x^5-1)
	[1] : [0, 0, 0, +Infinity, 0]
	sage: Genus2ClassificationLabel([1, 1/2, 1/2, 1, 0])
	[2] : 'D19'
	sage: ABGDEWithLabel(x^5+x^4-4*x^3-10*x^2+12*x)
	[3] : [1, 1/2, 1/2, 1, 0, 'D19']
	sage: StableModelLabelAndThicknesses(x^5+x^4-4*x^3-10*x^2+12*x)
	[4] : ['V', [1/6]]
\end{lstlisting}}
\end{Ex}

At present, the cluster picture package is often incapable of computing the cluster picture for polynomials defined over non-trivial extensions of $\mathbb{Q}_2$ due to SageMath not being able to factorize polynomials over these fields. Our code is structured in a way that any improvement in the coverage of the cluster picture package will automatically be reflected. 

\subsection{Curves with small conductor}\label{SubSectionCurvesWithSmallConductor}
In an upcoming article \cite{CurvesWithSmallConductor2025}, Booker and Sutherland list more than six million genus 2 curves over $\BQ$ with conductor at most $2^{20}$. Using our implementation, we compute $\alpha, \beta, \gamma, \delta', \epsilon$ and the label for each of these curves. Once the article \cite{CurvesWithSmallConductor2025} is published, we will add a file containing the curves together with the values we computed to \cite{WorksheetsG25}.

\begin{figure}[H]
	\centering
	\includegraphics[width=1\textwidth]{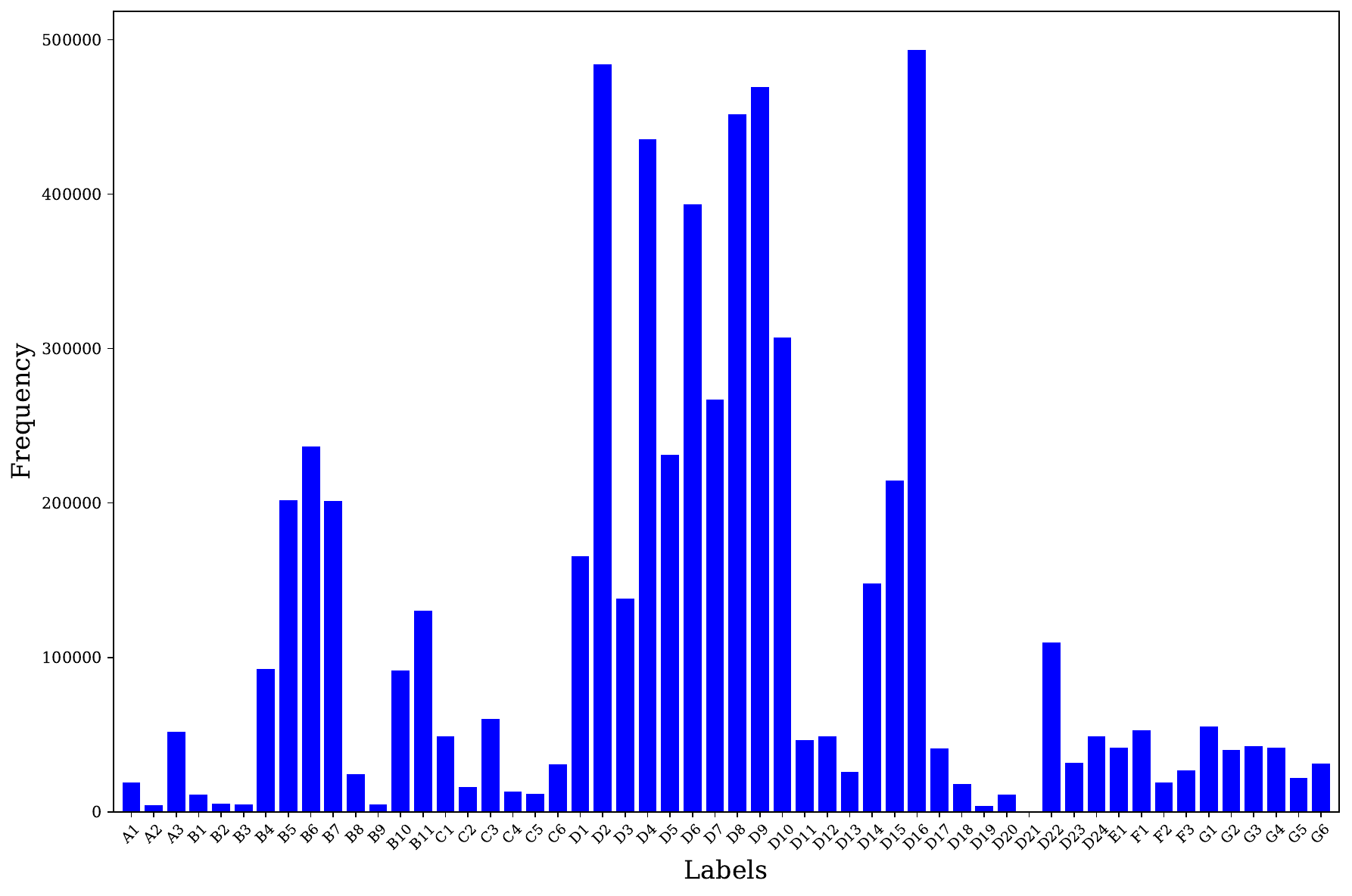}
	\caption{Label frequencies for the genus $2$ curves with conductor at most $2^{20}$}
	\label{FigureCaseDistribution}
\end{figure}

In Figure \ref{FigureCaseDistribution}, we plot the frequency of the labels from the genus $2$ classification. It turns out that 53  of the 54 cases of our classification are obtained by the curves in \cite{CurvesWithSmallConductor2025}, only missing case (D21). The case (D19), whose arithmetic conditions on $\alpha, \beta, \gamma$ and $\delta$ are very similar to those of (D21) is the rarest of the appearing cases, only appearing for 3961 curves and with large conductors. A possible explanation for this is that the conditions on $\alpha, \beta, \gamma, \delta$ for (D19) and (D21) force the conductor to be large for a curve defined over $\BQ$. We expect that there exist curves of case (D21) defined over $\mathbb{Q}$ but we have not yet found such a curve. In summary:

\begin{Conj}\label{RealisationConjecture}
	All 54 reduction types of the genus $2$ classification are realized over $\mathbb{Q}$.
\end{Conj}

\begin{figure}[h!]
	\centering
	\includegraphics[width=1\textwidth]{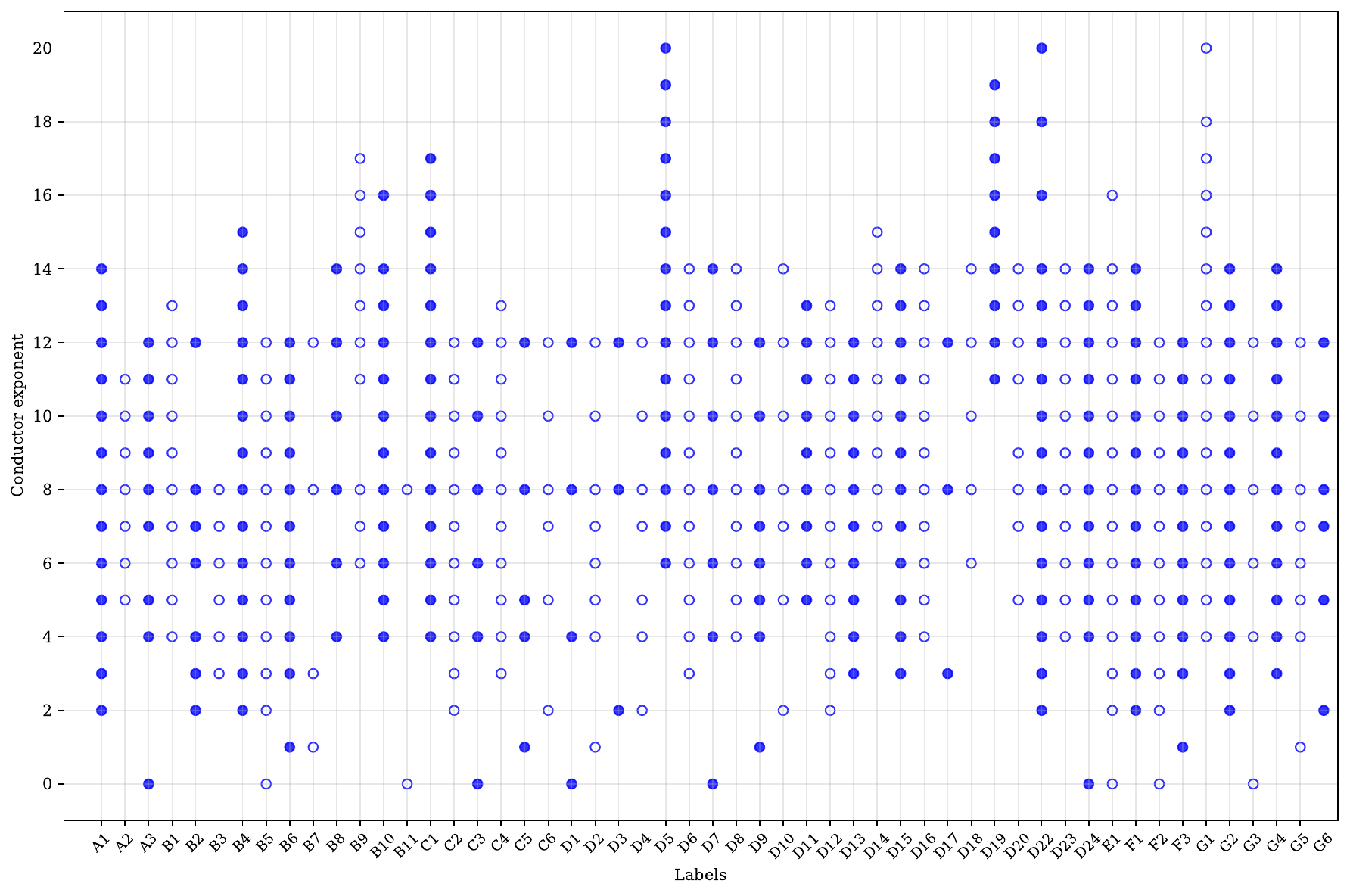}
	\caption{Conductor exponent values by cases for the genus $2$ curves with conductor at most $2^{20}$; the points are filled/unfilled solely for improved readability }
	\label{FigureConductorExponentsLabels}
\end{figure}

For a genus $2$ curve $C/\BQ$, we denote the $2$-adic valuation of its conductor $N$ as \emph{conductor exponent}.  Brumer and Kramer in \cite[Thm. 6.2]{ConductorBound} show that  the conductor exponent of a genus $2$ curve defined over $\mathbb{Q}$ is at most $20$. 
In Figure \ref{FigureConductorExponentsLabels}, we plot the conductor exponent values for the reduction types from the classification in Subsection \ref{Genus2SubSection}. It is likely that these are not all possible values for curves defined over $\BQ$. In particular, for curves in the data set, conductor exponent $20$ can only be obtained if the conductor is equal to $2^{20}$.
	
We proceed by plotting the conductor exponent  and the valuation of the minimal discriminant for some cases that only depend on one parameter. For a curve $C$, we denote its conductor by $N$ and its minimal discriminant by $D$. As explained in the introduction, in residue characteristic $p\neq 2$ the cluster picture as defined in \cite{DokchitserDokchitser} together with the action of $\Gal(L/\BQ_p)$ for the minimal extension $L/\mathbb{Q}_p$ over which the curve admits stable reduction determine the $p$-adic valuation of these. As our classification is finer than the cluster picture, we hope that the same might be true for $p=2$. 

\begin{figure}[h!]
	\centering
	\includegraphics[width=1\textwidth]{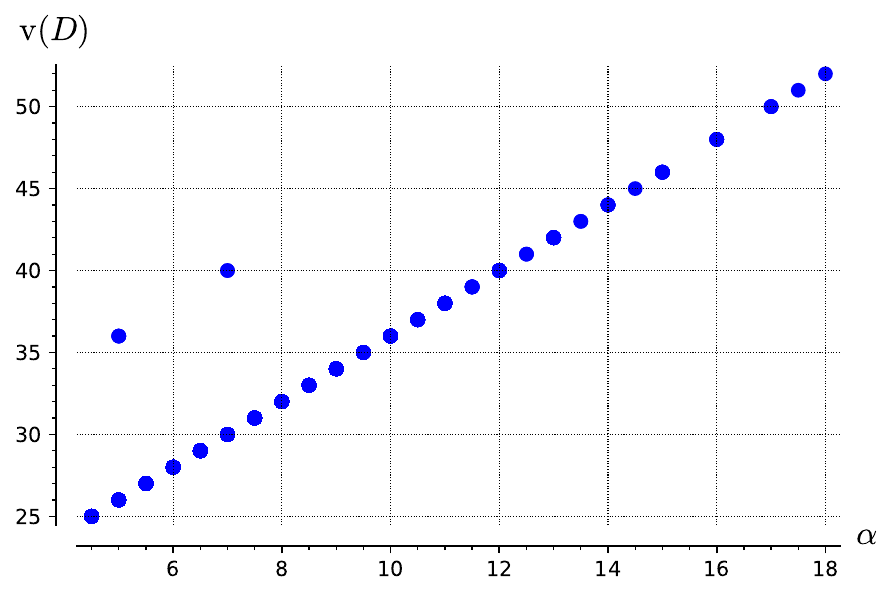}
	\caption{Valuation of the minimal discriminant for curves in case (C5) }
	\label{FigC5Discriminant}
\end{figure}

\begin{figure}[h!]
	\centering
	\includegraphics[width=1\textwidth]{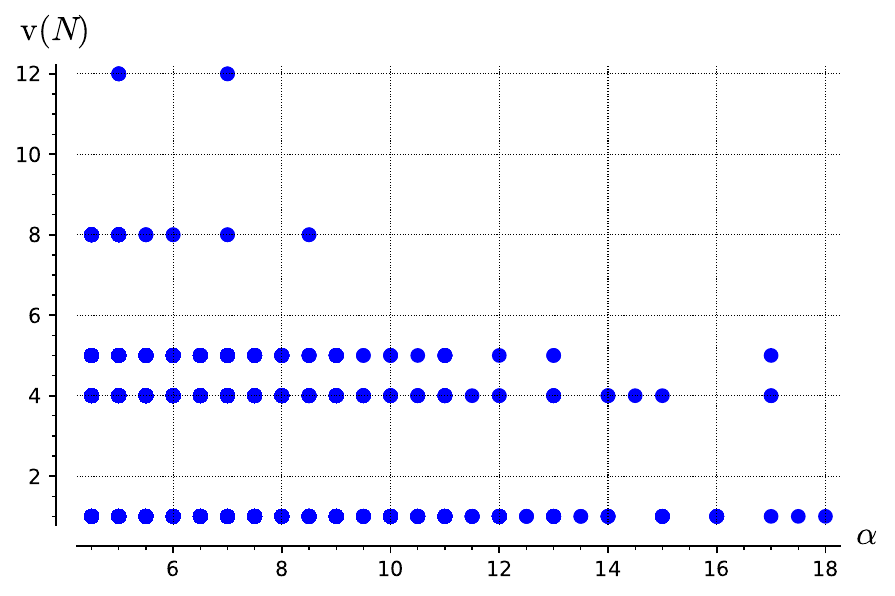}
	\caption{Conductor exponent values for curves in case (C5) }
	\label{FigC5Conductor}
\end{figure}

In Figures \ref{FigC5Discriminant} and \ref{FigC5Conductor}, we plot the valuation of the minimal discriminant and conductor exponent against $\alpha$ in case (C5). 
We chose (C5) because here we only have one parameter (here $\gamma=\delta=\epsilon=0$) and there are only two possibilities for the action of $\Gal(L/\BQ_2)$ on the dual graph of $C_0$. We see that Figure \ref{FigC5Discriminant} seems to suggest the formulas $16+2\alpha$ and $26+2\alpha$ for the minimal discriminant in the two cases. However, Figure \ref{FigC5Conductor} shows that our classification is not fine enough to yield a formula for the conductor exponent in this case. In particular, for the same value of $\alpha$ there are up to $5$ distinct conductor exponent values.

\begin{figure}[h!]
	\centering
	\includegraphics[width=1\textwidth]{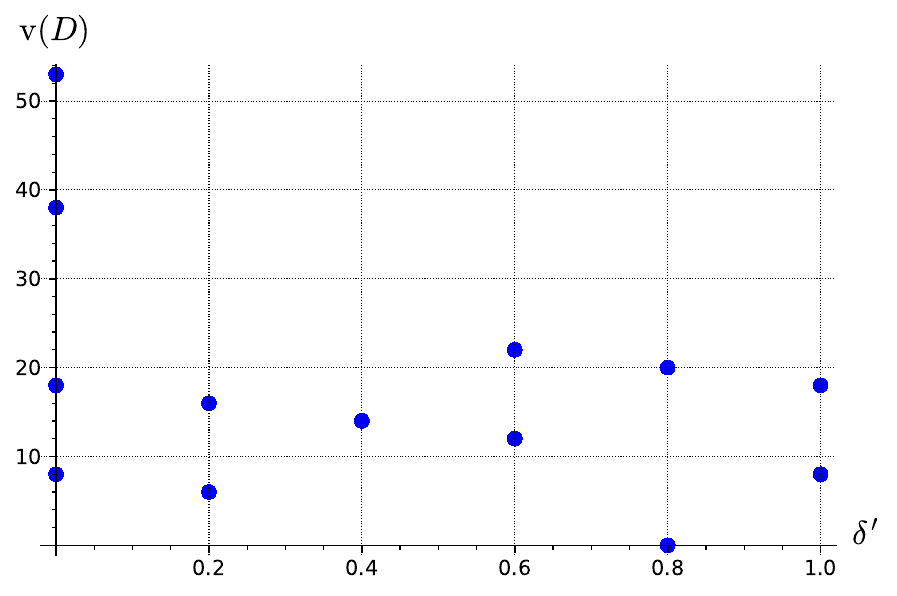}
	\caption{Valuation of the minimal discriminant for curves in case (A) with $\delta'\leq 1$ }
	\label{FigureADiscriminant}
\end{figure}

\begin{figure}[h!]
	\centering
	\includegraphics[width=1\textwidth]{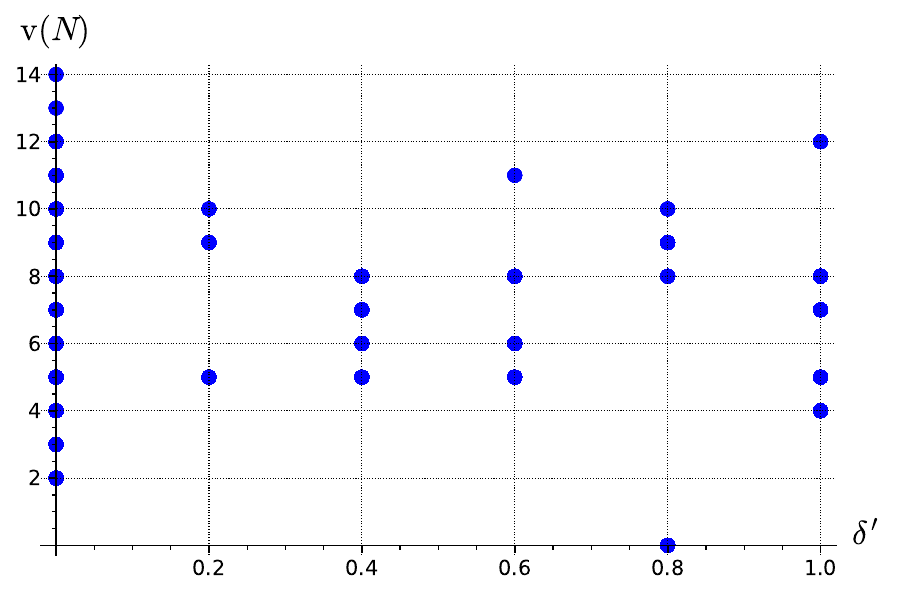}
	\caption{Conductor exponent values for curves in case (A) with $\delta'\leq 1$ }
	\label{FigureAConductor}
\end{figure}

\newpage
Similarly, in Figures \ref{FigureADiscriminant} and \ref{FigureAConductor}, we plot $\delta'$ against the valuation of the minimal discriminant and the conductor exponent in case (A)  for $\delta'\leq1$. We choose this bound on $\delta'$ since it implies $\delta=\delta'$. There are at most two possibilities for the action of $\Gal(L/\BQ_2)$ on the dual graph of $C_0$ and we again have only one varying parameter here as $\alpha=\beta=\gamma=\epsilon=0$. The data shows that  our classification is not fine enough to yield a formula for the conductor exponent in this case.

\newpage
\printbibliography
\end{document}